\def\nr#1{(\ref{#1})} \def\ekv#1#2{\begeq\label{#1}#2\endeq}
\def\eekv#1#2#3{\begin{eqnarray}\label{#1}#2 \\ #3
\nonumber\end{eqnarray}}

\def\iint{\int\hskip -2mm\int}
\def\iiint{\int\hskip -2mm\int\hskip
-2mm\int}  
scaled \magstep1  
   \def\3{\vert \hskip -1pt\vert\hskip -1pt\vert }
 \def\ably{arbitrarily} \def\an{analytic}
\def\asy{asymptotic} \def\bdd{bounded} \def\bdy{boundary}
\def\clas{classical analytic symbol} \def\coef{coefficient}
    \def\ctf{canonical
transformation} 
 
 \def\dop{differential operator}
 \def\ev{eigenvalue} \def\e{equation}
\def\fu{function} \def\fy{family} \def\F{Fourier} \def\fop{Fourier
integral operator} 
 \def\hol{holomorphic}
\def\indep{independent} \def\lhs{left hand side} \def\mfld{manifold}
\def\ml{microlocal} \def\neigh{neighborhood}
\def\nondeg{non-degenerate} \def\op{operator} 
\def\pb{problem}    
\def\plsh{plurisubharmonic} \def\pol{polynomial} \def\pro{proposition}
 \def\pop{pseudodifferential operator}
  \def\res{resonance} \def\rhs{right hand side}
\def\sa{selfadjoint} 
\def\sop{Schr{\"o}dinger operator} \def\st{strictly}

\def\sufly{sufficiently}
\def\tf{transformation}  
\def\tf{transform} \def\traj{trajectory} 
  \def\ufly{uniformly}
 \def\wrt{with respect to} 
\def\Re{{\mathrm Re\,}} \def\Im{{\mathrm Im\,}}
\documentclass[10pt]{article}
\usepackage{amssymb,latexsym}
\usepackage{graphics}
\usepackage{psfrag}

\title{Non-\sa{} perturbations of \sa{}
\op{}s in 2 dimensions IIIa. One branching point} \author{Michael
Hitrik\\Department of Mathematics \\ University of California
\\Los Angeles \\ CA 90095-1555, USA\\hitrik@math.ucla.edu \and
Johannes Sj{\"o}strand\\Centre de Math{\'e}matiques\\Laurent
Schwartz\\Ecole Polytechnique\\FR  91128
Palaiseau\\France\\johannes@math.polytechnique.fr} \date{}
\def\wrtext#1{\relax\ifmmode{\leavevmode\hbox{#1}}\else{#1}\fi}
\def\abs#1{\left|#1\right|} \def\begeq{\begin{equation}}
\def\endeq{\end{equation}} \def\Remark{\vskip 2mm \noindent {\em
Remark}}

\newcommand{\eps}{\epsilon}
\def\part#1{\frac{\partial}{\partial #1}}

\newcommand{\real}{\mbox{\bf R}}
\newcommand{\comp}{\mbox{\bf C}}
\newcommand{\z}{\mbox{\bf Z}}

\renewcommand{\exp}{\mbox{\rm exp\,}}

\newtheorem{dref}{Definition}[section]

\newtheorem{theo}[dref]{Theorem}
\newtheorem{prop}[dref]{Proposition}

\newenvironment{proof}{\medskip\noindent{{\em Proof:}}}{\hfill$\Box$
\medskip}

\begin{document}
\maketitle

\begin{abstract}
This is the third in a series of works devoted to spectral
asymptotics for non-\sa{}
perturbations of selfadjoint $h$-\pop{}s in dimension 2, having a
periodic classical flow. Assuming that the strength $\epsilon$
of the perturbation is in the range $h^2\ll \epsilon \ll h^{1/2}$
(and may sometimes reach even smaller values), we
 get an \asy{} description of the \ev{}s in rectangles
 $[-1/C,1/C]+i\epsilon [F_0-1/C,F_0+1/C]$, $C\gg 1$, when $\epsilon F_0$ is a saddle point
 value of the flow average of the leading perturbation.
\end{abstract}



\tableofcontents
\section{Introduction}\label{section0}
\setcounter{equation}{0}

This work is the third in a series devoted to non-\sa{}
perturbations of \sa{} semiclassical \pop{}s in two dimensions,
whose classical bicharacteristic flow is periodic on each energy
surface. The previous works in this series are~\cite{HiSj1,HiSj2},
and more recently, in collaboration with S. V\~u Ng\d{o}c, the
authors have begun a study of the case when the classical flow of
the unperturbed operator is no longer periodic but rather
possesses invariant Lagrangian tori with a Diophantine
property---see~\cite{HiSjVu} for the first work in this direction.

In this work, we continue with the perturbed periodic case. After switching on a
perturbation of size $\epsilon $, the spectrum will be confined to
a band of width ${\cal O}(\epsilon )$, and the more precise
distribution of \ev{}s is very much governed by the flow average
of the imaginary part of the leading symbol of the perturbation.
In the previous works, we studied the \ev{}s associated to
non-critical values of this flow average or to non-degenerate
maxima or minima in a suitable sense (after restriction to the
2-dimensional \mfld{} of trajectories in an energy surface). In
this paper we study the remaining generic case, namely that of
\ev{}s associated with a \nondeg{} saddle point.

\par We will work under the general assumptions of \cite{HiSj1},
\cite{HiSj2}, that we now recall. Let $M$ denote ${\bf R}^2$ or a
compact real-analytic \mfld{} of dimension 2. We shall let $M^{\bf
C}$ stand for a complexification of $M$, so that $M^{\bf
C}=\comp^4$ in the case when $M=\real^2$.

\par When $M={\bf R}^2$, let
\begeq
\label{Red.01}
P_\epsilon =P(x,hD_x,\epsilon ;h)
\endeq
be the Weyl quantization on ${\bf R}^2$ of a symbol $P(x,\xi
,\epsilon ;h)$ depending smoothly on $\epsilon \in{\rm
neigh\,}(0,{\bf R})$ with values in the space of \hol{} \fu{}s of
$(x,\xi )$ in a tubular \neigh{} of ${\bf R}^4$ in ${\bf C}^4$,
with \begeq\label{Red.02} \vert P(x,\xi ,\epsilon ;h)\vert \le
{\cal O}(1)m(\Re (x,\xi))
\endeq
there. Here $m$ is assumed to be an order \fu{} on ${\bf R}^4$, in
the sense that $m>0$ and for some $C_0$, $N_0>0$,
\begeq\label{Red.03} m(X)\le C_0\langle X-Y\rangle ^{N_0}m(Y),\
X,Y\in{\bf R}^4,\quad \langle X-Y\rangle :=(1+\vert X-Y\vert
^2)^{1\over 2}.
\endeq
We also assume that
\begeq\label{Red.04}
m\ge 1,\endeq
and
\begeq\label{Red.05}
P(x,\xi ,\epsilon ;h)\sim \sum_{j=0}^\infty  p_{j,\epsilon }(x,\xi)h^j,\ h\to 0, \endeq
in the space of such \fu{}s. We make the ellipticity assumption
\begeq\label{Red.06}
\vert p_{0,\epsilon }(x,\xi )\vert \ge {1\over C}m(\Re (x,\xi )),\
\vert (x,\xi )\vert \ge C, \endeq
for some $C>0$.

\par When $M$ is a compact \mfld{}, we let
\begeq\label{Red.07}
P_\epsilon =\sum_{\vert \alpha \vert \le m}a_{\alpha ,\epsilon
}(x;h)(hD_x)^\alpha ,\endeq
be a \dop{} on $M$, such that for every choice of local coordinates,
centered at some point of $M$, $a_{\alpha ,\epsilon }(x;h)$ is a smooth
\fu{} of $\epsilon $ with values in the space of \bdd{} \hol{} \fu{}s
in a complex \neigh{} of $x=0$. We further assume that
\begeq\label{Red.07.5}
a_{\alpha ,\epsilon }(x;h)\sim \sum_{j=0}^\infty  a_{\alpha ,\epsilon
,j}(x)h^j,\ h\to 0, \endeq
in the space of such \fu{}s. The semi-classical principal symbol in
this case is given by
\begeq\label{Red.08}
p_{0,\epsilon }(x,\xi )=\sum a_{\alpha ,\epsilon ,0}(x)\xi ^\alpha ,
\endeq
and we make the ellipticity assumption
\begeq\label{Red.09}
\vert p_{0,\epsilon}(x,\xi )\vert \ge {1\over C}\langle \xi \rangle ^m,\ (x,\xi
)\in T^*M,\,\vert \xi \vert \ge C,\endeq
for some large $C>0$. (Here we assume that $M$
has been equipped with some Riemannian metric, so that $\vert \xi \vert
$ and $\langle \xi \rangle =(1+\vert \xi \vert ^2)^{1/2}$ are
well-defined.)

\par Sometimes, we write $p_\epsilon $ for $p_{0,\epsilon }$ and
simply $p$ for $p_{0,0}$. Assume \begeq\label{Red.010} P_{\epsilon
=0} \hbox{ is formally \sa{}.}
\endeq
In the case when $M$ is compact, we let the under\-lying Hil\-bert
spa\-ce be $L^2(M,\mu (dx))$ for some po\-si\-ti\-ve real-ana\-ly\-tic den\-si\-ty
$\mu(dx)$ on $M$.

\par Under these assumptions, $P_\epsilon $ will have discrete spectrum
in some fixed \neigh{} of $0\in{\bf C}$, when $h>0,\epsilon \ge 0$
are \sufly{} small, and the spectrum in this region will be contained
in a band $\vert {\rm Im\,} z\vert \le {\cal O}(\epsilon )$.

\par Assume for simplicity that (with $p=p_{\epsilon =0}$)
\begeq\label{Red.011}
p^{-1}(0)\cap T^*M\hbox{ is connected.}\endeq
Let $H_p=p'_\xi \cdot {\partial \over \partial x}-p'_x\cdot {\partial
\over \partial \xi }$ be the Hamilton field of $p$. In this work, we
will always assume that for $E\in{\rm neigh\,}(0,{\bf R})$:
\begin{eqnarray}\label{H1}
&&\hbox{The }H_p\hbox{-flow is periodic on }p^{-1}(E)\cap T^*M\hbox{
with}\\ &&\hbox{period }T(E)>0 \hbox{ depending \an{}ally on }E.\nonumber
\end{eqnarray}
Let $q={1\over i}{({\partial \over \partial \epsilon })}_{\epsilon
=0}p_\epsilon $, so that
\begeq\label{Red.013}
p_\epsilon =p+i\epsilon q+{\cal O}(\epsilon ^2m),
\endeq
in the case when $M={\bf R}^2$, and $p_\epsilon =p+i\epsilon q+{\cal O}(\epsilon ^2\langle
\xi \rangle ^m)$ in the compact case. Let
\begeq\label{Red.014}
\langle q\rangle ={1\over T(E)}\int_{-T(E)/2}^{T(E)/2}q\circ \exp
tH_p dt\hbox{ on }p^{-1}(E)\cap T^*M.
\endeq
Notice that $p,\langle q\rangle $ are in involution,
$0=H_p\langle q\rangle =:\{ p,\langle q\rangle \} $. In \cite{HiSj1},
we saw how to reduce ourselves to the case when
\begeq
\label{Red.015}
p_\epsilon=p+i\epsilon \langle q\rangle +{\cal O}(\epsilon ^2),
\endeq
near $p^{-1}(0)\cap T^*M$. An easy consequence of this, also remarked
upon in~\cite{HiSj1}, is that the
spectrum of $P_\epsilon $ in $\{z\in {\bf C}; \vert \Re z\vert
<\delta \}$ is confined to $]-\delta ,\delta [+i\epsilon ]\langle
\Re q\rangle _{{\rm min},0}-o(1),\langle \Re q\rangle _{{\rm
max},0}+o(1)[$, when $\delta ,\epsilon ,h\to 0$, where $\langle
\Re q\rangle _{{\rm min},0}=\min_{p^{-1}(0)\cap T^*M}\langle \Re
q\rangle $ and similarly for $\langle \Re q\rangle _{{\rm max},0}$. We
shall mainly think about the case when $\langle q\rangle $ is
real-valued but will work under the more general assumption
that \begeq\label{H2} \Im \langle q\rangle \hbox{ is an \an{}
\fu{} of }p\hbox{ and }\Re \langle q\rangle ,\endeq
in a region of $T^*M$, where
$\vert p\vert \le 1/\vert {\cal O}(1)\vert $.

\par Let $\Lambda _{0,F_0}=\{ \rho \in T^*M; p(\rho )=0,\, \Re \langle
q\rangle (\rho )=F_0\}$. Assume
\begin{eqnarray}
\label{H3} &&T(0) \hbox{ is the minimal period for the $H_p$-flow}\\&&\hbox{at
every point of $\Lambda _{0,F_0}$ and $\Lambda _{0,F_0}$ is
connected.}\nonumber
\end{eqnarray}

\par The connectedness assumption is for convenience only and can easily
be removed. Then $\Sigma _0:=p^{-1}(0)/\exp ({\bf R}H_p)$ is a
symplectic 2-dimensional \mfld{} near the image
$\widetilde{\Lambda }_{0,F_0}$ of $\Lambda _{0,F_0}$. We consider
$\Re \langle q\rangle $ as an analytic \fu{} on ${\rm
neigh\,}(\widetilde{\Lambda }_{0,F_0},\Sigma _0)$. Assume
\begin{eqnarray}\label{H4}&&\hbox{This \fu{} has $F_0$ as critical value
and the
corresponding }\\ &&\hbox{critical point is unique, non-degenerate and of signature
0.}\nonumber
\end{eqnarray}

\par Then $\widetilde{\Lambda }_{0,F_0}$ is an $\infty $-shaped curve, and
$\langle q\rangle $ is an analytic \fu{} in a \neigh{} of that
curve (which is the level-curve of $\Re \langle q\rangle $
corresponding to $F_0$).

In the following, we may assume that $F_0=0$ for simplicity. In
Section \ref{SectionRed} we shall construct an $\epsilon
$-dependent \ctf{} $\kappa _\epsilon $ which is an $\epsilon
$-perturbation of a real \ctf{} $\kappa _0$, with \ekv{int.1} {
\kappa _\epsilon ,\kappa _0:\,{\rm neigh\,}(\{ \tau =0\}
,(T^*S^1)_{t,\tau }^{\bf C})\times {\rm neigh\,}(K_{0,0},{\bf
C}^2_{x,\xi })\to {\rm neigh\,}(\Lambda _{0,0},T^*M^{\bf C}), }
such that $p\circ \kappa _0=g(\tau )$, $\langle q\rangle \circ
\kappa _0=\langle q\rangle (\tau ,x,\xi )$ and \ekv{int.2} {
p_\epsilon \circ \kappa _\epsilon =g(\tau )+i\epsilon \langle
q\rangle (\tau ,x,\xi )+{\cal O}(\epsilon ^2), } where also the
${\cal O}(\epsilon ^2)$ is \indep{} of $t$. Here $K_{0,0}\subset
{\bf R}^2$ is an $\infty $-shaped curve with the self-crossing at
$(0,0)$ and  $(0,0)$ is the saddle point for the \fu{} $(x,\xi
)\mapsto \Re \langle q\rangle (0 ,x,\xi )$ with $\langle q\rangle
(0,x,\xi )=F_0\,\,(=0)$.

\par In the present work it seems quite essential to assume that
\ekv{int.3}{\hbox{The subprincipal symbol of }P_{\epsilon =0}\hbox{
vanishes.}}
(In \cite{HiSj1} this assumption was an optional one that permitted to get
improved results.)

\par After further reductions for the lower order symbols, described
in Section 2, we get a \ml{}
reduction of $P_\epsilon $ near $\Sigma _{0,0}$ to an \op{}
$\widehat{P}_\epsilon (hD_t,x,hD_x;h)$ with symbol
\begin{eqnarray}\label{int.4}
\widehat{P}_\epsilon (\tau ,x,\xi ;h)&=&g(\tau )+i\epsilon
\left(\langle q\rangle (\tau ,x,\xi )+{\cal O}(\epsilon )+{h^2\over
i\epsilon }p_2(\tau ,x,\xi )+{h\over
i}\widetilde{p}_1+...\right)\\& = &g(\tau )+i\epsilon Q(\tau ,x,\xi
,\epsilon ,{h^2\over \epsilon };h). \nonumber
\end{eqnarray}
The \op{} $\widehat{P}_\epsilon $ is only microlocally defined near
$\{ (t,\tau ,x,\xi )\in T^*S^1\times T^*{\bf R};\, \tau =0,\
(x,\xi )\in K_{0,0}\}$ but that allows us to define \asy{}ally its
\ev{}s in a rectangle $]-{1\over C},{1\over C}[+i\epsilon
]{-1\over C},{1\over C}[$ and they are of the form
\ekv{int.5} {
g\left(hk-{S_0\over 2\pi }-{k_0h\over 4}\right)+i\epsilon
w_{j,k},\ k\in {\bf Z}, }
where $w_{j,k}$ are the \ev{}s near $0$ of
$Q(hk-{S_0\over 2\pi }-{k_0\over 4}h,x,hD_x,\epsilon ,{h^2\over
\epsilon };h)$ in the \ml{} space $L^2_{\theta '}({\bf R})$
defined with Floquet conditions along the two loops of $K_{0,0}$
as in \cite{HiSj1}.  Here $\theta '=(\theta _1,\theta _2)\in{\bf
R}^2$ with $\theta _j=\frac{S_j}{2\pi }+{k_jh\over 4}$,
$k_j\in{\bf Z}$, and $(S_0,S_1,S_2)$ appear as action differences
when quantizing $\kappa _0$, while $k_0,k_1,k_2$ are Maslov
indices.

\par For $\tau \in{\rm neigh\,}(0,{\bf R})$, let $R(\tau )$ be the real
analytic curve formed by the values of $\langle q\rangle (\tau
,\cdot )$. That $R(\tau )$ is a curve follows from \nr{H2} and we
see that $R(\tau )$ is of the form \ekv{int.6} { \Im w=r(\tau ,\Re
w), } where $r$ is \an{} in a \neigh{} of $0$. Also, let $\rho
_c^0(\tau )\in{\rm neigh\,}((0,0),{\bf R}^2)$ be the critical
point of $\rho \mapsto \langle q\rangle (\tau ,\rho )$ (with $\rho
_c^0(0)=(0,0)$), and let $\rho _c(\tau )\in {\rm
neigh\,}((0,0),{\bf C}^2)$ be the critical point of the principal
symbol
\ekv{int.7} { (x,\xi )\mapsto Q^0(\tau ,x,\xi ,\epsilon
,\frac{h^2}{\epsilon })=\langle q\rangle +{\cal O}(\epsilon
)+{h^2\over i\epsilon }p_2(\tau ,x,\xi ) } of $Q(\tau
,x,\xi,\epsilon ,h^2/\epsilon ;h)$ appearing in \nr{int.4}.
Clearly, $\rho _c(\tau )=\rho _c^0(\tau )+{\cal O}(\epsilon
+{h^2\over \epsilon })$. Put $w_c(\tau )=\langle q\rangle (\tau
,\rho _c^0(\tau ))$ and introduce the exceptional boxes
\ekv{int.8} { {\cal B}(\tau )=\left\{ w ;\, \vert \Re
(w-w_c^0(\tau )\vert \le C_0(\epsilon +{h^2\over \epsilon }),\
\vert \Im w-r(\tau ,\Re w)\vert \le {{C_0}(\epsilon +{h^2\over
\epsilon })\over \vert \ln (\epsilon +{h^2\over \epsilon })\vert
}\right\} , } for some fixed \sufly{} large $C_0>0$.

\begin{theo}\label{ThInt.1}
We make the assumptions above and especially {\rm (\ref{Red.010})},
{\rm (\ref{H1})}, {\rm (\ref{H2})}, {\rm (\ref{H3})}, {\rm (\ref{H4})},
{\rm (\ref{int.3})}, and put $\tau _k=hk-{S_0\over 2\pi
}-{k_0h\over 4}$, $k\in \z$. Assume furthermore that $h^2\ll \eps \ll
h^{1/2}$. For $C>0$ \sufly{} large, the \ev{}s of
$P_\epsilon $ in $]-{1\over C},{1\over C}[+i\epsilon ]F_0-{1\over
C},F_0+{1\over C_0}[$ are of the form {\rm (\ref{int.5})}, where the following can
be said about the $w_{j,k}$:
\smallskip

\par The number of $w_{j,k}$ in ${\cal B}(\tau _k)$ is ${\cal O}({\epsilon \over
h}+{h\over \epsilon })\vert \ln (\epsilon +{h^2\over \epsilon })\vert .$
\smallskip

\par If $w_{j,k}\not\in {\cal B}(\tau _k)$, then $\vert\Re
(w_{j,k}-w_c^0(\tau _k))\vert >C_0(\epsilon
+{h^2\over \epsilon }) $, with $C_0$ as in {\rm (\ref{int.8})}.
\smallskip

\par There is a bijection $b_k$ between the set of these $w_{j,k}$
outside ${\cal B}(\tau _k)$ and the union of three sets of points away
from ${\cal B}(\tau _k)$: $E_{\rm ext}=E_e$, $E_{\rm leftint}=E_{li}$,
$E_{\rm rightint}=E_{ri}$ such that
$$
b_k(w)-w={\cal O}\left(e^{-\vert \Re(w-w_c^0(\tau _k))\vert\over Ch }{h\over
\vert \ln \vert \Re (w-w_c^0(\tau _k))\vert \vert }+h^\infty \right).
$$
\smallskip

\par Here $E_e$ is a subset of $\{ \Re (w-w_c^0(\tau _k))< -C_0 (\epsilon
+{h^2\over \epsilon })\}$ and $E_{li}$, $E_{ri}$ are subsets of
$\{ \Re (w-w_c^0(\tau _k))> C_0 (\epsilon
+{h^2\over \epsilon })\}$ (or vice versa but we only stick to the first
option for simplicity) that can be described by Bohr-Sommerfeld conditions
\ekv{int.9}
{
b_\Theta (w,\epsilon ,{h^2\over \epsilon },\tau _k;h)=2\pi (j+{1\over 2})h,\
j\in{\bf Z},\ \Theta =e,\, li,\, ri ,
}
where
\ekv{int.10}
{ b_\Theta (w,\epsilon ,{h^2\over \epsilon },\tau _k;h)\sim \sum_{\nu
=0}^\infty b_\Theta ^\nu (w,\epsilon ,{h^2\over \epsilon },\tau _k;h)h^\nu , } in
the space of \bdd{} \fu{}s of $w,\epsilon ,h^2/\epsilon ,\tau $, that are
smooth near $(0,0,0)$ in $(\epsilon ,h^2/\epsilon ,\tau )$ and \hol{} in
$w$ for
$$
\vert \Im (w-w_c(\tau ))\vert \le {1\over C}\vert \Re (w-w_c(\tau
))\vert ,\ \pm \Re (w-w_c(\tau _k))\ge C_0(\epsilon +{h^2\over \epsilon }),
$$
with a "$-$" when $\Theta =e$ and "$+$" for $\Theta =li,\, ri$.
\smallskip

\par Further,
\ekv{int.11}
{
b_\Theta ^1\hbox{ is \hol{} in a full \neigh{} of }w=w_c(\tau ),
}
\ekv{int.12}
{
b_\Theta ^\nu ={\cal O}(\vert w-w_c(\tau )\vert ^{1-\nu }),\ \nu \ge 2,
}
\begin{eqnarray}\label{int.13}
&&b_e^0(w,\epsilon ,{h^2\over \epsilon },\tau _k)-2\mu \ln (-\mu ),\ b_{li}-\mu \ln
\mu ,\ b_{ri}-\mu \ln
\mu \\&& \hbox{are \hol{} in a \neigh{} of }w=w_c(\tau ).
\nonumber
\end{eqnarray}
Here $\mu $ is a renormalized spectral parameter defined by $w =
K_{\epsilon ,h^2/\epsilon }(\tau _k,\mu ;h)$, with $K$ given in
Propositions {\rm (\ref{Prop1D2})}, {\rm (\ref{Prop1D1})}. Finally
$b_\Theta ^0$ can be described as actions along suitable cycles in
the complexified cotangent space, see Section {\rm
\ref{SectionEv}}.
\end{theo}

\par Inside the exceptional boxes the \ev{}s $w_{j,k}$ (for each fixed $k$)
continue to accumulate to roughly at most 5 curves where three of
the curves are the extensions of the curves carrying the $E_e$,
$E_{li}$, $E_{ri}$ (defined by replacing $2\pi hj $ \nr{int.9} by
a continuous real parameter) and one of the new curves, which
exists under certain conditions, can be related to barrier top
\res{}s in dimension 1. There are at most 2 and at least one point
(if we exclude degenerate cases) where three of the curves
terminate and form a "$Y$". Away from those points we may have
crossings of two of the curves (like for instance the ones
carrying $E_{li}$ and $E_{ri}$). Away from the "Y" points and with
some margin, the distribution of \ev{}s can be described by
Bohr-Sommerfeld rules as in the theorem, and near the "Y" points
as well as elsewhere, we can get quite detailed estimates for the
distribution of \ev{}s. Indeed, the \ev{}s can be identified with
zeros of quite explicitly given \hol{} \fu{}s which in most
regions appear as the sum of 4 exponential functions, and for such
functions it is possible to study the distribution of zeros quite
in detail. (See Davies ~\cite{Da} for inspiring results in this
direction and Hager~\cite{Ha} for quite elaborated results
obtained in parallel with the present work.) The appearance of
$Y$-shaped \ev{} distributions for non-\sa{} \op{}s in one
dimension seems to be quite well known and we refer to Shkalikov
\cite{Sh}, Redparth \cite{Re} and further references given there,
as well as to the recent works by L. Nedelec \cite{Ne} and E.
Servat and A. Tovbis~\cite{ServatTovbis}. The $Y$--shaped
eigenvalue distribution is also readily observed numerically---see
Figure 1 on this page, and also Figure 2 and Figure 3 on pages 8
and 9, respectively.

Unfortunately it turned out to be exceedingly difficult to give a
concise and precise description of what happens inside the
exceptional boxes in the form of a theorem in less than several
pages, so instead we refer the reader to the sections
\ref{SectionSk}--\ref{SectionEv} where this description can be
found.

\begin{figure}
\centering {\includegraphics{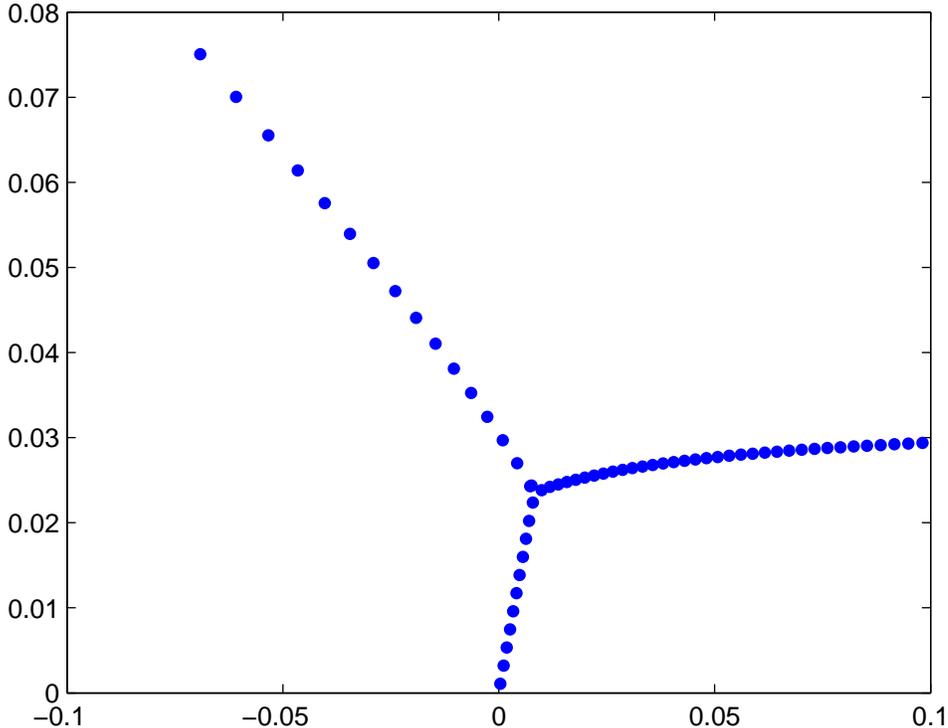}} \caption{Numerical
computation of the eigenvalues of the one-dimensional double well
Schr\"odinger operator perturbed by a complex even potential,
$P_{\eps}=(hD_x)^2-x^2+x^4+i\eps x^2$, in the case when $h=0.001$
and $\eps=0.8$. When computing the eigenvalues of $P_{\eps}$,
following~\cite{Tref}, we discretized the operator using the
Chebyshev spectral method.}
\end{figure}

\par In Section \ref{SectionBt} we apply our results to the study of
barrier top \res{}s for potentials of the form $-x^2+{\cal O}(x^4)$,
${\bf R}^2\ni x\to 0$. In the preceding section \ref{SectionPara}, we
make a remark that permits to improve the domain of validity in the
direction of small \res{}s. This gives an improvement also in the
applications to barrier tops in \cite{HiSj1,HiSj2} and allows us in the
present work to treat \res{}s $E$ with $h^{1-\delta }\ll \vert E\vert \ll
h^{1/3}$ for every fixed $\delta >0$, while a direct application of Theorem
\ref{ThInt.1} would only give the range $h^{2/3}\ll \vert E\vert \ll h^{1/3}$.
(In this special situation one can say that the lower bound $\epsilon \gg
h^2$, can be replaced by $\epsilon \gg h^{N_0}$ for every fixed $N_0>0$.)

\par When starting this project, we
underestimated the amount of ingredients needed, and in order to
keep the work within a reasonable size, we decided to exclude from
the paper the very interesting case when there are more than one
saddle point on the same connected component of $(\Re \langle
q\rangle )^{-1}(F_0)$ at real energy 0. The most important case
here is probably the one with 2 saddle points arising because of
an anti-symplectic involution (typically $(x,\xi )\mapsto (x,-\xi
)$ in the Schr{\"o}dinger case). We hope to take up at least the 2
saddle point case in a future work (having settled essentially all
heavy technicalities in the present work). We might then also
include the interesting case when \nr{H3} breaks down at isolated
points, leading to orbifolds. See Colin de Verdi{\`e}re--V\~u
Ng\d{o}c \cite{CoVu} in the \sa{} case.

\begin{figure}
\centering {\includegraphics{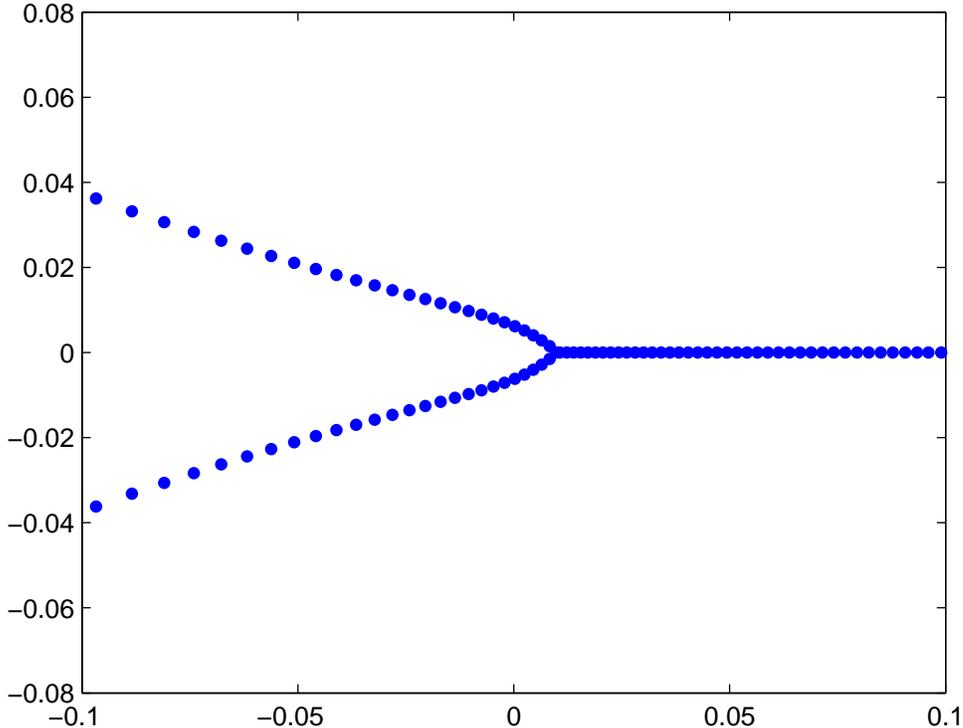}} \caption{Numerical
computation of the eigenvalues of the one-dimensional double well
Schr\"odinger operator perturbed by a complex odd potential,
$P_{\eps}=(hD_x)^2-x^2+x^4+i\eps x^3$, in the case when $h=0.001$
and $\eps=0.8$.}
\end{figure}

\par We also ran into a somewhat unexpected difficulty. Indeed, for the
one-dimensional \op{}s $Q$, we have a pseudospectral phenomenon
leading to an exponential growth of the resolvent norm in
important regions near the spectrum of these \op{}s. This makes it
very important to keep the errors in the reduction to the \op{}
$\widehat{P}_\epsilon $ in \nr{int.4} exponentially small, so that
the accumulated error in the global resolvent constructions
remains controlled. In \cite{HiSj1, HiSj2} we avoided that problem
by working in naturally adapted norms where the pseudospectral
problems disappeared, but that does not seem equally easy to do
here. This refined reduction is carried out in Section
\ref{SectionExp} using quasi-norms from the theory of analytic
\pop{}s originally due to Boutet de Monvel--Kr\'ee \cite{BoKr}, in
the simplified variant of \cite{Sj6}. The price to pay is the
apparent necessity to impose the condition \nr{int.3} and the
upper bound $\epsilon \ll h^{1/2}$, that should be compared to the
bound $\epsilon \ll h^\delta $ for every fixed $\delta >0$ in
\cite{HiSj1,HiSj2} or even $\epsilon \ll 1$ in \cite{Sj8,Sj2}.

\begin{figure}
\centering {\includegraphics{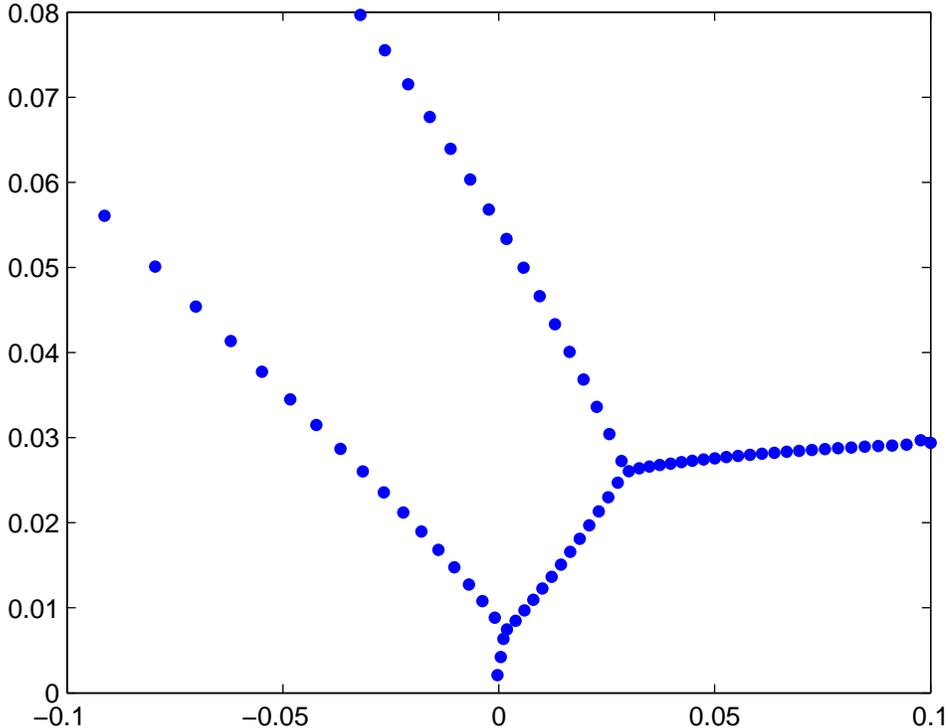}} \caption{Numerical
computation of the eigenvalues of the operator
$P_{\eps}=(hD_x)^2-x^2+x^4+i\eps\left(x^2+\delta x\right)$, in the
case when $h=0.001$, $\eps=0.8$, and $\delta=0.12$. Here $\delta$
is chosen so that the perturbation $x^2+\delta x$ is of the same
sign in both of the potential wells, although of different orders
of magnitude.}
\end{figure}

\medskip
\par \noindent
{\bf Acknowledgment.} We thank S. Fujiie and T. Ramond for
interesting discussions on the relationship with the complex
WKB-method. Those discussion may be important for future
generalizations. We also thank the MSRI, Berkeley, for the
possibility to advance on this project during the Spring of 2003.
The first author gratefully acknowledges the partial support
provided by the National Science Foundation under grant
DMS--0304970.

\section{Reduction to a one-dimensional
\pop{}}\label{SectionRed}
\setcounter{equation}{0}

\par Let $H_0\subset p^{-1}(0)$ be a hypersurface which is transversal to
the $H_p$-directions and such that $H_0$ can be identified with ${\rm
neigh\,}(\widetilde{\Lambda }_{0,0},\Sigma _0)$. We can then identify
$\widetilde{\Lambda }_{0,0}$ with a curve $\widetilde{H}_{0,0}$ in $H_0$.

\par Let $f=g^{-1}\circ p$, where $g$ is the unique increasing analytic \fu{}
with $g(0)=0$ such that the $H_f$-flow is $2\pi $ periodic with $2\pi $ as its minimal
period on $\Lambda _{0,0}$.

\par Let $\alpha :{\rm neigh\,}(K_{0,0},{\bf R}^2)\to {\rm
neigh\,}(\widetilde{H}_{0,0},H_0)$ be a real-\an{} \ctf{}, where
$K_{0,0}$ is an $\infty $-shaped curve, as in the introduction.
The existence of such a map with a suitable $K_{0,0}$ follows from
Theorem 5 in~\cite{DacMoser} (see also~\cite{Siburg}), according
to which if $\Omega_1$ and $\Omega_2$ are two compact
real-analytic symplectic manifolds of dimension 2, possibly with
boundary, which have the same area and such that there exists an
orientation preserving analytic diffeomorphism between them, then
$\Omega_1$ can be mapped onto $\Omega_2$ by an analytic canonical
transformation.

We extend $\alpha $ to a \ctf{}
$$\kappa :\,{\rm neigh\,}(\{ \tau =0\} ,
T^*S^1)_{t,\tau })\times {\rm neigh\,}(K_{0,0},{\bf R}^2_{x,\xi })
\to{\rm neigh\,}(\Lambda _{0,0},T^*M),$$ with $f\circ \kappa =\tau
$ in the following way: Extend $H_0$ to an analytic hypersurface
$H$ in the full phase space which intersects $p^{-1}(0)$
transversally along $H_0$. Let $\widetilde{t}$ be the
grad-periodic \fu{} defined near $\Lambda _{0,0}$ which solves
$$H_f\widetilde{t}=1, \hbox{ with }\widetilde{t}=0\hbox{ on }H.$$
Because of the $2\pi $-periodicity of the $H_f$-flow, we see that
$\widetilde{t}$ is well-defined up to an integer multiple of $2\pi
$. Using that $H_f$ and $H_{\widetilde{t}}$ commute, we notice
that if $\rho $ is a point close to $\Lambda _{0,0}$, then we can
write
$$\rho =\exp (tH_f+\tau H_{\widetilde{t}})(\rho _0),\ \rho _0\in H_0,$$
where $\tau \in {\bf R}$ is small and $t\in{\bf R}$ is well-defined modulo
a multiple of $2\pi $.

\par Put $\kappa (t,\tau ;x,\xi )=\exp (tH_f+\tau
H_{\widetilde{t}})(\alpha (x,\xi ))$. Then $\kappa $ has the desired
properties.

As in Section 2 of \cite{HiSj2}, we introduce the triple
$S=(S_0,S_1,S_3)\in {\bf R}^3$ of action differences, with $S_0$
corresponding to a closed $H_p$-orbit $\subset p^{-1}(0)$, and
$S_1$, $S_2$ corresponding to the left and right closed orbits of
the $\infty $-shaped set $\widetilde{H}_{0,0}$. Let $\theta
=(\theta _0,\theta _1,\theta _2)$, $\theta _j=S_j/(2\pi h)+k_j/4$,
where $k_j\in {\bf Z}$ is a suitable Maslov index. Let $L^2_\theta
(S^1\times {\bf R})$ be the space of \ml{}ly defined \fu{}s
$u(t,x)$ in ${\rm neigh\,}(\{ \tau =0\} ,(T^*S^1)_{t,\tau })\times
{\rm neigh\,}(K_{0,0},{\bf R}^2)$ that are multivalued but $\theta
$-Floquet periodic as in \cite{HiSj1}, Section 6, or as in
\cite{MeSj}. Let $U:L^2_\theta (S^1\times {\bf R})\to L^2(M)$ be a
\ml{}ly defined unitary \fop{} as in the cited works.

\par Repeating the argument in the beginning of Section 3 of~\cite{HiSj1}, we may and will assume
from now on that the leading perturbation $q$ in (\ref{Red.013})
has already been averaged along the $H_p$--flow so that the
leading symbol of $P_{\eps}$ becomes
\ekv{Red.2}{p_\epsilon
=p+i\epsilon \langle q\rangle +{\cal O}(\epsilon ^2).}
The operator
$$
\widetilde{P}_\epsilon :=U^{-1}P_\epsilon U,
$$
has the principal symbol
\begeq
\label{Red.3}
\widetilde{p}_\epsilon =g(\tau )+i\epsilon \langle q\rangle (\tau
,x,\xi )+{\cal O}(\epsilon ^2).
\endeq

\par At this stage, we get a complete analogue to the situation in Section 3 in
\cite{HiSj1}. Proposition 3.2 there extends, and we get a
reduction of $\widetilde{P}_{\epsilon}$ to an operator
$\widehat{P}_\epsilon =\widehat{P}_\epsilon (hD_t,x,hD_x;h)$, also
acting on $L^2_\theta $, with a complete symbol \indep{} of $t$,
and the principal symbol still of the form
\ekv{Red.4}{\widehat{p}_\epsilon =g(\tau )+i\epsilon \langle
q\rangle (\tau ,x,\xi )+{\cal O}(\epsilon ^2),} now completely
\indep{} of $t$.

\par As in Section 5 of \cite{HiSj1}, we use \F{} series expansions in the
$t$-variable and get a reduction of (\ref{Red.4}) to the family of \op{}s:
\ekv{Red.5}
{\widehat{P}_\epsilon \left(h(k-{k_0\over 4})-{S_0\over 2\pi },x,hD_x;h\right),
\ k\in{\bf Z},}
where $k_0\in{\bf Z}$ is a fixed Maslov index and it is understood that we
only consider such values of $k$ for which the first argument ($\tau $)
in $\widehat{P}_\epsilon $ is small.
\smallskip

\par\noindent a) In the general case, without any assumption on the
subprincipal symbol of $P_{\epsilon =0}$, we write the full symbol of
$\widehat{P}_\epsilon $ as
\ekv{Red.6}
{
\widehat{P}_\epsilon (\tau ,x,\xi ;h)=g(\tau )+\epsilon \left[i\langle
q\rangle (\tau ,x,\xi )+{\cal O}(\epsilon )+{h\over \epsilon }p_1(\tau
,x,\xi )+h{h\over \epsilon }p_2(\tau ,x,\xi )+...\right],
}
and consider $h/\epsilon $ as an additional small parameter.
\smallskip

\par\noindent b) When the subprincipal symbol of $P_{\epsilon =0}$
vanishes, we have the same fact for $\widehat{P}_\epsilon $ by the
improved Egorov property of $U$---see Section 2 of~\cite{HiSj1}.
Thus $(p_1)_{\epsilon =0}=0$ and we can write
$${h\over \epsilon }p_1(\tau ,x,\xi ,\epsilon )=h\widetilde{p}_1(\tau
,x,\xi ,\epsilon ).$$
Instead of (\ref{Red.6}), we get
\ekv{Red.8}
{\widehat{P}_\epsilon (\tau ,x,\xi ;h)=g(\tau )+\epsilon \left[i\langle
q\rangle (\tau ,x,\xi )+{\cal O}(\epsilon )+{h^2\over \epsilon }p_2(\tau
,x,\xi )+h\widetilde{p}_1+h^2\widetilde{p}_2+...\right],}
depending on the small parameters $\epsilon ,h^2/\epsilon $.

\par Summing up the discussion so far, we have
\begin{prop}\label{PropRed1}
Let $P_\epsilon $ be as above satisfying the assumptions {\rm
(\ref{Red.02})}, {\rm (\ref{Red.03})}, {\rm (\ref{Red.04})}, {\rm
(\ref{Red.05})}, {\rm (\ref{Red.06})}, {\rm (\ref{Red.09})}, {\rm
(\ref{Red.010})}, {\rm (\ref{Red.011})}, {\rm (\ref{H1})}, {\rm
(\ref{H2})}, {\rm (\ref{H3})}, {\rm (\ref{H4})}. Also assume that $0\le
\epsilon \le h^\delta $ for some fixed $\delta
>0$. Then there exist $G_0(x,\xi )$ \hol{} in some fixed \neigh{} of
$p^{-1}(0)$, an elliptic \fop{} $U$ of order 0, with the
associated \ctf{} $\kappa $ as above, and an $h$-\pop{}
$A=A(t,hD_t,x,hD_x;h)$ of order $0$ with principal symbol ${\cal
O}(\epsilon ^2)$, such that the \op{} \ekv{Red.9} {
\widehat{P}_\epsilon =e^{{i\over h}A}U^{-1}e^{-{\epsilon \over
h}G_0^w} P_\epsilon e^{{\epsilon \over h}G_0^w}Ue^{-{i\over h}A}
={\rm Ad}_{e^{{i\over h}A}U^{-1}e^{-{\epsilon \over
h}G_0}}P_\epsilon } has a symbol $\widehat{P}_\epsilon (\tau
,x,\xi ,\epsilon ;h)$ \indep{} of $t$, modulo ${\cal O}(h^\infty
)$. Here $G_0^w=G_0^w(x,hD_x)$.

\par In the general case, we have {\rm (\ref{Red.6})}, provided that $h\ll
\epsilon \le h^\delta $, and when the subprincipal symbol of
$P_{\epsilon =0}$ vanishes, we have {\rm (\ref{Red.8})}
provided
that $h^2\ll \epsilon \le h^\delta $.
\end{prop}

\par In this proposition all symbols and phase \fu{}s are \hol{} in fixed
$h,\epsilon $-\indep{} domains. The weight $G_0$ in (\ref{Red.9})
is used to get a first reduction of the principal symbol to the
form (\ref{Red.2})---see also (\ref{Sk.52}) in section 8.

\section{Exponential decoupling}\label{SectionExp}
\setcounter{equation}{0}

Since some solution operators to the localized 1-dimensional problems
later on will have some exponential growth, the decoupling result
of Section \ref{SectionRed} has
 to be sharpened in the sense that we get
some exponential smallness control over the remainders.

\par First we need to recall some notions about \clas{}s and their
associated quasi-norms. That was introduced in the pioneering work
by L. Boutet de Monvel and P. Kr\'ee~\cite{BoKr}, but here we
shall use the simplified quasi-norms of \cite{Sj6}. If $\Omega
\subset {\bf C}_{x,\xi }^{2n}$ is open, a \clas{} of order 0 is
given by the formal \asy{} expansion, \ekv{Exp.01} { a(x,\xi
;h)\sim \sum_{k=0} ^\infty h^ka_k(x,\xi ),} where $a_k$ are \hol{}
in $\Omega $ and satisfy the growth condition, \ekv{Exp.02} {
\forall K\subset\subset \Omega ,\, \exists C=C_K>0;\ \vert
a_k(x,\xi )\vert \le C^{k+1}k^k,\ (x,\xi )\in K. }

\par To such an $a$, we associate the formal \dop{} of infinite order,
\ekv{Exp.03}
{A(x,\xi ,D_x;h)=a(x,\xi +hD_x;h)\sim \sum_{k=0}^\infty h^kA_k(x,\xi ,D_x),}
where
$$A_k=\sum_{\ell+\vert \alpha \vert =k}{1\over \alpha !}\partial _\xi
^\alpha a_\ell (x,\xi )D_x^\alpha .
$$
Let $\Omega _t\subset\subset \Omega $, $t_0\le t\le t_1$ be an
increasing \fy{} of open subsets with $t_0<t_1$, such that
$${\rm dist}_\infty (\Omega _s,\partial \Omega _t)\ge t-s,\hbox{ for
}t_0\le s < t\le t_1,$$
with the distance associated to the $\ell^\infty$-norm.
Let $f_j(A_j)\ge 0$ be the smallest constant such that
$$\Vert A_j\Vert _{s,t}\le f_j(A_j)\big( {j\over t-s}\big)^j,$$
where $\Vert \cdot \Vert _{s,t}$ is the \op{} norm from the space of
\bdd{} \hol{} \fu{}s on $\Omega _t$ to the same space on $\Omega _s$. Then
$$\3 a\3_\rho :=\sum_0^\infty \rho ^j f_j(A_j)$$
is finite for $\rho >0$ small enough, and conversely, the finiteness of $\3
a\3_\rho $ for some fixed $\rho >0$ implies that $a=A(1)$ is a \clas{} on
$\Omega _{t_1}$. If $a,b$ are \clas{}s on $\Omega $ and we let $a(x,hD;h)$
and $b(x,hD;h)$ denote the associated $h$-\pop{}s for the classical
quantization, then the composition of these two \op{}s has the symbol
$$
a\# b\sim \sum_{k=0}^\infty h^k\sum_{\vert \alpha \vert =k}{1\over \alpha
!}\partial _\xi ^\alpha a(x,\xi ;h)D_x^\alpha b(x,\xi ;h).
$$
Using that the \dop{}s $A$, $B$ compose correspondingly, it was shown
very simply in \cite{Sj6} that
\ekv{Exp.04}
{
\3 a\# b\3_\rho \le \3 a\3_\rho \3 b\3_\rho ,
}
implying that the composed symbol is also a \clas{}.

\par If we prefer to work with the Weyl quantization, then the same result
and proof as in \cite{Sj6} remain valid, provided that we modify
the choice of the associated infinite order \dop{} to
\begin{equation}
\label{Exp.05} A(x,\xi,D_{x,\xi};h)u=\sum_{k=0}^\infty
\frac{h^k}{k!} \left(\left(\frac{i}{2} \sigma (D_{x,\xi
};D_{y,\eta})\right)^k a(x,\xi ;h)u(y,\eta)\right)\bigg
\vert_{y=x,\eta =\xi},
\end{equation}
so that
$$
A(x,\xi ,D_{x,\xi };h)=a(x-{h\over 2}D_\xi ,\xi +{h\over 2}D_x;h).
$$

\par
In the following, we will consider symbols depending on additional
parameters; \ekv{Exp.06} { a(x,\xi ,\epsilon ,h;h)=\sum
h^ka_k(x,\xi ,\epsilon  , h), } including $h$ (which is then
viewed as an \indep{} parameter) and consequently the admissible
values of $\rho $ will depend on these parameters. When defining
$\3 a\3_\rho $ in the case when $h$ is among the additional
parameters, we have in mind some specific representation of the
form (\ref{Exp.06}). To have $\rho $ tending to 0 as some power of
$h$---and that is what we will encounter, means roughly that we
consider Gevrey symbols.

\begin{prop}
\label{PropExp1}
Let $\ell (x,\xi )$ be affine and real and let $a(x,\xi ;h)$ be  an
analytic symbol of order 0. Then
$$\3 [\ell (x,hD),a(x,hD;h)]\3_\rho \le 2\rho \Vert \nabla \ell \Vert _\infty
\3 a\Vert _\rho .$$
\end{prop}

\par Here and in the following we shall consider $[\ell ,a]$ as an
$h$-\pop{} of order 0.

\begin{proof}
The symbol of $[ \ell(x,hD),a(x,hD;h)]$ is equal to
$${h\over i}\{ \ell ,a\}={h\over i}\nu (\partial _{x,\xi })a,$$
where $\nu =H_\ell$ is the Hamilton field of $\ell$. To $a$ we associate
the infinite order \dop{} $A(x,\xi ,D_x;h)$ as in (\ref{Exp.03}). Similarly, we have
$$[\ell (x,hD),a(x,hD;h)]\leftrightarrow \sum_{j=1}^\infty  h^jB_j(x,\xi
,D_x) ={h\over i}\sum_{k=0}^\infty  h^kC_k,$$
with $C_k=\nu (\partial _{x,\xi })(A_k)$ in the sense that $\nu $ acts as
a \dop{} on the \coef{}s of $A_k$. Thus,
\ekv{Exp.1}
{
B_j={1\over i}C_{j-1}={1\over i}\nu (\partial _{x,\xi })(A_{j-1}).
}
We can also express this as
\ekv{Exp.1.5}{B_j={1\over i}\Big( {\partial \over \partial r}\Big)_{r=0} \tau
_{-r\nu }\circ A_{j-1}\circ \tau _{r\nu },}
where $\tau _{r\nu }$ denotes translation in the complex domain by the
vector $r\nu $ and $r$ may be complex.

\par Assume for simplicity that $\vert \nu \vert =\Vert \nu \Vert _\infty \le 1$.
Then for $(x,\xi )\in \Omega _s$, $t_0\le s+\vert r\vert <\widetilde{t}\le
t_1-\vert r\vert $ and for $u$ \hol{} in a
larger domain, we get
\begin{eqnarray*}
\vert (\tau _{-r\nu }\circ A_{j-1}\circ \tau _{r\nu }u)(x,\xi )\vert
&=&\\
\vert (A_{j-1}\tau _{r\nu }u)((x,\xi )+r\nu )\vert &\le&
{f_{j-1}(A_{j-1})(j-1)^{j-1}\over (\widetilde{t}-(s+\vert r\vert
))^{j-1}}\sup_{\Omega _{\widetilde{t}}} \vert \tau _{r\nu }u\vert \\
&\le& f_{j-1}(A_{j-1}){(j-1)^{j-1}\over (\widetilde{t}-(s+\vert r\vert
))^{j-1}}\sup_{\Omega _{\widetilde{t}+\vert r\vert }}\vert u\vert .
\end{eqnarray*}
For $t_0\le s<t\le t_1$,
choose $\widetilde{t}=t-\vert r\vert $, $2\vert r\vert <t-s$:
$$\vert \tau _{-r\nu }\circ A_{j-1}\circ \tau _{r\nu }u(x,\xi )\vert \le
f_{j-1}(A_{j-1}){(j-1)^{j-1}\over (t-s-2\vert r\vert )^{j-1}}\sup_{\Omega
_t}\vert u\vert  .$$
In other words,
$$
\Vert \tau _{r\nu }\circ A_{j-1}\circ \tau _{r\nu }\Vert _{s,t}\le
{f_{j-1}(A_{j-1})(j-1)^{j-1}\over (t-s-2\vert r\vert )^{j-1}},
$$
and from the Cauchy inequality and (\ref{Exp.1.5}), we get
$$\Vert B_j\Vert _{s,t}\le {f_{j-1}(A_{j-1})(j-1)^{j-1}\over \delta
(t-s-2\delta )^{j-1}}={2f_{j-1}(A_{j-1})\over (t-s)^j}{(j-1)^{j-1}\over
{2\delta \over (t-s)}(1-{2\delta \over t-s})^{j-1}},\ 0<2\delta <t-s.$$
Here we choose $\delta $
so that $\theta :={2\delta \over t-s}$ minimizes ${1\over \theta (1-\theta
)^{j-1}}$. We find
$$\theta ={1\over j},\quad {1\over {1\over j}(1-{1\over j})^{j-1}}={j^j\over (j-1)^{j-1}}.$$
Hence,
$$\Vert B_j\Vert _{s,t}\le {2f_{j-1}(A_{j-1})\over (t-s)^j}j^j,$$
so
$$f_j(B_j)\le 2f_{j-1}(A_{j-1}),$$
and
$$\3 B\3 _\rho \le 2\sum_{j=1}^\infty  f_{j-1}(A_{j-1})\rho ^j =2\rho\3 A\3
_\rho  .$$
\end{proof}

\par In the following, we allow a finite number of families $(\Omega
_t^\nu )_{t_1^\nu \le t\le t_2^\nu }$, $\Omega _t^\nu
\subset\subset \Omega $, $\nu =1,...,N$ as above. If $\3 a\3^\nu
_\rho $ denotes the $\rho $-quasi-norm, defined with the help of
the $\nu $-th \fy{}, we define \ekv{Exp.1.6} {\3 a\3_\rho
=\sum_\nu \3 a\3_\rho ^\nu . } Here it is understood that if $a$
is parameter dependent with $h$ among the parameters, then we use
the \underline{same} representation (\ref{Exp.06}) when defining
each of the quasi-norms $\3 a\3_\rho ^\nu $. Notice that we still
have (\ref{Exp.04}).

\begin{prop}\label{PropExp2}
Let $g(x,\xi ;h)$ be a \clas{} of order 0, defined in a finite
union $D=\cup_{\nu =1}^N D_\nu $ of polydiscs with
$\overline{\Omega }\subset\subset D$, and let $g(x,hD;h)$ be the
corresponding $h$-pseudo\-diffe\-ren\-tial operator. Let $\3\cdot \3_\rho $ be a quasi-norm of
the form {\rm (\ref{Exp.1.6})} with the corresponding family
$\Omega _t^\nu \subset\subset D_\nu $. Then for $\rho $ small
enough, we have \ekv{Exp.2} { \3 [g(x,hD;h),a(x,hD;h)]\3 _\rho \le
C(g) \rho \3 a\3 _\rho . }
\end{prop}
\begin{proof}
Write $g=g_0(x,\xi )+hg_1(x,\xi ;h)$, where $g_1$ is a \clas{} of order 0. We
notice that
$$\3 hg_1\3 _\rho \le C\rho \3 g_1\3_\rho ,$$
(where $hg_1$ is viewed as a symbol of order 0), so on the \op{} level, we
have
\ekv{Exp.3}
{
\3 [hg_1,a]\3_\rho \le 2\3 hg_1\3_\rho \3 a\3_\rho  \le 2C\rho \3g_1\3_\rho
\3a\3_\rho .
}
Hence it only remains to treat the contribution to the commutator from
$g_0(x,hD)$. We may assume we work in a polydisc centered at
(0,0) with the radii $r_1,r_2,..,r_n,s_1,..,s_n$. Then,
$$
g_0(x,\xi )=\sum_{\alpha ,\beta \in{\bf N}^n}g_0^{\alpha ,\beta }x^\alpha
\xi ^\beta ,
$$
where
$$\sum \vert g_0^{\alpha ,\beta }\vert r^\alpha s^\beta <\infty .$$
Now choose the classical quantization for simplicity.
On the \op{} level,
\ekv{Exp.4}
{
g_0(x,hD)=\sum g_0^{\alpha ,\beta }x^\alpha (hD)^\beta ,
}
where the sum converges in the space of \an{} symbols, since
$\3 x_j\3_\rho \le r_j,\quad \3 hD_{x_j}\3_\rho \le s_j+\rho
=:\widetilde{s}_j$,
and we can allow some shrinking in $r_j,s_j$ and choose $\rho >0$ small
enough.

\par Using that $\3 [x_j,a]\3_\rho , \3 [hD_{x_j},a]\3_\rho \,\le 2\rho
\3a\3_\rho $, in view of Proposition 3.1, we see that
$$
\3 [x^\alpha (hD_x)^\beta ,a]\3_\rho \le (\alpha _1r^{\alpha
-e_1}\widetilde{s}^\beta +\alpha _2r^{\alpha -e_2}\widetilde{s}^\beta
+..+\alpha _nr^{\alpha -e_n}\widetilde{s}^\beta +\beta _1r^\alpha
\widetilde{s}^{\beta -e_1}+..+\beta _nr^\alpha \widetilde{s}^{\beta
-e_n})2\rho \3 a\3_\rho ,
$$
which can be written more briefly as
$$
\3 [x^\alpha (hD_x)^\beta ,a]\3_\rho \le (\partial _{r_1}+..+\partial
_{r_n}+\partial _{\widetilde{s}_1}+..+\partial
_{\widetilde{s}_n})(r^\alpha \widetilde{s}^\beta )2\rho \3 a\3_\rho .
$$
Hence
$$
\3 [g_0(x,hD),a\3_\rho \le \big[ (\partial _{r_1}+..+\partial
_{r_n}+\partial _{\widetilde{s}_1}+..+\partial
_{\widetilde{s}_n})(\sum_{\alpha ,\beta }\vert g_0^{\alpha ,\beta }\vert
r^\alpha \widetilde{s}^\beta )\big]2\rho \3 a\3_\rho .
$$
\end{proof}

\par\noindent \it Remark. \rm The \pro{} remains valid
for the Weyl quantization. Indeed, only (\ref{Exp.4}) has to be modified, by
adding a term $h\widetilde{g}_1(x,hD;h)$ to the \rhs{}, where $g_1$ is an
\an{} symbol of order 0.\medskip

%\par We need to go a step further in a special case:
%\begin{prop}\label{PropExp3}
%Let $g(\tau )$ be an \an{} \fu{} of one variable. Then if
%$a(t,x;hD_t,hD_x;h)$ is an \an{} $h$-\pop{}, we have
%\ekv{Exp.5}
%{
%[g(hD_t),a]={h\over i}\{ g,a\} (t,x,hD_{t,x})+r, \hbox{where }\3 r\3_\rho
%\le C(g) \rho ^2\3 a\3_\rho ,
%}
%for $\rho $ small enough.
%\end{prop}
%\begin{proof}
%Write $g(\tau )=\sum g_j\tau ^j$ as the sum of a convergent power series.
%Then $g(hD_t)=\sum g_j(hD_t)^j$ with convergence in $\rho $ norm for
%\sufly{} small $\rho $. Now,
%\begin{eqnarray}\label{Exp.6}
%[g(hD_t),a]&=& \sum_{j=1}^\infty  g_j\sum_{k_1+k_2=j-1\atop k_1,k_2\ge
%0}(hD_t)^{k_1}[hD_t,a](hD_t)^{k_2}\\
%&=&\sum_{j=1}^\infty  g_jj[hD_t,a](hD_t)^{j-1}+\sum_{j=2}^\infty g_j\sum_{k_1+k_2=j-1\atop k_1,k_2\ge
%0} [(hD_t)^{k_1},[hD_t,a]](hD_t)^{k_2}.
%\nonumber
%\end{eqnarray}
%The first term in the last expression is equal to $(h/i)\{ g,a\}
%(t,x,hD_{t,x};h)$.
%
%\par Using \Pro{} \ref{PropExp1} twice, we have
%$$\3 [(hD_t)^{k_1}[hD_t,a]](hD_t)^{k_2}\3_\rho \le k_1\3 hD_t\3_\rho
%^{k_1+k_2-1}\3 [hD_t,[hD_t,a]]\3_\rho \le 4k_1\3 hD_t\3_\rho ^{k_1+k_2-1}\rho
%^2\3a\3_\rho ,$$
%so denoting the last term in (\ref{Exp.6}) by $r$, we get
%\begin{eqnarray}\label{Exp.7}
%\3 r\3_\rho &\le &(\sum_{j=2}^\infty  g_j\sum_{k_1+k_2=j-1\atop k_1\ge 1,
%k_2\ge 0}4k_1\3 hD_t\3 ^{j-2}_\rho )\rho ^2\3 a\3_\rho \\
%&\le& (\sum_{j=2}^\infty  g_j4(j-1)^2\3 hD_t\3_\rho ^{j-1})\rho ^2\3
%a\3_\rho .\nonumber
%\end{eqnarray}
%\end{proof}

\par Returning to the considerations of Section 2, let us consider the \an{} symbol,
\ekv{Exp.8} { P=g(\tau )+\epsilon q(\tau ,x,\xi )+hr(t,\tau ,x,\xi
,\epsilon ;h)=g(\tau )+\epsilon \widetilde{q}, } defined in ${\rm
neigh\,}(t\in S^1,\tau =0;((S^1+i{\bf R})\times {\bf C})\times
\Omega)$, where $\Omega \subset {\bf C}^2_{x,\xi }$ is open. Here
we assume either that $r$ is a \clas{} of order 0 or simply that
$r$ has an \asy{} expansion in integral powers of $h$ in the space
of \hol{} \fu{}s. We have already seen in Section 2 that after a
finite number of conjugations, we may assume that $r$ is \indep{}
of $t$ modulo ${\cal O}(h^N)$.  We also assume that $g'\ne 0$.

\par In the general case, we have
\ekv{Exp.8.1}
{\widetilde{q}={\cal O}(1+{h\over \epsilon }),}
and when
the subprincipal symbol of $P_{\epsilon =0}$ vanishes, we have $hr={\cal
O}(h^2+\epsilon h+\epsilon ^2)$,
\ekv{Exp.8.2}
{\widetilde{q}={\cal O}(1+h^2/\epsilon +h+\epsilon ).} In the two cases, we shall
assume respectively that
\ekv{Exp.8.3}
{h\ll \epsilon \le h^\delta ,}
and
\ekv{Exp.8.4}
{
h^2\ll \epsilon \le h^\delta,
}
for some $\delta >0$. Then in both cases, we have $\widetilde{q}={\cal O}(1)$.

\par We shall see how to eliminate the $t$-dependence by
conjugation with a \pop{} up to an exponentially decaying error.
 The problem can be attacked directly, but it seems
that we get better remainder estimates if we first reduce
ourselves to the case, when \ekv{Exp.8.5} {g(\tau )=\tau .} This
is possible by means of a \hol{} functional calculus. Indeed, let
$f=g^{-1}$ be the inverse of the map $g$. It is easy to see that
$f(P)$ is well-defined in the sense of formal analytic $h$-\pop{}s
or equivalently in the sense of composition of \clas{}s. (When $r$
is merely assumed to have an \asy{} expansion in powers of $h$, we
consider those $h$s as additional parameters.) We also see that
$f(P)$ as a symbol has the same properties as $P$ above, but now
with $g$ given by (\ref{Exp.8.5}). It will also be easy to return
to the original $P$, for if ${\rm Ad}_AP=e^APe^{-A}$, then at
least formally, ${\rm Ad}_Af(P)=f({\rm Ad}_AP)$, and to say that a
symbol is \indep{} of $t$ is equivalent to saying that the
corresponding \op{} commutes with translations in $t$ and this
latter property is stable under taking \hol{} \fu{}s of the \op{}.
Until further notice $g$ will be given by (\ref{Exp.8.5}).

\par
Using Proposition \ref{PropExp2}, we get
\begin{eqnarray}\label{Exp.15}
[P,A]&=&h{1\over i}\{ g,A\}(t,x,hD_{t,x};h)+R(A),\\
R(A)&=&\epsilon
[\widetilde{q},A],
\nonumber
\end{eqnarray}
where
\begin{equation}\label{Exp.16}
\3 R(A)\3_\rho \le  C(\widetilde{q})\epsilon \rho \3 A\3_\rho ,
\end{equation}
assuming that the $\3\cdot\3_\rho $-quasi-norm is chosen as in
Proposition \ref{PropExp2}.

\par Consider the map
\ekv{Exp.19}
{
A\mapsto {\rm Ad}_A(P)=e^APe^{-A},
}
where $\3 A\3_\rho $ is supposed to be small.
At least formally,
\ekv{Exp.19.5}
{
{\rm Ad}_A(P)=e^{{\rm ad}_A}(P).
}
We expand
\ekv{Exp.20}
{
{\rm Ad}_A(P)=\sum_{0}^\infty  {1\over k!}({\rm ad}_A)^k(P),
}
and get the expression for the differential
\ekv{Exp.21}
{
\delta A\mapsto {\rm ad}_{\delta A}(P)+\sum_{k=2}^\infty  {1\over
k!}\sum_{\nu +\mu =k-1\atop \nu ,\mu \ge 0}({\rm ad}_A)^\nu {\rm ad}_{\delta
A}({\rm ad}_A)^\mu (P).
}
An application of Proposition 3.2 shows that the $\rho $ quasi-norm of the last term can be estimated by
\begin{eqnarray*}
C(P)\rho \sum_{k=2}^\infty {1\over k!}\sum_{\nu +\mu =k-1}(2\3 A\3_\rho )
^{\nu +\mu }\3 \delta A\3_\rho &=& C(P)\rho (\sum_{k=2}^\infty {1\over
(k-1)!}(2\3 A\3_\rho ) ^{k-1})\3 \delta A\3_\rho \\
&=& C(P)\rho (e^{2\3 A\3_\rho }-1)\3 \delta A\3_\rho .
\end{eqnarray*}
So, if we assume some fixed upper bound on $\3 A\3_\rho $, we can represent
the differential of (\ref{Exp.19}) as
\ekv{Exp.22}
{
\delta A\mapsto {\rm ad}_{\delta A}(P)+K(A,\delta A),\quad \3 K(A,\delta
A)\3_\rho \le \widetilde{C}(P)\rho \3 A\3 _\rho \3\delta A\3_\rho ,
}
and combining this with (\ref{Exp.15}), (\ref{Exp.16}), we get
the expression for the differential:
\begin{eqnarray}\label{Exp.23}
&\delta A\mapsto -{h\over i}\{ g,\delta A\} +\widetilde{K}(A,\delta A)=
-{h\over i}g'(\tau )\partial _t\delta A +\widetilde{K}(A,\delta A),&\\
&\3 \widetilde{K}(A,\delta A)\3_\rho \le C\rho (\3
A\3_\rho +\epsilon )\3 \delta A\3_\rho .&\nonumber
\end{eqnarray}

\par Consider the linear \pb{}
\ekv{Exp.12}
{{1\over i}g'(\tau )\partial _tA=B-\langle B\rangle ,
}
where $B$ is a \clas{} of order 0 and
\ekv{Exp.11}
{
\langle B\rangle (\tau ,x,\xi ;h)={1\over 2\pi }\int_0^{2\pi }B(t,\tau
,x,\xi ;h)dt.
}
It has the solution
\ekv{Exp.13}
{
A={\cal L}(B)={\cal L}(B-\langle B\rangle ),\quad {\cal L}(B)=-{i\over
g'(\tau )}\int_0^{2\pi }({s\over 2\pi }-{1\over 2})B(t-s,\tau ,x,\xi
)ds.
}
Clearly (with a convenient choice of the families $\Omega _s$),
\ekv{Exp.14}
{
\3 {\cal L}(B)\3_\rho \le C_0\3 B-\langle B\rangle \3_\rho \le
\widetilde{C}_0\3 B\3_\rho .
}

Choose
\ekv{Exp.17}
{
\rho \ll {h\over \epsilon }.
}
We look for $A$ of the form $\sum _0^\infty A_j$ such that ${\rm Ad}_A(P)$
is \indep{} of $t$. We shall do the construction by successive
approximations in a such a way that uniformly during all the steps,
\ekv{Exp.24}
{
\3 A\3_\rho \le {\cal O}(1) ,\ \rho \3 A\3_\rho \ll h,\hbox{ and hence }\3
\widetilde{K}(A,\delta A)\3_\rho \le h\theta \3 \delta A\3_\rho ,\hbox{
where }\theta \ll 1.
}
To start with, assume, as we may, that $r-\langle r\rangle ={\cal O}(h^2)$
in the sense of ordinary symbols. Choose
$$A_0={\cal L}(r),\quad \3 A_0\3_\rho  \le C_0\3 r-\langle r\rangle
\3_\rho ={\cal O}(h^2).$$
Then, using (\ref{Exp.23}), (\ref{Exp.24}),
$$
{\rm Ad}_{A_0}(P)=P-hr+h\langle r\rangle +hr_1,
$$
where $P-hr+h\langle r\rangle $ is \indep{} of $t$ and $\3 r_1\3_\rho \le
\theta \3 r-\langle r\rangle \3_\rho $.

\par Put
$$
A_1={\cal L}(r_1),\quad \3 A_1\3_\rho \le C_0\theta \3 r-\langle
r\rangle \3_\rho .
$$
Then
$${\rm Ad}_{A_0+A_1}(P)={\rm Ad}_{A_0}(P)-h(r_1-\langle r_1\rangle )+hr_2=
g+\epsilon q+h\langle r\rangle +h\langle r_1\rangle +hr_2,$$
 and
$$
\3 r_2\3_\rho \le \theta \3 r_1\3_\rho \le \theta ^2\3 r-\langle r\rangle
\3_\rho .
$$
Since $\theta \ll 1$ the procedure will converge geometrically and we get
a formal solution $A$ with $\3 A\3_\rho \le C_1\3 r-\langle r\rangle
\3_\rho $.
\medskip

\par By construction $\3 A\3_\rho ={\cal O}(h^2)$ for $0<\rho \ll \min
(1,h/\epsilon  ) $ and we have defined ${\rm Ad}_A(P)$ by (\ref{Exp.20}).
We define $e^{tA}$
as a formal \an{} symbol of order 0, by
$$e^{tA}=\sum_{k=0}^\infty {t^kA^k\over k!},$$
so that
$$
\3 e^{tA}\3_\rho \le e^{\3 A\3 _\rho },\quad \partial
_te^{tA}=Ae^{tA}=e^{tA}A,\ e^{0A}=1
$$
in the space of formal symbols.
Similarly, we see that
$$
F_t={\rm Ad}_{tA}=e^{t{\rm ad}_A}=\sum _0^\infty  {t^k{\rm ad}_A^k\over k!}
$$
satisfies $\partial _tF_t={\rm ad}_A\circ F_t$, $F_0={\rm id}$. We then
verify (\ref{Exp.19.5}) simply by noticing that
$$
\partial _t(e^{tA}Pe^{-tA})={\rm ad}_A(e^{tA}Pe^{-tA})
$$

\par Notice that from the fact that
$$
\3 A\3_\rho ={\cal O}(h^2),\quad A=\sum_0^\infty  h^\nu A_\nu ,\ a_\nu
=A_\nu (1),
$$
we infer that locally $\vert a_\nu \vert \le f_\nu (A_\nu )(C\nu )^\nu ,$
with $\sum f_\nu (A_\nu )\rho ^\nu <{\cal O}(h^2) $. Hence $\vert a_\nu \vert \le
C_0h^2(C\nu /\rho )^\nu $, so
$$\vert h^\nu a_\nu \vert \le C_0h^2 (C\nu )^\nu ({h\over \rho })^\nu ,$$
so $A$ is an \an{} symbol with $h$ replaced by $h/\rho $ and can
be realized with an uncertainty ${\cal O}(h^2)e^{-1/(Ch/\rho
))}={\cal O}(h^2)e^{-\rho /(Ch)}$---see~\cite{Sj6}. Taking $\rho $
as large as possible respecting (\ref{Exp.17}), and recalling that
we also assume that $\rho $ is \bdd{}, we get the uncertainty
\ekv{Exp.25} { {\cal O}(h^2) e^{-1/C(\epsilon +h)}. }

\par This discussion can also be applied with $A$ replaced by $e^A$. Let
$B$ be such a realization of $e^A$. Then we can view $A$ as an \an{} symbol of
order $0$ by declaring that $h$ is an \indep{} parameter. Let $B^{-1}$ be a
parametrix, so that $B\# B^{-1}=1+{\cal O}(e^{-1/(Ch)})$, if we also denote
by $B^{-1}$ a realization. From the construction, it follows that $B\circ
P\circ B^{-1}$ has a symbol which is $t$-\indep{} up to an error of the
size (\ref{Exp.25}). With this in mind, we can state

\begin{prop}\label{PropExp3}
We can construct $A$ in Proposition {\rm \ref{PropRed1}}, such
that the symbol $\widehat{P}_\epsilon $ there is \indep{} of $t$
up to an error which is ${\cal O}(1)\exp\left(-1/C(\epsilon
+h)\right)$. Here we assume {\rm (\ref{Exp.8.3})} in the general
case and {\rm (\ref{Exp.8.4})} in the case when the subprincipal
symbol of $P_{\epsilon =0}$ vanishes.
\end{prop}

\par
We end this section with a heuristic discussion explaining why we
eventually will assume that the subprincipal symbol of
$P_{\epsilon =0}$ is zero. After decoupling by eliminating the
$t$-dependence as above, we get a \fy{} of 1-dimensional \op{}s
(\ref{Red.5}) that we consider at the branching level. If we first
consider the general case without any assumptions on the
subprincipal symbol, we can expect to have an estimate on the
inverse of these \op{}s or on the associated Grushin \pb{}s
roughly of the order \ekv{Exp.26} { \exp \left(C\left({\epsilon
\over h}+{h/\epsilon \over h}\right)\right)=\exp
\left(C\left({\epsilon \over h}+{1\over \epsilon }\right)\right).
} In order to combine everything, we would like this quantity
times (\ref{Exp.25}) to be $\ll 1$. This is clearly not the case.

\par In the case when the subprincipal symbol of $P_{\epsilon =0}$ vanishes
we expect to improve the estimate on
the 1-dimensional resolvents or inverses of Grushin problems to roughly
\ekv{Exp.27}
{\exp \left(C\left({\epsilon \over h}+{h^2/\epsilon \over
h}\right)\right)=\exp \left(C\left({\epsilon \over
h}+{h\over \epsilon }\right)\right).}
This leads to the condition
$$
{\epsilon \over h}+{h\over \epsilon }\ll {1\over \epsilon +h},$$
which simplifies to the condition \ekv{Exp.28} { \epsilon ^2\ll
h,\ h^2\ll \epsilon , } which in addition to (\ref{Exp.8.4}) gives
\ekv{Exp.29}{h^2\ll \epsilon \ll h^{1/2}.} This is the condition
on $\epsilon$ stated in Theorem 1.1.

\section{Transition matrix at a
branching level.}\label{SectionTr}
\setcounter{equation}{0}

In this section and the following one, we study certain model
\pb{}s. Much of the material is standard and close for instance to
\cite{HeSj2},~\cite{Ra}, and \cite{FuRa} (see also \cite{CoPa} for
the $C^\infty $-case), but our setup is somewhat different, and we
need to recollect some of the basic facts before returning to our
\op{} $P_\epsilon $.

\par Consider
\ekv{tr.1}{
({1\over 2}(xD_x+D_xx)-\alpha )u=0,
}
or equivalently the \e{}
$$x{\partial \over \partial x}u=(i\alpha
-{1\over 2})u.$$ From Proposition 11 in~\cite{Ra} we recall that
the solutions of (\ref{tr.1}) in ${\cal D}'(\real)$ form a
2-dimensional subspace of ${\cal S}'({\real})$.
\par\noindent For $x>0$, we express $u$ as
$u_1x^{i\alpha -{1\over 2}}$,
\par\noindent For $x<0$, we express $u$ as
$u_3\vert x\vert ^{i\alpha -{1\over 2}}$,
\par\noindent For $\xi >0$, we express
$\widehat{u}(\xi )$ as
$u_2\xi^{-i\alpha -{1\over 2}}$,
\par\noindent For $\xi <0$, we express
$\widehat{u}$ as
$u_4\vert \xi \vert ^{-i\alpha -{1\over 2}}$.
\par Here $\widehat{u}(\xi )={\cal F}u(\xi )={1\over
\sqrt{2\pi }}\int e^{-ix\xi }u(x)dx$ is the \F{} \tf{} and we
observe that (\ref{tr.1}) is equivalent to
\ekv{tr.2}
{({1\over 2}(\xi D_\xi +D_\xi \xi )+\alpha )\widehat{u}=0.}
If $\vert \Im \alpha\vert  < 1/2 $, we have two
solutions $u=U_\pm$ of (\ref{tr.1}), given by
\ekv{tr.3}
{
\widehat{U}_+(\xi )=H(\xi )\xi ^{-{1\over 2}-i\alpha },\
\widehat{U}_-(\xi )=H(-\xi ) \vert \xi \vert ^{-{1\over
2}-i\alpha }, }
where $H=1_{[0,+\infty [}$ is the Heaviside \fu{}, and the
general solution  to (\ref{tr.1}) becomes a linear combination,
\ekv{tr.4}
{u=u_2U_++u_4U_-.}
We see that $U_+$ is the \bdy{} value of a \hol{} \fu{} in
the upper half-plane that we also denote by $U_+$, and for
$x=iy$, $y>0$, we get
$$U_+(iy)={1\over \sqrt{2\pi }}\int_0^\infty e^{-y\xi }\xi
^{-{1\over 2}-i\alpha }d\xi ={1\over \sqrt{2\pi }}\Gamma
({1\over 2}-i\alpha )y^{i\alpha -{1\over 2}}.$$
Thus for real $x$,
$$U_+(x)={1\over \sqrt{2\pi }}\Gamma ({1\over 2}-i\alpha
)\big( {x+i0\over i}\big) ^{i\alpha -{1\over 2}},$$
which gives
\ekv{tr.5}
{
U_+(x)={1\over \sqrt{2\pi }}\Gamma ({1\over 2}-i\alpha
)\times \cases{e^{{\pi \over 2}\alpha +i{\pi \over
4}}x^{i\alpha -{1\over 2}},\ x>0,\cr
e^{-{\pi \over 2}\alpha -i{\pi \over
4}}\vert x\vert ^{i\alpha -{1\over 2}}, \ x<0.}
 }
Similarly, $U_-(x)=U_-(x-i0)$, with
$$U_-(-iy)={1\over {\sqrt{2\pi}}}\int_{-\infty }^0 e^{y\xi
}\vert \xi \vert ^{-i\alpha -{1\over 2}}d\xi ={1\over
\sqrt{2\pi }}\Gamma ({1\over 2}-i\alpha )y^{i\alpha -{1\over
2}},\ y>0,$$
so
$$U_-(x)={1\over \sqrt{2\pi }}\Gamma ({1\over 2}-i\alpha
)(i(x-i0))^{i\alpha -{1\over 2}},$$
\ekv{tr.6}
{
U_-(x)={1\over \sqrt{2\pi }}\Gamma ({1\over 2}-i\alpha
)\times \cases{e^{-{\pi \over 2}\alpha -i{\pi \over
4}}x^{i\alpha -{1\over 2}},\ x>0,\cr
e^{{\pi \over 2}\alpha +i{\pi \over
4}}\vert x\vert ^{i\alpha -{1\over 2}}, \ x<0.
} }

\par Now let $u$ be a solution of (\ref{tr.1}), so that
\ekv{tr.7} {u=u_1H(x) x^{-{1\over 2}+i\alpha }+u_3 H(-x)\vert
x\vert ^{-{1\over 2}+i\alpha }=u_2U_++u_4U_-.} Using (\ref{tr.5})
and (\ref{tr.6}), we get \eekv{tr.8} { u_1={1\over \sqrt{2\pi
}}\Gamma ({1\over 2}-i\alpha )e^{{\pi \over 2}\alpha +i{\pi \over
4}}u_2+{1\over \sqrt{2\pi }}\Gamma ({1\over 2}-i\alpha )e^{-{\pi
\alpha \over 2}-i{\pi \over 4}}u_4, } { u_3={1\over \sqrt{2\pi
}}\Gamma ({1\over 2}-i\alpha )e^{-{\pi \over 2}\alpha -i{\pi \over
4}}u_2+{1\over \sqrt{2\pi }}\Gamma ({1\over 2}-i\alpha )e^{{\pi
\alpha \over 2}+i{\pi \over 4}}u_4. } Here we want to express
$u_2,u_1$ in terms of $u_3,u_4$. From (\ref{tr.8}), we get
\ekv{tr.9} { \pmatrix{u_2\cr u_1}=\pmatrix{{\sqrt{2\pi }\over
\Gamma ({1\over 2}-i\alpha )}e^{{\pi \over 2}\alpha +i{\pi \over
4}} &-e^{\pi \alpha +i{\pi \over 2}}\cr e^{\pi \alpha +i{\pi \over
2}} &{\sqrt{2\pi }\over \Gamma ({1\over 2}+i\alpha )}e^{{\pi \over
2}\alpha -i{\pi \over 4}}}\pmatrix{u_3\cr u_4}, } where we also
used the reflection identity, \ekv{tr.10} {\Gamma ({1\over
2}+iz)\Gamma ({1\over 2}-iz)={\pi \over \cosh \pi z}.}

\par Recall that $\Gamma (z)$ is meromorphic with simple
poles at $-k$, for $k\in{\bf N}$, and no other poles. The
reflection identity above can also be written
$$\Gamma (z)\Gamma (1-z)={\pi \over \sin \pi z},$$
and implies that if $\Gamma (z)=0$, then $1-z$ has to be
pole, so $1-z=-k$ for some $k\in{\bf N}$, $z=k+1$, which is
impossible since we also know that $\Gamma (k+1)=k!\ne 0$.
Hence $\Gamma (z)$ has no zeros, and ${1\over \Gamma (z)}$ is
entire. The transition matrix in (\ref{tr.9}) has determinant 1 and
is an entire holomorphic function of $\alpha$. The relation
(\ref{tr.9}) remains valid therefore for all $\alpha \in \comp$.

\par We next compute the transition matrix analogous to the
one in (\ref{tr.9}) in the semi-classical case, for solutions of
\ekv{tr.11}
{({1\over 2}(xhD_x +hD_x x)-\mu )u=0.}
This is of course the same equation as (\ref{tr.1}) with $\alpha
=\mu /h$. We now require
\begin{eqnarray*}
u&=&u_1x^{i{\mu \over h}-{1\over 2}},\ x>0,\\
u&=&u_3\vert x\vert ^{i{\mu \over h}-{1\over 2}},\ x<0,\\
{\cal F}_hu(\xi )&=&u_2\xi ^{-i{\mu \over h}-{1\over 2}},\
\xi >0\\ {\cal F}_hu(\xi )&=&u_4\vert \xi \vert ^{-i{\mu \over
h}-{1\over 2}},\ \xi <0
\end{eqnarray*}
where \begeq \label{tr.11.5} {\cal F}_hu(\xi )={1\over \sqrt{2\pi
h}}\int e^{-ix\xi /h}u(x)dx={1\over \sqrt{h}}\widehat{u}({\xi
\over h}).
\endeq
A simple computation gives \ekv{tr.12} { \pmatrix{u_2\cr
u_1}=\pmatrix{{\sqrt{2\pi }\, h^ {i{\mu \over h}}\over \Gamma
({1\over 2}-i{\mu \over h})}e^{{\pi \over 2}{\mu \over h} +i{\pi
\over 4}} &-e^{\pi {\mu \over h}+i{\pi \over 2}}\cr e^{\pi {\mu
\over h}+i{\pi \over 2}} &{\sqrt{2\pi }\, h^{-i{\mu \over h}}\over
\Gamma ({1\over 2}+i{\mu \over h})}e^{{\pi \over 2}{\mu \over
h}-i{\pi \over 4}}}\pmatrix{u_3\cr
u_4}=\pmatrix{a_{2,3}&a_{2,4}\cr a_{1,3}&a_{1,4}}\pmatrix{u_3\cr
u_4} .}

We summarize the discussion above in the following proposition.

\begin{prop}
\label{PropTr1} Let $\mu\in \comp$ be such that $\abs{\Im
\mu}<h/2$. If $u\in {\cal D}'(\real)$ is a solution of {\rm
(\ref{tr.11})} then $u\in {\cal S}'(\real)$ and there exist $u_1$,
$u_2$, $u_3$, $u_4\in \comp$ such that
$$
u=u_1 H(x) x^{i\frac{\mu}{h}-\frac{1}{2}}+u_2 {\cal F}_h^{-1}
\left(H(\xi) \xi^{-i\frac{\mu}{h}-\frac{1}{2}}\right)=u_3 H(-x)
\abs{x}^{i\frac{\mu}{h}-\frac{1}{2}}+u_4 {\cal
F}_h^{-1}\left(H(-\xi)\abs{\xi}^{-i\frac{\mu}{h}-\frac{1}{2}}\right).
$$
Here ${\cal F}_h$ is the semiclassical Fourier transform defined
in {\rm (\ref{tr.11.5})} and the coefficients $u_2$, $u_1$ can be
expressed in terms of $u_3$, $u_4$ by {\rm (\ref{tr.12})}. The
transition matrix which occurs in {\rm (\ref{tr.12})} is entire
holomorphic in $\mu$ and has determinant $1$.
\end{prop}

\par We finish this section by the following observation, which will be useful in section 6.
The \op{} $P_0={1\over 2}(xhD_x+hD_xx)$ has the principal symbol
$p_0(x,\xi )=x\xi $. For $\mu \in {\bf C}$, $\vert \mu \vert \ll
1$, define $\rho _j\in p_0^{-1}(\mu )$ by $\rho _1=(1,\mu )$,
$\rho _2=(\mu ,1)$, $\rho _3=(-1,-\mu )$, $\rho _4=(-\mu ,-1)$.
Working in the semiclassical limit, define the microlocal null
solutions $e_j$ of $P_0-\mu$, for $j=1,..,4$ by $e_1=x^{i{\mu
\over h}-{1\over 2}}$ near $\rho _1$, $e_1=0$ near $\rho _3$,
${\cal F}_he_2=\xi ^{-i{\mu \over h}-{1\over 2}}$ near $\rho
_1\approx \kappa _{{\cal F}_h}\rho _2$, $e_2=0$ near $\rho _4$,
$e_j(x)=e_{j-2}(-x)$, $j=3,4$. Here $\kappa_{{\cal F}_h}$ is the
map $(x,\xi)\mapsto (\xi,-x)$ associated to ${\cal F}_h$. Then a
general null solution of $P_0-\mu $ can be written either as
$u_2e_2+u_1e_1$ or as $u_3e_3+u_4e_4$, where (\ref{tr.12}) holds.

\section{Asymptotics of the
transition matrix}\label{SectionAs}
\setcounter{equation}{0}
\medskip
\par In this section, we shall derive asymptotic formulas for the entries of the transition matrix (\ref{tr.12}).
We shall use the following version of Stirling's formula
(\cite{Ol}), \ekv{as.1} {{\Gamma (z)\over \sqrt{2\pi
}}=e^{-z}z^{z-{1\over 2}}(1+{1\over 12z}+{1\over 288 z^2}+...),\
\vert z\vert \to \infty ,\, \vert {\rm arg\,}z\vert \le \pi
-\delta ,} for every fixed $\delta >0$.  We apply this in two
cases:
\smallskip

\noindent {\it Case A.} We have $\vert \mu \vert /h\gg 1$ and $\mu
\not\in$ a conic \neigh{} of the negative imaginary axis. Then
${1\over 2}-i{\mu \over h}$, $-i{\mu \over h}$ avoid a conic
\neigh{} of ${\bf R}_-$ and we can apply (\ref{as.1}), to get:
\begin{eqnarray*} {1\over \sqrt{2\pi }}\Gamma ({1\over 2}-i{\mu \over
h})&=&e^{-{1\over 2}+i{\mu \over h}}({1\over 2}-i{\mu \over
h})^{-i{\mu \over h}}e^{{\cal O}({h\over \mu })} \\ &=&\exp
\left[-{1\over 2}+i{\mu \over h}-i{\mu \over h} \ln ({1\over
2}-i{\mu \over h})+{\cal O}({h\over \mu } )\right]\\ &=&\exp
\left[-{1\over 2}+i{\mu \over h}-i{\mu \over h}(\ln (-i{\mu \over
h})-{h\over 2i\mu })+{\cal O}({h\over \mu })\right]\\ &=&\exp
\left[-{1\over 2}+i{\mu \over h}-i{\mu \over h}\ln (-i{\mu \over
h})+{1\over 2}+{\cal O}({h\over \mu })\right],\end{eqnarray*} so
in this case we have \ekv{as.2_-} {{1\over \sqrt{2\pi }}\Gamma
({1\over 2}-i{\mu \over h})=\exp \left[{i\mu \over h}-i{\mu \over
h}\ln (-i\mu )+{i\mu \over h}\ln h+{\cal O}_-({h\over \mu
})\right].} Here and in what follows $\ln$ always stands for the
principal branch of the logarithm.

\smallskip \noindent
{\it Case B.} We have $\vert \mu \vert /h\gg 1$ and $\mu $ avoids
some conic \neigh{} of the positive imaginary axis. Then we can
apply the earlier results with $\mu $ replaced by $-\mu $ and get:
\ekv{as.2_+} {{1\over \sqrt{2\pi }}\Gamma ({1\over 2}+i{\mu \over
h})=\exp \left[-{i\mu \over h}+i{\mu \over h}\ln (i\mu )-{i\mu
\over h}\ln h+{\cal O}_+({h\over \mu })\right].}
\smallskip

\par If $\Re
\mu \ge {1\over C}\vert \Im \mu \vert $, we combine this with the
reflection identity:
$${1\over 2\pi }\Gamma ({1\over 2}-i{\mu \over h})\Gamma ({1\over 2}+i{\mu
\over h})={1\over 2\cosh {\pi \mu \over h}}= e^{-{\pi \mu \over
h}+{\cal O}(e^{-2\pi \Re \mu /h})},$$
and the fact that
$$\ln (i\mu )-\ln
(-i\mu )=i\pi $$ in this region, to conclude that \ekv{as.3}
{{\cal O}_+({h\over \mu })+{\cal O}_-({h\over \mu })= {\cal
O}(e^{-2\pi \Re \mu /h}).} It is straightforward to establish a
corresponding estimate in the region $\Re \mu \le -{1\over C}\vert
\Im \mu \vert $, and we can summarize both cases in: \ekv{as.4}
{{\cal O}_+({h\over \mu })+{\cal O}_-({h\over \mu })={\cal
O}(e^{-2\pi \vert \Re \mu \vert /h}),\ \vert \Re \mu \vert \ge
{1\over C}\vert \Im \mu \vert .}
\medskip
\par\noindent \it Remark. \rm When $\mu $ is real we have
\ekv{as.4.1}
{\Re {\cal O}_+({h\over \mu })=\Re {\cal O}_-({h\over \mu })={\cal
O}(e^{-2\pi \vert \mu \vert /h}).}
In fact, if we first assume that $\mu \gg h$, we get from (\ref{as.2_-}),
(\ref{as.2_+}):
$${\Gamma ({1\over 2}-i{\mu \over h})\over \sqrt{\pi }}=\exp \left[i{\mu \over
h}-i{\mu \over h}\ln \mu +i{\mu \over h}\ln h-{\pi \mu \over
2h}+{\cal O}_-({h\over \mu })\right],$$
$${\Gamma ({1\over 2}+i{\mu \over h})\over \sqrt{\pi }}=\exp \left[-i{\mu \over
h}+i{\mu \over h}\ln \mu -i{\mu \over h}\ln h-{\pi \mu \over
2h}+{\cal O}_+({h\over \mu })\right],$$ and using that the second
quantity is equal to the complex conjugate of the other, we
conclude that (\ref{as.4.1}) holds in this case. In the case $\mu
\ll -h$, we can use the same argument. In this case
(\ref{as.2_-}), (\ref{as.2_+}) give
$${\Gamma ({1\over 2}-i{\mu \over h})\over \sqrt{\pi }}=\exp \left[i{\mu \over
h}-i{\mu \over h}\ln \vert \mu\vert  +i{\mu \over h}\ln h+{\pi \mu
\over 2h}+{\cal O}_-({h\over \mu })\right],
$$
$${\Gamma ({1\over 2}+i{\mu \over h})\over \sqrt{\pi }}=\exp [-i{\mu \over
h}-i{\mu \over h}\ln \vert \mu\vert  -i{\mu \over h}\ln h+{\pi \mu \over 2h}+{\cal
O}_+({h\over \mu })],$$
and we can again conclude that (\ref{as.4.1}) holds.
\medskip

\par In any closed sector away from $i{\bf R}_-$, we can use (\ref{tr.12}) and (\ref{as.2_-}) to get
\ekv{as.5} {a_{2,3}=\exp \left({i\over h}\left(\mu \ln (-i\mu
)-i{\pi \mu \over 2}-\mu +{h\pi \over 4}+i h{\cal O}_-({h\over \mu
})\right)\right).}

\par  In any closed sector away from $i{\bf
R}_+$, we can use (\ref{tr.12}), the reflection identity and
$(\ref{as.2_+})$, to get \ekv{as.6} {a_{2,3}=2\cosh ({\pi\mu \over
h})e^{{i\over h}(\mu \ln (i\mu )-{i\pi \mu \over 2}-\mu +{\pi
h\over 4}-i h{\cal O}_+({h\over \mu }))}}

\par  Using (\ref{as.2_+}), (\ref{tr.12}), we
get for $\mu $ away from a sector around $i{\bf R}_+$: \ekv{as.7}
{a_{1,4}=e^{{i\over h}(-\mu \ln (i\mu )-i{\pi\mu  \over 2}+\mu
-{h\pi \over 4}+i h{\cal O}_+({h\over \mu }))}} In a sector $\Im
\mu >-C\vert \Re \mu \vert $, we use the reflection identity
$${\sqrt{2\pi }\over \Gamma ({1\over 2}+i{\mu \over h})}={\Gamma
({1\over 2}-i{\mu \over h})\over \sqrt{2\pi }}2\cosh \pi {\mu
\over h},$$ to get \ekv{as.8} { a_{1,4}= 2\cosh ({\pi\mu \over h}
) e^{{i\over h}(-\mu \ln ({\mu \over i})-i{\pi \mu \over 2}+\mu
-{h\pi \over 4}-i h{\cal O}_-({h\over \mu }))}.}

Combining the asymptotic formulae (\ref{as.5})--(\ref{as.8}), we
may summarize the discussion in this section in the following
proposition.

\begin{prop}
We have the following tableau for the coefficients of the
transition matrix {\rm (\ref{tr.12})}, when $\vert \mu \vert /h\gg
1$. In all cases:
$$ a_{2,4}=-e^{\pi {\mu \over h}+i{\pi \over 2}},\ a_{1,3}=e^{\pi
{\mu \over h}+i{\pi \over 2}}$$
For $\Re \mu >C^{-1}\vert \Im \mu
\vert $:
$$a_{2,3}=
e^{{i\over
h}(\mu \ln \mu -i\pi \mu -\mu +{\pi h\over 4} +ih{\cal O}_-({h\over
\mu }))},$$
$$a_{1,4}= e^{{i\over h}(-\mu \ln \mu -i\pi \mu +\mu -{\pi h\over 4}
+ih{\cal O}_+({h\over \mu }))}.$$ For $\Im \mu >-C\vert \Re \mu \vert $: $$
a_{2,3}= e^{{i\over h}(\mu \ln {\mu\over i} -i{\pi\mu\over 2 } -\mu +{\pi
h\over 4}+ih{\cal O}_-({h\over \mu }))},$$
$$a_{1,4}=
e^{{i\over h}(-\mu \ln {\mu\over i}-i{\pi
\mu\over 2}\pm i\pi \mu +\mu -{\pi h\over 4} -ih{\cal O}_-({h\over \mu }))}.
$$
For $\Re \mu
<-C^{-1}\vert \Im \mu \vert $:
$$ a_{2,3}=
e^{{i\over h}(\mu \ln (-\mu ) -\mu +{\pi h\over 4}+ih{\cal
O}_-({h\over \mu }))},$$
$$a_{1,4}=
e^{{i\over h}(-\mu \ln (-\mu )+\mu -{\pi h\over 4}+ih{\cal
O}_+({h\over \mu }))}.  $$
For $\Im \mu <C\vert \Re \mu \vert $:
$$
a_{2,3}=
e^{{i\over h}(\mu \ln (i\mu )
-i{\pi \mu \over 2}\pm i\pi \mu -\mu +{\pi h\over 4}-ih{\cal O}_+ ({h\over
\mu }))}, $$
$$ a_{1,4}=
e^{{i\over h}(-\mu \ln (i\mu )-i{\pi\mu \over 2}+\mu -{\pi h\over
4}+ih{\cal O}_+({h\over \mu }))}.  $$

Here the terms with $\pm i\pi \mu $ in the exponents indicate that
we should take the sum of the two possible terms with the same
remainders ${\cal O}_+$ or ${\cal O}_-$ in the exponent for each
of the two terms.
\end{prop}

\section{The one-dimensional spectral \pb }\label{Section1D}
\setcounter{equation}{0}

We now return to Proposition \ref{PropRed1}, which shows (when
combined with the exponential remainder estimates of Section
\ref{SectionExp}), that the study of $P_\epsilon $ near $\Lambda
_{0,0}$ (considered in a suitable weighted space) can be reduced
by conjugation to that of $\widehat{P}_\epsilon $ acting on the
space $L^2_\theta (S^1\times {\bf R})$ of \fu{}s defined \ml{}ly
in some fixed \neigh{} of $\tau =0$, $(x,\xi )\in K_{0,0}\subset
{\bf R}^2$ in $T^*S^1\times T^*{\bf R}$, with a $\theta =(\theta
_0,\theta _1,\theta _2)$ Floquet periodicity condition. Here
$K_{0,0}$ is an $\infty $-shaped curve as in the introduction, and
$\theta $ was defined in the beginning of Section
\ref{SectionRed}. In order to fix the ideas, we assume that we are
in the case when the subprincipal symbol of $P_{\epsilon =0}$
vanishes, so that the symbol of $\widehat{P}_\epsilon $ is given
by (\ref{Red.8}), and we are then in the parameter range $h^2\ll
\epsilon \le h^\delta $. At least formally, the study of
$\widehat{P}_\epsilon $ can be reduced to a \fy{} of
one-dimensional \pb{}s by a \F{} series expansion in the
$t$-variable and (as noted in (\ref{Red.5})) we get the
one-dimensional \op{}s \ekv{1D.1} {
\widehat{P}_{\epsilon}(h(k-\theta _0),x,hD_x;h),\ k\in{\bf Z}, }
where $\theta _0=S_0/(2\pi h)+k_0/4$, and we restrict the range of
$k$ by requiring that $h(k-\theta _0)$ be small. The \op{}s
(\ref{1D.1}) should be considered as acting on the \ml{} space
$L^2_{\theta '}({\bf R})$ defined similarly to $L^2_\theta
(S^1\times {\bf R})$, with $\theta '=(\theta _1,\theta _2)$. Using
(\ref{Red.8}), we see that (\ref{1D.1}) becomes \ekv{1D.2} {
g(\tau )+i\epsilon Q(\tau ,x,hD_x,\epsilon ,{h^2\over \epsilon
};h),\ \tau =h(k-\theta _0), } where
$$Q(\tau ,x,\xi ,\epsilon ,{h^2\over \epsilon } ;h)\sim
Q_0(\tau ,x,\xi ,\epsilon ,{h^2\over \epsilon })+hQ_1+h^2Q_2+...$$
is \hol{} \wrt{} $(\tau ,x,\xi )$ in a fixed complex \neigh{} of
$\{ 0\} \times K_{0,0}$ and depends smoothly on the other
parameters. We further notice that \ekv{1D.3} { Q_0(\tau ,x,\xi
,0,0)=\langle q\rangle (\tau ,x,\xi ) } is equal to the \traj{}
average of $q$, expressed in suitable real symplectic coordinates,
and we know by construction that \ekv{1D.4} { \langle q\rangle
(0,x,\xi )=0\hbox{ on }K_{0,0}. }

\par We also recall the assumptions (\ref{H2}), (\ref{H4}), which say that
\ekv{1D.5} { \langle q\rangle (\tau ,x,\xi )=f(\tau ,\Re \langle
q\rangle(\tau,x,\xi)),\ \Re f(\tau ,r)=r,\,\, f(0,0)=0} \ekv{1D.6}
{ \Re \langle q\rangle ''_{(x,\xi ),(x,\xi )}(0,0,0)\hbox{ is
\nondeg{} of signature }0. } Here we assume for simplicity that
$(0,0)\in K_{0,0}$ is the branching point.

\par In the following we shall discuss the spectrum of the 1-dimensional
\op{} \ekv{1D.7} {Q=Q_\tau =Q(\tau ,x,hD_x,\epsilon ,{h^2\over
\epsilon };h).} Since this \op{} is only defined \ml{}ly and up to
an error ${\cal O}\left(e^{-{1\over C(\epsilon +h)}}\right)$, we
have to keep in mind that for the moment the \ev{}s will be
defined only formally and up to errors of at least the same size.
\smallskip
\par\noindent \it A first localization of the spectrum. \rm Assume first
that $\langle q\rangle $ is real-valued. Then from the sharp G\aa{}rding
inequality, we see that the spectrum of the \op{} (\ref{1D.7}) in the
band $\vert \Re z\vert <1/{\cal O}(1)$ is contained in
\ekv{1D.8}
{
\left\{ z\in{\bf C};\, \vert \Re z\vert <1/{\cal O}(1),\,\, \vert \Im z\vert
\le {\cal O}(1)(h+\epsilon +{h^2\over \epsilon })\right\} .
}
Under the more general assumption (\ref{1D.5}), we see that (\ref{1D.8})
can be applied to $g(\tau ,Q_\tau )$, where $g(\tau ,\cdot )=f^{-1}(\tau
,\cdot )$. So in the general case, we see that the spectrum of the \op{}
(\ref{1D.7}) in the band $\vert \Re z\vert <1/{\cal O}(1)$ is contained in
\ekv{1D.9}
{
\Sigma _\tau :=\left\{ z\in {\bf C};\, \vert \Re z\vert <1/{\cal O}(1),\, \vert
z-f(\tau ,\Re z)\vert \le {\cal O}(1) (h+\epsilon +{h^2\over \epsilon })\right\} .
}
We also have
\ekv{1D.10}
{
\Vert (Q_\tau -z)^{-1}\Vert \le {{\cal O}(1)\over {\rm dist\,}(z,\Sigma
_\tau ) },\hbox{ for }\vert \Re z\vert <1/{\cal O}(1),\,\, z\not\in \Sigma
_\tau .
}
\smallskip
\par\noindent {\it Normal forms near the branching points.}
Let $Q_0=Q_0(\tau ,x,\xi ,\epsilon ,{h^2\over \epsilon })$ be the
principal symbol of $Q$. Following \cite{HiSj1},
\cite{HeSj}, we get the following adaptation of Proposition 5.3 of
\cite{HiSj1}:
\begin{prop}\label{Prop1D1}
We can find a \ctf{}: $(x,\xi )\mapsto \kappa _{\tau ,\epsilon
,{h^2\over \epsilon }}(x,\xi )$ depending \an{}ally on $\tau $ and
smoothly on $\epsilon ,h^2/\epsilon $ with values in the \hol{}
\ctf{}s: ${\rm neigh\,}((0,0),{\bf C}^2)\to {\rm
neigh\,}((0,0),{\bf C}^2)$, and an \an{} \fu{} $k_{\epsilon
,h^2/\epsilon }(\tau ,q)$ depending smoothly on $\epsilon
,h^2/\epsilon $, such that
$$
\kappa _{\tau ,\epsilon ,h^2/\epsilon }(0,0)=(x(\tau ,\epsilon
,h^2/\epsilon ),\xi (\tau ,\epsilon ,h^2/\epsilon ))
$$
is the unique critical point close to {\rm (0,0)} of $(x,\xi
)\mapsto Q_0(\tau ,x,\xi ,\epsilon ,h^2/\epsilon )$ and with
\ekv{1D.11} { Q_0(\tau ,\kappa _{\tau ,\epsilon ,h^2/\epsilon
}(x,\xi ),\epsilon ,{h^2\over \epsilon })=k_{\epsilon ,{h^2\over
\epsilon }}(\tau ,x\xi ). } Moreover, $\kappa _{\tau ,0,0}$ is
real when $\tau $ is real and
$$
{\partial \over \partial q}\Re k_{\epsilon ,{h^2\over \epsilon }}(\tau
,0)>0.
$$
\end{prop}

\par After a conjugation by an elliptic \fop{} associated to the \ctf{}
$\kappa _{\tau ,\epsilon ,h^2/\epsilon }$ we may assume that the
leading symbol of $Q_{\tau}$, \ekv{1D.12} { Q_0=Q_0(\tau ,x, \xi
,\epsilon ,h^2/\epsilon ) } is a \fu{} of $\tau$, $\epsilon$,
$h^2/\epsilon $ and $x\xi $. We can get a complete normal form by
making further conjugations by analytic \pop{}s of order 0 in such
a way that the complete symbol also becomes a \fu{} of $\tau$,
$\epsilon$, $h^2/\epsilon $ and $x\xi $. This is carried out in
Appendix A. We get the following result which is very close to one
of the main results of the appendix b in \cite{HeSj2}.

\begin{prop}\label{Prop1D2}
We can quantize $\kappa _{\tau ,\epsilon ,h^2/\epsilon }$ by an
elliptic \fop{} $U=U_{\tau ,\epsilon ,h^2/\epsilon }$ with an
analytic symbol, depending \hol{}ally on $\tau $ and smoothly on
$\epsilon ,h^2/\epsilon $, such that \ekv{1D.31} {
U^{-1}QU=K_{\epsilon ,h^2/\epsilon }(\tau ,I;h)+{\cal
O}(e^{-{1\over Ch}}),\ I=P_0={1\over 2}(x\circ hD_x+hD_x\circ x),
} where $K_{\epsilon ,h^2/\epsilon }(\tau ,\iota ;h)$ is a
classical analytic symbol of order $0$ depending \hol{}ally on
$\tau $ and smoothly on $\epsilon ,h^2/\epsilon $.  The leading
part of $K$ is $k_{\epsilon ,h^2/\epsilon }(\tau ,\iota )$
appearing in Proposition {\rm \ref{Prop1D1}}.
\end{prop}

\par\noindent \it The quantization condition. \rm We start with a side
remark about normalization. When $P(x,hD_x)$ is a \sa{} one-dimensional
$h$-\pop{} of  principal type, and $z\in{\bf R}$, then we can normalize
\ml{}ly defined solutions of $(P-z)u=0$, $z\in {\bf R}$, by imposing
that
$$
({i\over h}[P,\chi ]u\vert u)=1.
$$
here $\chi =\chi (x,\xi )$ is defined near a piece of the real
characteristics and has the property that $\nabla \chi $ is of compact
support near the characteristics of $P-z$ and $\chi $ increases from 0 to 1
when we progress in the Hamilton flow direction. (See \cite{HeSj2}.) It is easy to check that
if we view $u$ as a solution of $(f(P)-f(z))u=0$, then we get the
corresponding normalization
$$
({i\over h}[f(P),\chi ]u\vert u)=g(z,z)({i\over h}[P,\chi ]u\vert u)=g(z,z),
$$
where
$$
f(P)-f(z)=(P-z)g(P,z).
$$
If we drop the requirement that $P$ be \sa{} or just let $z$
become complex, there is no obvious normalization of
null-solutions of $P-z$, but we still have a well-defined
sesqui-linear form on ${\cal N}(P-z)\times {\cal
N}(P^*-\overline{z})$, given by
$$
({i\over h}[P,\chi ]u\vert v).
$$
If we have some additional information allowing us to identify the two
null-spaces, then this can still be used to normalize null-solutions of
$P-z$. In the following we abandon the attempt to normalize completely the
null-solutions, since already the \op{} $Q_{\epsilon =0,h^2/\epsilon =0}$
is not necessarily \sa{}.

By Proposition \ref{Prop1D2} we have an analytic symbol $f(\cdot
;h)$ depending analytically on $\tau $ and smoothly on $\epsilon
$, $h^2/\epsilon $, such that \ekv{1D.33} { U^{-1}f(Q;h)U=P_0. }
Notice that if $u$ is a null-solution of $P_0-\mu $ in a full
\neigh{} of $(0,0)$, then $(Q-z)Uu=0$ near $(0,0)$, where the
spectral parameters are related by \ekv{1D.34} { f(z;h)=\mu . }
Recall from the end of Section \ref{SectionTr} that $P_0-\mu $ has
the four characteristic points $\rho _j$, $j=1,2,3,4$ and that
this \op{} has the \ml{} null solutions $e_j$ described after
Proposition \ref{PropTr1}. When $\mu $ is real, we check that
$e_j$ is normalized near $\rho _j$. If $v$ is a global
null-solution of $P_0-\mu $, with $v=v_je_j$ near $\rho _j$, then
by Proposition 4.1 we have: \ekv{1D.35} {\pmatrix{v_2\cr
v_1}=\pmatrix{a_{2,3}&a_{2,4}\cr a_{1,3} &a_{1,4}}\pmatrix{v_3 \cr
v_4}.}

Assume for simplicity that $\kappa=\kappa_{\tau,\eps,h^2/\eps}$ is
defined in a suitable domain, containing $\rho _j$, $j=1,..,4$.
Let $\alpha _j=\kappa (\rho _j)$, $f_j=Ue_j$. Then if $u$ is
null-solution of $Q-z$ near the branching point $\kappa((0,0))$,
equal to $v_jf_j$ near $\alpha _j$, then (\ref{1D.35}) still
holds. We may wish to renormalize the $f_j$, by putting
\ekv{1D.38} { f_j=e^{{i\over h}d_j}g_j. } Then a straightforward
calculation shows that if $u$ is a null-solution of $Q-z$ near the
branching point, and $u=u_jg_j$ near $\alpha _j$, then \ekv{1D.41}
{ \pmatrix{u_2\cr u_1}=\pmatrix{c_{2,3}&c_{2,4}\cr
c_{1,3}&c_{1,4}}\pmatrix{u_3\cr u_4},\ c_{j,k}=e^{-{i\over
h}(d_j-d_k)}a_{j,k}. }

\par Here is a natural example of such a renormalization. Assume for
simplicity that near $\alpha _j$ the set $Q_0^{-1}(z)$ takes the
form $\xi =\lambda _j(x)$, where $\lambda _j$ is \an{} and depends
\an{}ally on the parameters $\epsilon ,h^2/\epsilon ,\tau ,z$.
Then choose $g_j$ so that \ml{}ly near $\alpha _j$ we have the
standard WKB-form: \ekv{1D.38'} { g_j=b_j(x;h)e^{{i\over h}\psi
_j(x)}, } where $b_j$ is a classical elliptic analytic symbol of
order 0. The function $\psi _j$ solves the eikonal \e{}
$$
{\partial \psi _j\over \partial x}-\lambda _j(x)=0,\hbox{ with the
extra condition }\psi _j(\pi _x(\alpha _j))=0,
$$
and $b_j,\psi _j$ depend \an{}ally on the additional parameters
$\tau ,z$ and smoothly on $\eps, h^2/\eps$. Using an explicit
representation of $U$ we write near $\alpha _j$ for $j=1,3$:
\ekv{1D.36} { f_j(x)=h^{-{1+N\over 2}}\iint e^{{i\over h}(\psi
(x,y,\theta )+\phi _j(y))}A(x,y,\theta )a_j(y)dyd\theta
=e^{{i\over h}\widetilde{\psi}_j}\widetilde{b}_j(x;h). } Here the
last equality follows from stationary phase, $\psi $ is a
\nondeg{} phase \fu{} generating $\kappa $ and near $\rho _j$ we
write $e_j$ in the WKB-form
$$
e_j(y)=a_j(y)e^{\frac{i}{h}\phi_j(y)}=\vert y\vert ^{{i\mu \over
h}-{1\over 2}}=\vert y\vert ^{-{1\over 2}}e^{{i\over h}\mu \ln
\vert y\vert }.
$$
The function $\widetilde{\psi }_j(x)$ in (\ref{1D.36}) appears as
the critical value in the stationary phase expansion of
(\ref{1D.36}) and solves the same eikonal equation as $\psi _j$.

\par For $j=2,4$, we get near $\alpha _j$:
\ekv{1D.36'}
{
f_j(x)=h^{-{2+N\over 2}}\iiint e^{{i\over h}(\psi (x,y,\theta )+y\eta +\phi
_j(\eta ))}A(x,y,\theta )a_j(\eta )dyd\eta d\theta =e^{{i\over h}\widetilde{\psi
}_j}\widetilde{b}_j(x;h),
}
where $a_j$ and $\phi _j$ appear when writing ${\cal F}_he_j$ on WKB-form
near $\kappa _{\cal F}(\rho _j)$.

\par In this case we see that $d_j=d_j(h)$ is a classical analytic symbol
of order 0, depending \an{}ally on the additional parameters and with the
imaginary part of the
leading symbol vanishing when $\Im \mu =0$, $\epsilon =0$, $h^2/\epsilon =0$. The
leading part of $d_j$ can be further described in terms of symplectic geometry.

\par Put
\ekv{1D.42}
{
\theta _{j,k}:=d_j-d_k.
}
We have the obvious relation
\ekv{1D.43}
{
\theta _{2,3}+\theta _{1,4}=\theta _{1,3}+\theta _{2,4}.
}

\par Now we work in a full \neigh{} of $K_{0,0}$ (introduced after (\ref{int.2})).
Recall that we have the points $\alpha_1$, $\alpha_2$, $\alpha_3$,
$\alpha_4$ on the four crossing branches of $K_{0,0}$ distributed
with positive orientation around the branching point. We may
assume that $\alpha_3$, $\alpha_4$ are situated close to the left
closed curve $\gamma^1$ of $K_{0,0}$ and that $\alpha_1$,
$\alpha_2$ are situated close to the right closed curve $\gamma^2$
of $K_{0,0}$. Start with a \ml{} null-solution to $Q-z$ near
$\kappa((0,0))$, of the form $u_4g_4$ near $\alpha _4$ and of the
form $u_3g_3$ near $\alpha _3$. Here for the moment $u_3,u_4$ can
be prescribed \ably{}, and we then know that $u=u_jg_j$ near
$\alpha _j$ for $j=2,1$, where $u_2$, $u_1$ are given by
(\ref{1D.41}). We require temporarily that $u$ is a well-defined
single-valued null-solution along the whole left closed component
$\gamma ^1$ of $K_{0,0}$. Then if we follow $u$ around the
exterior part of $\gamma^1$ from $\alpha _4$ to $\alpha _3$, we
get \ekv{1D.44} { u_3=e^{{i\over h}S_{3,4}}u_4,\hbox{ where
}S_{3,4}=\int_{\gamma _{3,4}}\xi dx+{\cal O}(h)=S^0_{3,4}+{\cal
O}(h), } with $\gamma _{3,4}$ denoting the exterior part of
$\gamma ^1$ which joins $\alpha _4$ to $\alpha _3$.

\par Now recall that we really want $u\in L^2_{\theta '}$, $\theta
'=(\theta _1,\theta _2)$, where $\theta _j={S_j\over 2\pi h}+{k_j\over 4}$
where $S_j$ is a real action difference related to the reduction in
Proposition \ref{PropRed1} and $k_j\in{\bf Z}$ a corresponding Maslov
index. This means that $u$ should be multivalued, but Floquet periodic
along $\gamma ^1$ in the sense that
\ekv{1D.45}
{
\gamma _*^1u=e^{-2\pi i\theta _1}u,
}
where $\gamma _*^1u$ denotes the extension of $u$ along one loop of
$\gamma ^1$ which we assume to be oriented in the following way: $\alpha
_4\to \alpha _3\to (0,0)\to \alpha _4$. Starting near $\alpha _3$, we
get $\gamma _*^1u$ near the same point in two steps:
$$u_3g_3\to u_4g_4\to e^{{i\over h}S_{3,4}}u_4g_3.$$ The Floquet condition
(\ref{1D.45}) therefore becomes $e^{-2\pi i\theta
_1}u_3=e^{iS_{3,4}/h}u_4$, or equivalently
\ekv{1D.46}
{
u_3=e^{2\pi i\theta _1+{i\over h}S_{3,4}}u_4,
}
instead of (\ref{1D.44}).

\par Similarly, let $\gamma ^2$ be the right hand loop in $K_{0,0}$ with the
orientation: $\alpha _2\to \alpha _1\to (0,0)\to \alpha _2$. Then,
if we want $u$ to extend to a null-solution in $L_{\theta '}^2$
near $\gamma ^2$, we get the analogue of (\ref{1D.46}):
\ekv{1D.47} { u_1=e^{2\pi i\theta _2 +{i\over h}S_{1,2}}u_2, }
with $S_{1,2}$ defined as in (\ref{1D.44}) with $\gamma _{3,4}$
there replaced by $\gamma _{1,2}$, the exterior segment in $\gamma
^2$ that joins $\alpha _2$ to $\alpha _1$.

\par Start near $\alpha _4$ with $u_4g_4$, use (\ref{1D.46}) to get $u_3$
and then (\ref{1D.41}) to get $u_2,u_1$: \ekv{1D.48} {
\cases{u_2=(c_{2,3}e^{2\pi i\theta _1+{i\over
h}S_{3,4}}+c_{2,4})u_4,\cr u_1=(c_{1,3}e^{2\pi i\theta _1+{i\over
h}S_{3,4}}+c_{1,4})u_4}, } and in order to get a global solution
in $L^2_{\theta '}$, we also need to apply (\ref{1D.47}), which
gives our global one-dimensional quantization condition
\ekv{1D.49} { 0=c_{2,3}e^{2\pi i(\theta _1+\theta _2)+{i\over
h}(S_{3,4}+S_{1,2})}+c_{2,4}e^{2\pi i\theta _2+{i\over
h}S_{1,2}}-c_{1,3}e^{2\pi i\theta _1+{i\over h}S_{3,4}}-c_{1,4}, }
where we took $u_4=1$.

\par In this relation, we substitute (\ref{1D.41}), (\ref{1D.42}) and get
after multiplication with $e^{i\theta _{1,4}/h}$:
$$ 0=a_{2,3}e^{{i\over
h}(\widehat{S}_{3,4}+\widehat{S}_{1,2})}+a_{2,4}e^{{i\over
h}\widehat{S}_{1,2}}-a_{1,3}e^{{i\over h}\widehat{S}_{3,4}}-a_{1,4}, $$
with
\begin{eqnarray}
\label{1D.49.1} \widehat{S}_{1,2}&=&S_{1,2}+\theta _{1,2}+2\pi
h\theta _2=S_{1,2}+\theta _{1,2}+S_2+hk_2{\pi \over 2}\\
\nonumber
\widehat{S}_{3,4}&=&S_{3,4}+\theta _{3,4}+2\pi h\theta
_1=S_{3,4}+\theta _{3,4}+S_1+hk_1{\pi \over 2},
\end{eqnarray}
where we recall that
$$
\theta _j={S_j\over 2\pi }+{k_j\over 4}.
$$
With \begeq \label{1D.49.2}
\widetilde{S}_{j,k}=\widehat{S}_{j,k}+h{\pi \over 2},
\endeq
we get
\ekv{1D.50}
{
0=a_{2,3}e^{{i\over
h}(\widetilde{S}_{3,4}+\widetilde{S}_{1,2})-i{\pi \over 2}}+a_{2,4}e^{{i\over
h}\widetilde{S}_{1,2}}-a_{1,3}e^{{i\over
h}\widetilde{S}_{3,4}}-a_{1,4}e^{i{\pi \over 2}}.
}

\begin{prop}
\label{Prop1D3} Assume that
$$
Q(\tau,x,\xi,\eps,\frac{h^2}{\eps};h)\sim
Q_0(\tau,x,\xi,\eps,\frac{h^2}{\eps})+h
Q_1(\tau,x,\xi,\eps,\frac{h^2}{\eps})+\ldots
$$
is holomorphic in $(\tau,(x,\xi))\in {\rm neigh}(0,\comp)\times
{\rm neigh}(K_{0,0},\comp^2)$ and depends smoothly on
$\eps,\frac{h^2}{\eps}\in {\rm neigh}(0,\real)$. Here $K_{0,0}$ is
an $\infty$--shaped curve with the self-crossing at $(0,0)$.
Assume furthermore that
$$
Q_0(\tau,x,\xi,0,0)=\langle{q}\rangle(\tau,x,\xi)=f(\tau,\Re
\langle{q}\rangle(\tau,x,\xi)),
$$
where $f$ is an analytic function with $f(0,0)=0$. We assume next
that along $K_{0,0}$, $\Re \langle{q}\rangle(0,x,\xi)=0$ and that
$\Re \langle{q}\rangle(0,x,\xi)$ has a unique critical point on
$K_{0,0}$, $(0,0)$, which is a non-degenerate saddle point. When
$z\in {\rm neigh}(0,\comp)$, put
$$
\mu=f(z;h)
$$
where $f(z;h)$ is an analytic symbol introduced in Propositions
{\rm 6.1} and {\rm 6.2}, and in {\rm (\ref{1D.34})}. Then $z$ is a
quasi-eigenvalue of the operator
$Q\left(\tau,x,hD_x,\eps,\frac{h^2}{\eps};h\right)$  acting on
$L^2_{\theta'}(\real)$ if and only if the corresponding $\mu$
satisfies {\rm (\ref{1D.50})}. In {\rm (\ref{1D.50})}, the
coefficients $a_{1,3}$, $a_{1,4}$, $a_{2,3}$, and $a_{2,4}$ are
introduced in Proposition {\rm 4.1}, and the quantities
$\widetilde{S}_{1,2}$ and $\widetilde{S}_{3,4}$ are defined in
{\rm (\ref{1D.49.1})} and {\rm (\ref{1D.49.2})}. They depend
\hol{}ally on $\mu$ with \ekv{1D.51} {
\partial _\mu \widetilde{S}_{j,k}={\cal O}(1),
} and when $\mu $ is real, we have
 \ekv{1D.52} { \Im
\widetilde{S}_{j,k}={\cal O}\left(\epsilon +{h^2\over \epsilon
}\right). }
\end{prop}

In the formulation of the proposition, we leave the notion of a
quasi-eigenvalue undefined and refer the reader to Section 11 for
a complete justification of this terminology.

\section{Zeros of sums of exponential \fu{}s}\label{SectionZe}
\setcounter{equation}{0}
\par

Here we elaborate on arguments in \cite{Da}, and a related and
even more general discussion can be found in Hager~\cite{Ha}. The
results established in this section will be used in Section 10.

Let $\gamma _1,\gamma _2,..., \gamma _N$ be compact $C^1$ segments
in ${\bf C}$ such that $\gamma _j$ starts at $s_{j-1}\in
D(z_{j-1},r_{j-1}/2)$ and ends at $e_j\in D(z_j,r_j/2)$, where we
use the cyclic convention and view the index $j$ as an element of
${\bf Z}/N{\bf Z}$. We assume that $N\in\{ 1,2,...\}$ is fixed,
but allow $\gamma _j,z_j,r_j,s_j,e_j$ to vary with the
semi-classical parameter $h$ while all estimates below will be
uniform in $h$. Let $f$ be a \hol{} \fu{} defined in
$\cup_{j=0}^{N-1}(D(z_j,r_j)\cup{\rm neigh\,}(\gamma _{j+1}))$,
such that
$$f=e^{{i\over h}S_j(z)+{\cal O}(1)}\hbox{ on }\gamma _j,$$
$$\vert f\vert \le e^{{1\over h}(-\Im S_j(z))+{\cal O}(1)}\hbox{ on
}D(z_j,r_j),$$
where $S_j$ is \hol{} in ${\rm neigh\,}(\gamma _j)\cup
D(z_{j-1},r_{j-1})\cup D(z_j,r_j)$ and
$$\Im (S_{j+1}-S_j)={\cal O}(h)\hbox{ on }D(z_j,r_j).$$

\par In $D(z_j,r_j)$ we can write
$$f(z)=e^{{i\over h}S_j(z)}g_j(z),\ \vert g_j(z)\vert \le {\cal O}(1).$$
We further know that $\vert g_j(e_j)\vert \ge 1/{\cal O}(1)$.
Standard arguments (see for instance \cite{Sj7}), including
Jensen's formula, imply that the number of zeros of $g_j$ in
$D(z_j,r_j/2)$ is ${\cal O}(1)$ and if $\alpha _j$ is a segment in
$D(z_j,r_j/2)$ from $e_j$ to $s_j$ which avoids the zeros
$w_1,...,w_M$ of $g_j$ in $D(z_j,r_j/2)$ such that $\vert {\rm
var\, arg}_{\alpha _j}(z-w_k)\vert <2\pi $ for every $k$, then
\ekv{ze.1} {\Re {1\over 2\pi i}\int_{\alpha _j}{g_j'\over
g_j}dz={\cal O}(1),} and consequently \ekv{ze.2} {\Re {1\over 2\pi
i}\int_{\alpha _j} {f'\over f}dz={\cal O}(1)+{1\over 2\pi
h}\int_{\alpha _j}\Re S_j'(z)dz={\cal O}(1)+{1\over 2\pi h}\Re
(S_j(s_j)-S_j(e_j)).}

\par Let $\gamma $ be the closed contour given by $\gamma _1\cup \alpha
_1\cup \gamma _2\cup ... \cup \gamma _N\cup \alpha _0$. We want to study
$$N(f,\gamma ):={1\over 2\pi i}\int_\gamma  {f'(z)\over f(z)}dz= {1\over
2\pi }{\rm var\, arg}_\gamma (f).$$ When $\gamma $ is the oriented
\bdy{} of a \bdd{} domain $\Gamma $, where $f$ is \hol{}, then
$N(f,\gamma )$ is the number of zeros of $f$ inside $\Gamma $.

\par Along $\gamma _j$, we write $f=e^{{i\over h}\widetilde{S}_j(z)}$,
$\widetilde{S}_j(z)=S_j(z)+{\cal O}(h)$. Then,
$${1\over 2\pi i}\int_{\gamma _j}{f'\over f}dz={1\over 2\pi h}\int_{\gamma
_j}\widetilde{S}_j'(z)dz={1\over 2\pi
h}(\widetilde{S}_j(e_j)-\widetilde{S}_j(s_{j-1})),$$
so
$$\sum_{j=1}^N{1\over 2\pi i}\int_{\gamma _j}{f'\over f}dz={1\over 2\pi
h}\sum_{j=0}^{N-1}(\widetilde{S}_j(e_j)-\widetilde{S}_{j+1}(s_j))={1\over
2\pi h}\sum_{j=0}^{N-1}(S_j(e_j)-S_{j+1}(s_j))+{\cal O}(1),$$
and hence in view of (\ref{ze.2}) and the uniform \bdd{}ness of $N$:
\ekv{ze.3}
{
N(f,\gamma )={1\over 2\pi h}\sum_{j=0}^{N-1}(S_j(s_j)-S_{j+1}(s_j))+{\cal
O}(1).
}
Here we recall that $\Im (S_j-S_{j+1})={\cal O}(h)$ in $D(z_j,r_j)$. It
follows that $\nabla (S_j-S_{j+1})={\cal O}(h/r_j)$ in $D(z_j,r_j/2)$ and
consequently that
$$S_j(s_j)-S_{j+1}(s_{j+1})=S_j(z)-S_{j+1}(z)+{\cal O}(h),$$
for any other point $z\in D(z_j,r_j/2)$. Thus finally,
\ekv{ze.4}
{N(f,\gamma )={1\over 2\pi h}\sum_{j=0}^{N-1}(S_j(w_j)-S_{j+1}(w_j))+{\cal
O}(1),}
with $w_j\in D(z_j,r_j/2)$ chosen \ably{}. Here we can further replace
$S_j(w_j)-S_{j+1}(w_j)$ by its real part, since $\Im S_j-S_{j+1}={\cal O}(h)$
in $D(z_j,r_j)$.

\section{Skeleton in the region $\vert \mu \vert \gg h$.}\label{SectionSk}
\setcounter{equation}{0}

\par We now return to the situation in Section \ref{Section1D}. We are
interested in the solutions $\mu $ of (\ref{1D.50}). In the
following, we will write $S_{j,k}$ instead of
$\widetilde{S}_{j,k}$, so we are interested in the zeros of the
\fu{} $F_0(\mu ;h)$ appearing in (\ref{1D.50}), given by
$$F_0(\mu ;h)=e^{{i\over h}(S_{1,2}+S_{3,4})-i{\pi \over
2}}a_{2,3}+e^{{i\over h}S_{1,2}}a_{2,4}-a_{1,3}e^{{i\over
h}S_{3,4}}-a_{1,4}e^{{i\pi \over 2}}.$$ Pulling out a factor
$e^{-i\pi /2}$, we get the new equivalent function \ekv{Sk.1} {
F(\mu ;h)=e^{{i\over h}(S_{1,2}+S_{3,4})}a_{2,3}+e^{{i\over
h}S_{1,2}+\pi {\mu \over h}}+e^{{i\over h}S_{3,4}+\pi {\mu \over
h}}+a_{1,4},} which has the same zeros as $F_0$. Here we have also
used the explicit formulae for $a_{2,4}$, $a_{1,3}$ in
(\ref{tr.12}).

\par Using the results of Section 5, we shall now look at the \asy{}s of $F(\mu ;h)$, when $\vert
\mu \vert /h\gg 1$.  \smallskip

\noindent {\it Case 1:} Assume that \ekv{Sk.6} { Ch\le \vert \mu
\vert \ll 1,\ \abs{{\rm arg\,}\mu  -{\pi \over 2}} \le \pi
-{1\over C}.} (The case 2, given by $\vert {\rm arg\,}\mu +{\pi
\over 2}\vert <{\pi }-{1\over {\cal O}(1)}$ will be reduced to the
case 1 by a symmetry argument.)  In this region, we have
(\ref{as.8}):
$$a_{1,4}=
e^{{\cal O}_-({h\over\mu })}\left(e^{{i\over h}(-\mu \ln {\mu
\over i} +\mu -{\pi h\over 4}+i{\pi \mu \over 2}) }+e^{{i\over
h}(-\mu \ln {\mu \over i} +\mu -{\pi h\over 4}-i{3\pi \mu \over
2}) }\right),$$ and we get using also (\ref{Sk.1}) and
(\ref{as.5}),
\begin{eqnarray*}
F(\mu ;h)&=&e^{{i\over h}(S_{1,2}+S_{3,4})}
e^{{i\over h}(\mu \ln {\mu \over i}-
\mu +{\pi h\over 4}
-i{\pi \mu \over 2})-{\cal O}_-({h\over \mu })}+e^{{i\over h}S_{1,2}+\pi
{\mu \over h}} +e^{{i\over h}S_{3,4}+\pi {\mu \over h}}\\ &&+ e^{{\cal
O}_-({h\over \mu })+{i\over h}(-\mu \ln {\mu
\over i}+\mu -{\pi h\over 4}+i{\pi \mu \over 2})} +e^{{\cal O}_-({h\over
\mu })+{i\over h}(-\mu \ln {\mu
\over i}+\mu -{\pi h\over 4}-i{3\pi \mu \over 2})}\\ &=& e^{\pi \mu \over
2h}G(\mu ;h), \end{eqnarray*}
where
\begin{eqnarray}\label{Sk.1.5} G(\mu ;h)&=&a_1+a_2+a_3+a_4,\ a_4=a_{4^+}+a_{4^-}\\
a_1&=&e^{{i\over
h}(S_{1,2}+S_{3,4}+\mu \ln {\mu \over i}-\mu +{\pi h\over 4})-{\cal
O}_-({h\over \mu })}\nonumber\\
a_2&=&e^{{i\over h}S_{1,2}+{\pi \mu \over
2h}},\nonumber\\
a_{3}&=&e^{{i\over h}S_{3,4}+{\pi \mu \over 2h}}\nonumber\\
a_{4^\pm}&=&e^{{\cal O}_-({h\over \mu })+
{i\over h}(-\mu \ln {\mu \over i}+\mu -{\pi h\over 4})\pm{\pi \mu \over
h}}.\nonumber
\end{eqnarray}

\par We have
\ekv{Sk.2}{\vert a_j\vert =e^{r_j/h},\,\, j=1,2,3,4^\pm ,}
where
\begin{eqnarray*} r_1:&=&-\Im S_{1,2}-\Im S_{3,4}+(\Im \mu
)\ln {1\over \vert \mu \vert }-\Re \mu \,{\rm arg\,}{\mu \over i}+
\Im \mu  -h\Re {\cal
O}_-({h\over \mu }),\cr r_2:&=&-\Im
S_{1,2}+{\pi \over 2}\Re \mu ,\cr r_3:&=&-\Im S_{3,4}+{\pi \over
2}\Re \mu ,\cr r_4^\pm:&=& -(\Im \mu )\ln {1\over \vert \mu
\vert }+\Re \mu \,{\rm arg\,}{\mu \over i}-\Im \mu \pm \pi \Re \mu
+h\Re {\cal O}_-({h\over \mu }). \end{eqnarray*}
Notice that $a_{4^\pm}$
is dominating over $a_{4^\mp}$  when $\pm \Re \mu \ge 0$, and in each
half-plane $\pm \Re \mu
>0$, we may associate $a_4$ to the dominating term, modulo an error which
is ${\cal O}(e^{-2\pi \vert \Re \mu \vert /h})$ times the leading term. Also notice that the last
equations take the form
\begin{eqnarray}\label{Sk.3}
r_1&=& (\Im \mu )\ln {1\over \vert \mu \vert }-\Im S_{1,2}-\Im
S_{3,4}-Y(\mu ),\\
r_2&=&-\Im S_{1,2}+{\pi \over 2}\Re \mu ,\nonumber\\
r_3&=&-\Im S_{3,4}+{\pi \over 2}\Re \mu ,\nonumber\\
r_{4^\pm}&=&-(\Im \mu )\ln {1\over \vert \mu \vert }\pm \pi \Re \mu
+Y(\mu ),\nonumber
\end{eqnarray}
where
\ekv{Sk.3.1}
{Y(\mu )=(\Re \mu ){\rm arg\,}({\mu \over i})-\Im \mu +h\Re
{\cal O}_-({h\over \mu }).}

\medskip
Following the general principles, as explained for example
in~\cite{Da} (see also~\cite{BKW}), we shall now look for the
curves $\Gamma_{j,k}$, $j,k=1,2,3,4^{\pm}$, where
$\abs{a_j}=\abs{a_k}$, and we shall especially be interested in
those parts of $\Gamma_{j,k}$ where $\abs{a_j}=\abs{a_k}$ is
dominating over the other $\abs{a_{\nu}}$. In doing so, let us
remark first that we will not see any zeros of $G$ generated by
the zeros of $a_2+a_3$ as a dominating part of $G$, for when $\Re
\mu>0$, then $r_4^+$ dominates over $r_4^-$ and
$r_1+r_4^+=r_2+r_3$ and clearly we cannot have $r_2=r_3\ge \max
(r_1,r_4^+)+{\rm Const}$. Now when $\Re \mu <0$, $r_4^-$ dominates
over $r_4^+$ and
$$r_2+r_3 -\pi \Re \mu =r_1+r_4^-,$$ so that $r_2+r_3 < r_1+r_4^-$ in
this case, leading to an even stronger conclusion.

\medskip We now begin look at the location of zeros of $a_4$ and of sums
of two of the $a_j$.  \smallskip

\par\noindent Zeros of $a_4$: They are of
the form $\mu =i(k+{1\over 2})h$, $k=0,1,2,...$.  \smallskip

\par\noindent Zeros of
$a_3+a_{4^\pm}$: They are contained in the region $\Gamma _{3,4^\pm}$:
\ekv{E_{4+3=0}} { (\Im \mu )\ln {1\over \vert \mu \vert } = \Im
S_{3,4}(\mu )+Y(\mu )-{\pi \over 2}\Re \mu \pm \pi \Re \mu . }

\par Similarly the
zeros of $a_2+a_{4^\pm}$ are contained in $\Gamma _{2,4^\pm}$:
\ekv{E_{4+2=0}} { (\Im \mu )\ln
{1\over \vert \mu \vert } = \Im S_{1,2}(\mu )+Y(\mu )-{\pi \over 2}\Re \mu
\pm \pi \Re \mu . } \smallskip

\par\noindent The zeros of $a_1+a_3$ are
contained in $\Gamma _{1,3}$:
\ekv{E_{1+3=0}} { (\Im \mu )\ln {1\over \vert \mu \vert
}=\Im S_{1,2}(\mu )+ {\pi \over 2}\Re \mu +Y(\mu ). } \smallskip

\par\noindent The zeros of $a_1+a_2$ are
contained in $\Gamma _{1,2}$:
\ekv{E_{1+2=0}} { (\Im \mu )\ln {1\over \vert \mu \vert
}=\Im S_{3,4}(\mu )+ {\pi \over 2}\Re \mu +Y(\mu ).  } \smallskip

\par\noindent The zeros of $a_1+a_{4^\pm}$ are
contained in $\Gamma _{1,4^\pm}$:
\ekv{E_{1+4=0}} {(\Im \mu )\ln {1\over \vert \mu \vert
}={1\over 2}(\Im S_{1,2}+\Im S_{3,4}) \pm {\pi \over 2}\Re \mu +Y(\mu ).}
\smallskip

\par Put
$$X={\pi \over 2}\Re \mu +(\Re \mu ){\rm arg\,}({\mu\over i})-\Im \mu +h\Re {\cal
O}_-({h\over \mu} )={\pi \over 2}\Re \mu +Y(\mu ).$$

\par When $\Re \mu >0$, $a_{4^+}$ dominates over $a_{4^-}$ and we shall
only consider $\Gamma _{3,4^+}=\Gamma _{1,2}$, $\Gamma _{2,4^+}=\Gamma
_{1,3}$, $\Gamma _{1,4^+}$ given by
\ekv{Sk.4}{(\Im \mu )\ln {1\over \vert \mu \vert }=\cases{\Im
S_{3,4}+X,\hbox{ on }\Gamma
_{3,4^+}=\Gamma _{1,2}\cr \Im S_{1,2}+X\hbox{ on }\Gamma _{2,4^+}=\Gamma _{1,3},
\cr {1\over 2}(\Im S_{1,2}+\Im
S_{3,4})+X\hbox{ on
}\Gamma _{1,4^+}.}}

\par Recall now from Section 5 that ${\cal O}_-(h/\mu )$ in (\ref{Sk.3.1})
appears as a remainder in Stirling's formula, so that $\partial
_\mu {\cal O}_-(h/\mu )={\cal O}(h/\mu ^2)$, and hence $X$ is
\ufly{} Lipschitz for $\vert \mu \vert \ge h$. Proposition
\ref{PropA3.1} can therefore be applied to get the approximate
behavior of the $\Gamma _{j,k}$. This will be exploited later.

\par In the left half-plane, $a_{4^-}$ dominates over $a_{4^+}$ and we consider
all the 5 curves $\Gamma _{3,4^-}$, $\Gamma _{2,4^-}$, $\Gamma
_{1,4^-}$, $\Gamma _{1,2}$, and $\Gamma _{1,3}$, given by:
\ekv{Sk.5}{(\Im \mu )\ln {1\over \vert \mu \vert }=\cases{\Im
S_{3,4}-2\pi \Re \mu +X,\hbox{ on }\Gamma _{3,4^-}\cr \Im
S_{1,2}-2\pi \Re \mu +X, \hbox{ on }\Gamma _{2,4^-}\cr \Im
S_{1,2}+X\hbox{ on }\Gamma _{1,3},\cr \Im S_{3,4}+X\hbox{ on
}\Gamma _{1,2},\cr {1\over 2}(\Im S_{1,2}+\Im S_{3,4})-\pi \Re \mu
+X\hbox{ on }\Gamma _{1,4^-}.}} Again Proposition \ref{PropA3.1}
can be applied to give the approximate shape of $\Gamma _{j,k}$.
Recall that we are in the case 1 with $\vert \mu \vert \gg h$, so
that (\ref{Sk.6}) holds.

\smallskip
\par\noindent \it 1. Skeleton in the region $\Re \mu \ge 0$. \rm (We will
implicitly use that $\vert a_2\vert \vert a_3\vert =\vert a_1\vert
\vert a_{4^+}\vert $.) The region $\vert \Re \mu \vert \le {\cal
O}(h)$ will require a special discussion. In the region $\Re \mu
\ge Ch$, we have $\vert a_4^+\vert =e^{2\pi \Re \mu /h}\vert
a_{4^-}|\gg \vert a_4^-\vert $, so $\vert a_4\vert \sim \vert
a_4^+\vert $ and in this region we see from an earlier observation
that the zeros of $a_2+a_3$ will not play any essential role. In
this region we shall therefore use the $\Gamma _{j,k}$ appearing
in (\ref{Sk.4}) and, as pointed out earlier, we are interested
here in the part of each such $\Gamma _{j,k}$, where $\vert
a_j\vert =\vert a_k\vert $ dominates over the other $\vert a_\nu
\vert $.

It follows from Proposition \ref{PropA3.1} that the curves in
(\ref{Sk.4}) (as well as the ones in (\ref{Sk.5})), are of the
form \ekv{Sk.8} {\Im \mu =\gamma _{j,k}(\Re \mu ),\hbox{ with
}\vert \gamma _{j,k}\vert ,\vert \gamma '_{j,k}\vert \ll 1.}
Notice that every crossing point of two of the curves $\Gamma
_{1,4^+},\Gamma _{1,2},\Gamma _{2,4^+}$ is also a crossing point
for all three. This follows from (\ref{Sk.4}) or even more
trivially from the observation that 2 of the 3 \e{}s $\vert
a_1\vert =\vert a_{4^+}\vert $, $\vert a_1\vert =\vert a_2\vert $,
$\vert a_2\vert =\vert a_{4^+}\vert $, imply the third one. Also
notice that if we draw the two curves $\Gamma _{1,2}=\Gamma
_{3,4^+}$, $\Gamma _{1,3}=\Gamma _{2,4^+}$, then $\Gamma _{1,4^+}$
is between the two---see Figure 2 on page 35.

\par For $\mu >0$ we have $X(\mu )=h\Re {\cal O}_-({h\over \mu })=h{\cal
O}(e^{-2\pi \mu /h})$ by (\ref{as.4.1}) and hence $\gamma _{j,k}(\Re \mu )$
is described as in Proposition \ref{PropA3.1} with
\ekv{Sk.8.5}
{
F=F_{j,k}(\Re \mu )=h{\cal O}(e^{-2\pi \Re \mu /h})+\cases{\Im S_{3,4}(\Re
\mu ),\ (j,k)=(3,4^+),\, (1,2),\cr \Im S_{1,2}(\Re \mu ),\ (j,k)=(2,4^+),\,
(1,3),\cr {1\over 2}(\Im S_{1,2}+\Im S_{3,4})(\Re \mu ),\ (j,k)=(1,4^+).}
}

\par In the region $\Im \mu <\min (\gamma _{2,4^+}(\Re \mu ),\gamma
_{3,4^+}(\Re \mu ))$ we have $\vert a_{4^+}\vert \ge \max (\vert
a_1\vert ,\vert a_2\vert ,\vert a_3\vert ,\vert a_{4^-}\vert )$,
and if we restrict further to \ekv{Sk.9} {\Im \mu <\min (\gamma
_{2,4^+}(\Re \mu ), \gamma _{3,4^+}(\Re \mu ))-{Ch\over \ln
{1\over \vert \mu \vert }},\ \Re \mu >Ch,} with $C\gg 1$, we see
that $a_{4^+}$ is dominating in the sense that \ekv{Sk.10} { \vert
a_{4^+}\vert \ge 2\vert a_1+a_2+a_3+a_{4^-}\vert , } and hence
$G(\mu ;h)$ has no zeros in that region. Similarly in the region
$\Im \mu
>\max (\gamma _{1,2},\gamma _{1,3})(\Re \mu )$, we have $\vert
a_1\vert \ge \vert a_2\vert ,\vert a_3\vert, \vert
a_{4^\pm}\vert$, and if \ekv{Sk.11} { \Im \mu >\max (\gamma
_{1,2}, \gamma _{1,3})(\Re \mu )+{Ch\over \ln {1\over \vert \mu
\vert }},\ \Re \mu \ge 0, } with $C\gg 1$, then $a_1$ is
dominating in the sense that \ekv{Sk.12} { \vert a_1\vert \ge
2\vert a_2+a_3+a_4\vert , } and again $G(\mu ;h)$ has no zeros
there.

\par Now consider a point $\mu \in \Gamma _{3,4^+}$, where $\gamma
_{3,4^+}\le \gamma _{2,4^+}$, so that $\Im S_{3,4}\le \Im
S_{1,2}$, with $\Re \mu \gg h$.  Going down (ie.~decreasing $\Im
\mu $ while keeping $\Re \mu $ constant) by a distance $\gg h/\ln
{1\over \vert \mu \vert }$, we reach the region (\ref{Sk.9}),
where $a_{4^+}$ dominates.
\smallskip
\par\noindent \it Case a. \rm $0\le \Im S_{1,2}-\Im S_{3,4}\le {\cal
O}(h)$. Going up by a distance $\gg h/\ln {{1\over \vert \mu \vert }}$, we
cross $\Gamma _{2,4^+}=\Gamma _{1,3}$ and reach the region (\ref{Sk.11}), where
$a_1$ is dominating.
\smallskip
\par\noindent \it Case b. \rm $\Im S_{1,2}-\Im S_{3,4}\ge Ch$ for $C\gg
1$. Going up by a distance $\sim h/\ln {1\over \vert \mu \vert }$, we reach
the region, where $a_3$ is dominating:
\ekv{Sk.13}
{\vert a_3\vert \ge 2\vert a_1+a_2+a_4\vert ,}
and continuing to go up, $a_3$
remains dominating until we reach a $h/\ln{1\over \vert \mu \vert
}$-\neigh{} of $\Gamma _{1,3}=\Gamma _{2,4^+}$. After crossing that curve
and going up by another amount $Ch/\ln {1\over \vert \mu \vert }$, we
reach the region, where $a_1$ is dominating.

\begin{figure}
\centering
    \psfrag{1}[l][l]{$\Im \mu$}
    \psfrag{2}[l][l]{$\Re \mu$}
    \psfrag{3}[1][1]{$Ch$}
    \psfrag{4}[1][1]{$\Gamma_{3,4^+}=\Gamma_{1,2}$}
    \psfrag{5}[1][1]{$\Gamma_{2,4^+}=\Gamma_{1,3}$}
    \psfrag{6}[1][1]{$\Gamma_{1,4^+}$}
    \psfrag{7}[1][1]{$a_{4^+}$ dominates}
    \psfrag{8}[1][1]{$a_1$ dominates}
    \psfrag{9}[1][1]{$a_2$ dominates}
    \psfrag{10}[1][1]{$a_3$ dominates}
\scalebox{0.95} {\includegraphics{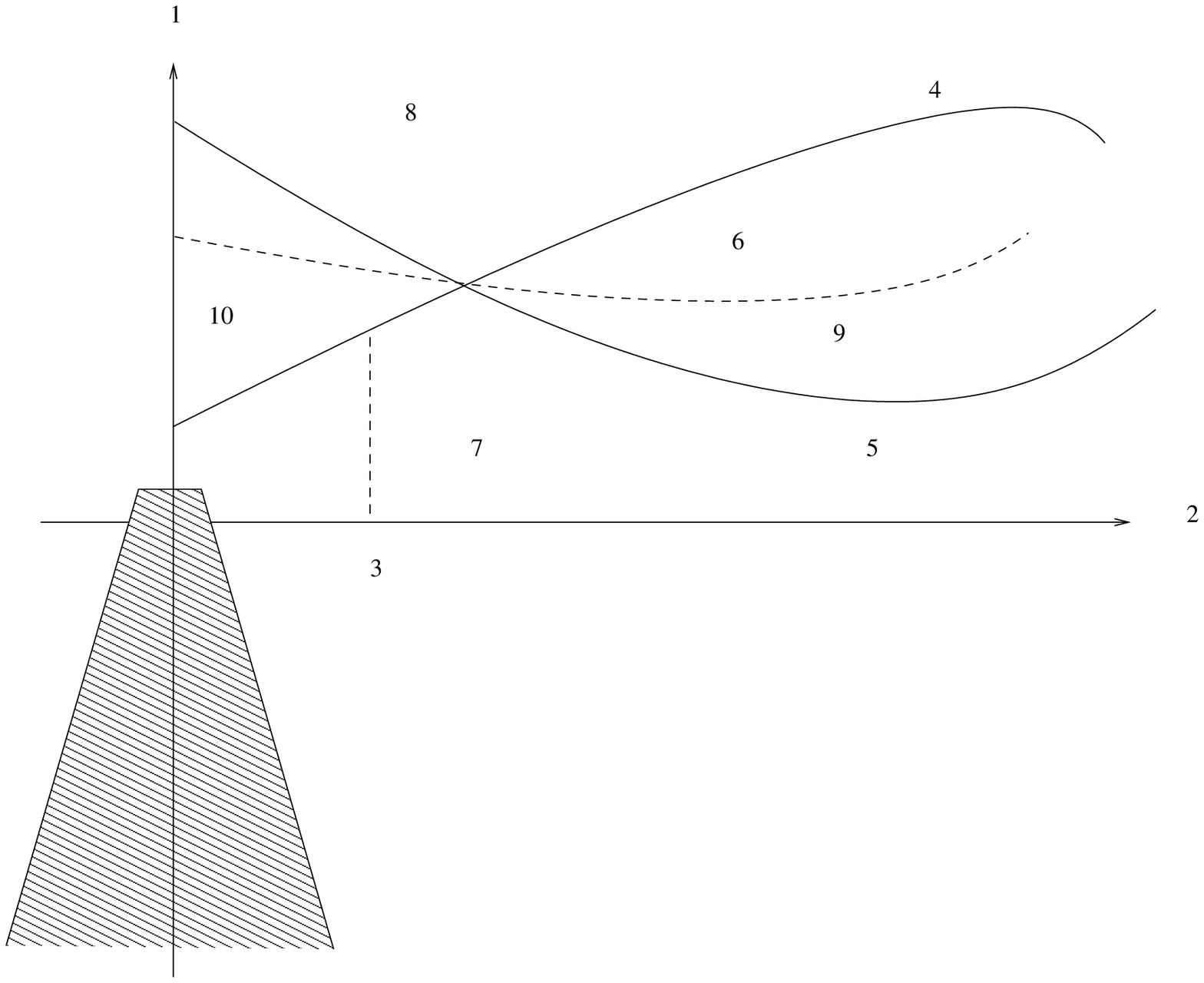}} \caption{The
union of the solid curves in the figure gives a schematic
representation of the skeleton $S'$ in the right half-plane
intersected with the region (\ref{Sk.6}). Proposition 8.1 shows
that the zeros of $G$ in this region are inside the union of the
thickened skeleton, obtained by placing a disc of radius
$Ch/\abs{\ln \abs{\mu}}$ around each point $\mu\in S'$, and the
set of all $\mu$ with $0\leq \Re \mu \leq Ch$ below $S'$.
Proposition 8.3 gives a more precise description of the location
of the zeros of $G$ with $\abs{\Re \mu}\leq {\cal O}(h)$.}
\end{figure}
\par Our discussion shows the following
\begin{prop}\label{PropSk.1}
We work in the region {\rm (\ref{Sk.6})} and assume in addition
that $\Re\mu >0$.
\par\noindent If $\Im \mu \le \min (\gamma _{2,4^+},\gamma _{3,4^+})(\Re
\mu )$, then
\ekv{Sk.14}{\vert a_{4^+}\vert \ge \max (\vert a_1\vert , \vert a_2\vert ,
\vert a_3\vert , \vert a_{4^-}\vert ),}
\smallskip
\par\noindent  If $\Im \mu \ge \max (\gamma _{1,2},\gamma _{1,3})({\Re
\mu )}$, then
\ekv{Sk.15}
{
\vert a_1\vert \ge \max (\vert a_2\vert , \vert a_3\vert , \vert
a_{4^\pm}\vert ).
}
\smallskip
\par\noindent If $\gamma _{3,4^+}(\Re \mu )(=\gamma _{1,2}(\Re \mu ))\le
\Im \mu  \le \gamma _{1,3}(\Re \mu )(=\gamma _{2,4^+}(\Re \mu ))$, then
\ekv{Sk.16}
{\vert a_3\vert \ge \max (\vert a_1\vert ,\vert a_2\vert ,\vert
a_{4^\pm}\vert ).}
\smallskip
\par\noindent If $\gamma _{2,4^+}(\Re \mu )(=\gamma _{1,3}(\Re \mu ))\le
\Im \mu \le \gamma _{3,4^+}(\Re \mu )(=\gamma _{1,2}(\Re \mu ))$, then
\ekv{Sk.17}
{\vert a_2\vert \ge \max (\vert a_1\vert ,\vert a_3\vert ,\vert
a_{4^\pm}\vert ).}
\smallskip

\par If the distance from $\mu $ to $(\Gamma _{1,2}=\Gamma
_{3,4^+})\cup (\Gamma _{1,2}=\Gamma _{3,4^+})$ is $\ge Ch/\ln
{1\over \vert \mu \vert }$, with $C\gg 1$, then in the respective
cases {\rm (\ref{Sk.15})}--{\rm (\ref{Sk.17})} can be sharpened to
the dominance in the sense explained above. In particular, $G$ has
no zeros in this region. If, in addition, $\Re \mu \ge Ch$ with
$C\gg 1$, then we have the same conclusion in the case of {\rm
(\ref{Sk.14})}.
\end{prop}

\par In the region (\ref{Sk.6}), intersected with the right half-plane
$\Re \mu >0$ we define the skeleton to be the union of the curves
$\Im \mu =\max (\gamma _{1,2},\gamma _{1,3})(\Re \mu )$ and $\Im
\mu =\min (\gamma _{2,4^+},\gamma _{3,4^+})(\Re \mu )$. The
proposition shows that the zeros of $G$ in the region under
consideration are contained in the union of all discs $D(\mu
,Ch/\ln {1\over \vert \mu \vert })$ with $\mu $ in the skeleton
just defined, and the set of all $\mu$ below the skeleton, with
$0\le \Re \mu <Ch$, for $C\gg 1$.\smallskip

\par\noindent \it 2. Skeleton in the region $\Re \mu \le 0$. \rm  Again the
region $\vert \Re \mu \vert \le {\cal O}(h)$ will require a separate
discussion so we restrict the attention to $\Re \mu \le -Ch$ and we will
use $\vert a_{4^-}\vert =e^{-2\pi \mu /h}\vert a_{4^+}\vert \gg \vert
a_{4^+}\vert $, so that $\vert a_4\vert \sim \vert a_{4^-}\vert $. We
therefore concentrate the attention to the curves in (\ref{Sk.5}). As
before we notice that every crossing point of two of the three curves
$\Gamma _{1,2},\Gamma _{2,4^-},\Gamma _{1,4^-}$ is a crossing point of all
three.  The same holds for $\Gamma _{1,3},\Gamma _{3,4^-},\Gamma _{1,4^-}$.

\par Now use that $\vert \Im \mu \vert $ is small and hence that $\Im
S_{1,2},\Im S_{3,4}$ and their derivatives \wrt{} $\Re \mu $
are small. We can therefore consider the two curves:
$$A:\, -2\pi \Re \mu =\Im S_{1,2}-\Im S_{3,4},\
B:\, -2\pi \Re \mu =\Im S_{3,4}-\Im S_{1,2}.$$
They are of the form
$$-\Re \mu =\gamma _A(\Im \mu ),\ -\Re \mu =\gamma _B(\Im \mu ),$$
where $\gamma _A$, $\gamma_B$
are small with small derivatives and satisfy:
\ekv{Sk.18}
{\vert \gamma _A(\Re \mu )\vert \sim \vert \gamma _B(\Re \mu )\vert ,\
\gamma _A\gamma _B\le 0.}

\par The curve $\Gamma _{1,4^-}$ will play a central role. It crosses $A,B$
at the points $\mu _A,\mu _B$ (unless these points are hidden in the
forbidden region), and we have
\ekv{Sk.19}
{
(\Re \mu _A)(\Re \mu _B)\le 0,\ \vert \Re \mu _A\vert \sim \vert \Re \mu
_B\vert .
}
We notice that $\mu _A$
is the unique crossing point for $\Gamma _{1,3},\Gamma _{3,4^-},\Gamma
_{1,4^-}$ while $\mu _B$
is the unique crossing point for $\Gamma _{1,2},\Gamma _{2,4^-},\Gamma
_{1,4^-}$. More precisely,  $\gamma _{1,3}(t)-\gamma _{1,4^-}(t)$,
$\gamma _{1,4^-}(t)-\gamma _{3,4^-}(t)$ vanish precisely for $t=\Re \mu _A$ and have the
same sign as $t-\Re \mu _A$. Similarly $\gamma _{1,2}(t)-\gamma _{1,4^-}(t)$,
$\gamma _{1,4^-}(t)-\gamma _{2,4^-}(t)$ vanish precisely for $t=\Re \mu
_B$ and have the same sign as $t-\Re \mu _B$. We also notice that if $\mu$
belongs to one of the three curves $\Gamma _{1,3},\Gamma _{1,4^-},\Gamma
_{3,4^-}$, then the distance from $\mu $ to any of the two other curves
among these three is $\ge C^{-1}\vert \Re \mu -\Re \mu _A\vert /\ln
{1\over \vert \mu \vert }$. The same observation holds for $\Gamma _{1,2},\Gamma _{1,4^-},
\Gamma _{2,4^-}$ with $\mu _B$ instead of $\mu _A$.

\par For $\mu <0$, we have $X-\pi \Re \mu =h\Re {\cal O}_-({h\over
\mu })=h{\cal O}(e^{-2\pi \vert \mu \vert /h})$ by (\ref{as.4.1}),
and hence $\gamma _{1,4^-}$ is described as in Proposition
\ref{PropA3.1} with \ekv{Sk.19.5} {F=F_{1,4^-}(\Re \mu )={1\over
2}(\Im S_{1,2}+\Im S_{3,4})(\Re \mu )+h{\cal O}(e^{-2\pi \vert \Re
\mu \vert /h}).}

\par Now assume to fix the ideas that $\Re \mu _A\le 0$ (the case $\Re \mu
_B\le 0$ can be treated similarly). Considering the three curves $\Gamma _{1,3},\Gamma _{1,4^-},\Gamma
_{3,4^-}$, we see that for $\Re \mu \le 0$:
\smallskip
\par\noindent If $\Im \mu \le \min (\gamma _{1,4^-},\gamma _{3,4^-})(\Re
\mu )$, then
$$\vert a_{4^-}\vert \ge \vert a_1\vert ,\vert a_2\vert ,\vert a_3\vert ,\vert
a_{4^+}\vert .$$
(In this case, we also have $\Im \mu\le \gamma _{2,4^-}(\Re \mu )$.)
\smallskip
\par\noindent If
$\Im \mu  \ge \max (\gamma _{1,4^-},\gamma _{1,3})(\Re \mu )$,
then
$$\vert a_1\vert \ge \vert a_2\vert ,\vert a_3\vert ,\vert a_{4^\pm}\vert .$$
(In this case, we also have $\Im \mu\ge \gamma _{1,2}(\Re \mu)$.)
\smallskip
\par\noindent If $\gamma _{3,4^-}(\Re \mu )\le \Im \mu \le \gamma
_{1,3}(\Re \mu )$, then
$$\vert a_3\vert \ge \vert a_1\vert ,\vert a_2\vert ,\vert a_{4^\pm}\vert .$$
This covers all possible cases with $\Re \mu \le 0$. Notice that the last
case can appear only when $\Re \mu _A\le \Re \mu \le 0$. When $\Re \mu
_B\le 0$, we get the analogous discussion after a permutation of the
indices 2 and 3:
\par\noindent If $\Im \mu \le \min (\gamma _{1,4^-},\gamma _{2,4^-})(\Re
\mu )$, then
$$\vert a_{4^-}\vert \ge \vert a_1\vert ,\vert a_2\vert ,\vert a_3\vert ,\vert
a_{4^+}\vert .$$
\smallskip
\par\noindent If $\Im \mu \ge \max (\gamma _{1,4^-},\gamma _{1,2})(\Re \mu )$,
then
$$\vert a_1\vert \ge \vert a_2\vert ,\vert a_3\vert ,\vert a_{4^\pm}\vert .$$
\smallskip
\par\noindent If $\gamma _{2,4^-}(\Re \mu )\le \Im \mu \le \gamma
_{1,2}(\Re \mu )$, then
$$\vert a_2\vert \ge \vert a_1\vert ,\vert a_3\vert ,\vert a_{4^\pm}\vert .$$
\smallskip

\par Moreover, if in addition all the inequalities for $\Im \mu $ are valid
with an extra margin $Ch/\ln {1\over \vert \mu \vert }$, $C\gg 1$,
then in the various cases, we have dominance of $a_1,a_2,a_3$
respectively, in the sense explained before. For the dominance of
$a_{4^-}$, we also need the assumption that $\Re \mu \le -Ch$ with
$C\gg 1$.

\begin{figure}
\centering
    \psfrag{1}[l][l]{$\Im \mu$}
    \psfrag{2}[l][l]{$\Re \mu$}
    \psfrag{3}[1][1]{$-Ch$}
    \psfrag{4}[1][1]{$\Gamma_{1,3}$}
    \psfrag{5}[1][1]{$\Gamma_{1,2}$}
    \psfrag{6}[1][1]{$\Gamma_{1,4^-}$}
    \psfrag{7}[1][1]{$\Gamma_{3,4^-}$}
    \psfrag{8}[1][1]{$\Gamma_{2,4^-}$}
    \psfrag{9}[1][1]{$a_{1}$ dominates}
    \psfrag{10}[1][1]{$a_{4^-}$ dominates}
    \psfrag{11}[1][1]{$a_3$ dominates}
\scalebox{1} {\includegraphics{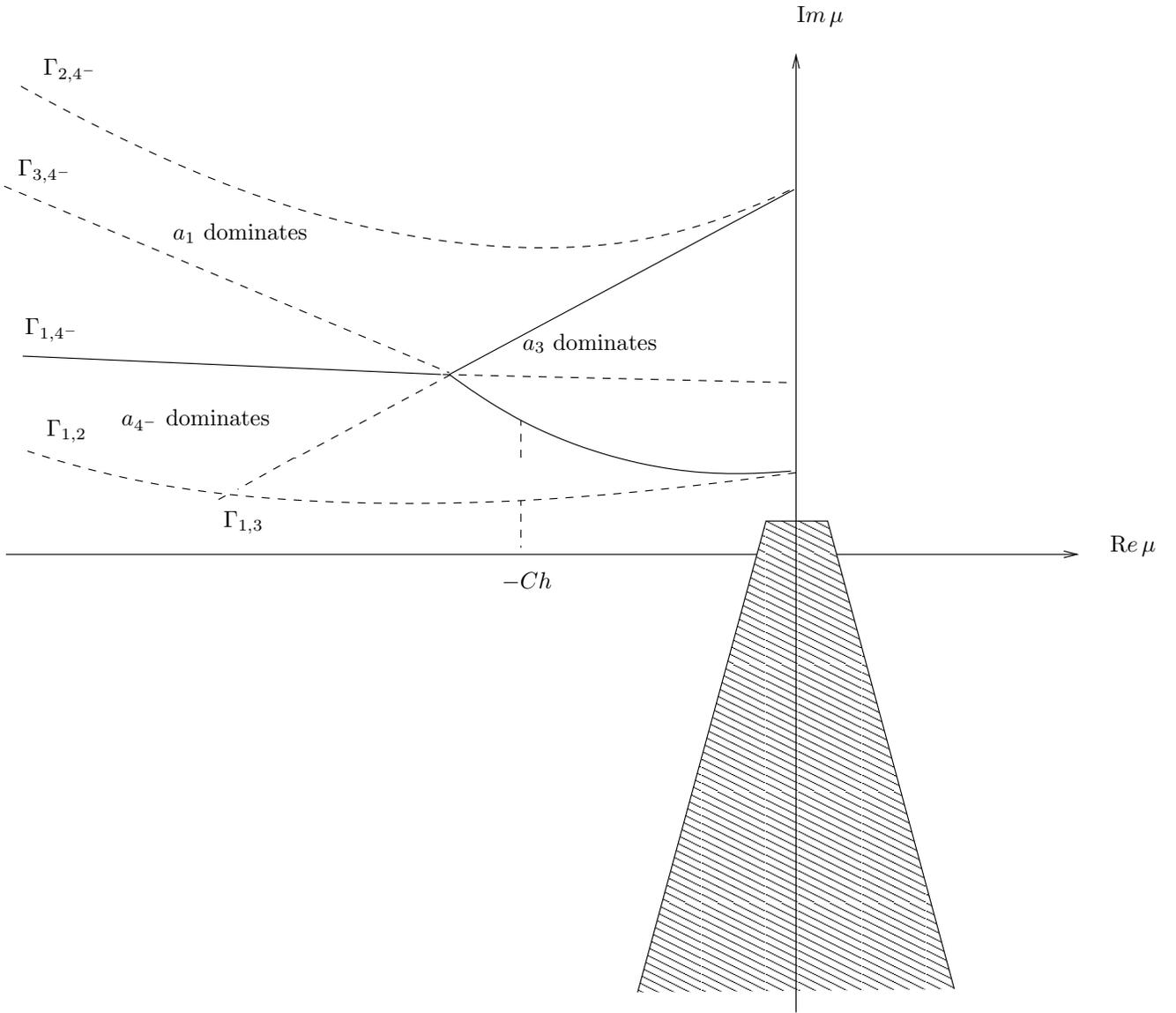}} \caption{The union
of the solid curves in the figure gives a schematic representation
of the skeleton $S'$ in the left half-plane intersected with the
region (\ref{Sk.6}). Proposition 8.2 shows that the zeros of $G$
are inside the union of the thickened skeleton, obtained by
placing a disc of radius $Ch/\abs{\ln \abs{\mu}}$ around each
point $\mu\in S'$, and the set of all $\mu$ with $\abs{\Re \mu}
\leq Ch$ below $S'$.}
\end{figure}
\smallskip
\par\noindent \it Exponential localization to the skeleton. \rm Recall that
we are still working in the region (\ref{Sk.6}). In this region we define
the skeleton to be:
\ekv{Sk.20}{S=S'\cup \Gamma _4,}
where we define $S'$ to be the union of the following two sets in
respectively the left and the right half-planes: \smallskip
\par\noindent -- In the closed left half-plane (intersected with
(\ref{Sk.6})), assume $\Re \mu _A\le 0$ to fix the ideas, then this part
of $S'$ is given by all points of the form $\Im \mu =\gamma _{1,4^-}(\Re
\mu )$ with $\Re \mu \le \Re \mu _A$, all points of the form $\Im \mu
=\gamma _{3,4^-}(\Re \mu )$ or of the form $\Im \mu =\gamma _{1,3}(\Re
\mu )$, with $\Re \mu _A\le \Re \mu \le 0$.
\smallskip
\par\noindent -- In the closed right half-plane the corresponding part
of $S'$ is defined to be the union of the two curves:
$\Im \mu =\max (\gamma _{1,3},\gamma _{1,2})(\Re \mu )$ and
$\Im \mu =\min (\gamma _{2,4^+},\gamma _{3,4^+})(\Re \mu )$.
\smallskip
\par $\Gamma _4$ is defined to be the part of the imaginary axis given by
$$
\Re \mu =0,\ Ch\le \Im \mu \le \min (\gamma _{2,4^+},\gamma
_{3,4^+})(0),
$$
where $C$ is the same constant as in (\ref{Sk.6}). Notice that this
part may be empty.
The earlier discussion shows that we have
\begin{prop}\label{PropSk.2}
The zeros of $G$ in the domain {\rm (\ref{Sk.6})} are contained in
the set \ekv{Sk.21} {\left(\bigcup_{\mu \in S'}D\left(\mu
,{Ch\over \ln {1\over \vert \mu \vert }}\right)\right)\bigcup\{\mu
\hbox{ below }S';\, \vert \Re \mu \vert <Ch\} . }
\end{prop}

\Remark. The localization result of Proposition 8.2 improves if we
are far from the branch points of the skeleton. Thus for instance,
if $\mu $ is below $S'$ and ${\rm dist\,}(\mu ,S')\ge Ch/\ln
{1\over \vert \mu \vert }$, then
$$\max \vert a_{4^\pm}\vert \ge
e^{{\ln {1\over \vert \mu \vert } \over Ch} {\rm dist\,}(\mu ,S')}
(\vert a_1\vert
+\vert a_2\vert +\vert a_3\vert ),$$
so the zeros of $G$
are exponentially small perturbations of those of $a_4=a_{4^+}+a_{4^-}$: In
this region there is a bijection $b$ between the zeros of $a_4$ and those
of $G$, with
$$\vert b(\mu )-\mu \vert \le {\cal O}(h)\exp
[-(Ch)^{-1}\ln({1\over \vert \mu \vert }){\rm dist\,}(\mu ,S')].$$

\par Similarly let $\mu _0$ be a point of $S'$ in the right half-plane,
say $\mu _0\in \Gamma _{1,3}$ with $a_1,a_3$ dominating above and below
this point respectively. If
$${\rm dist\,}(\mu _0,\Gamma _{3,4^+}=\Gamma _{1,2})\ge {C^2h\over \ln
{1\over \vert \mu _0\vert }},$$
then in $D(\mu _0,Ch/(\ln (1/\vert \mu \vert _0\vert )))$, we have
$$\max (\vert a_1\vert ,\vert a_3\vert )\ge e^{(Ch)^{-1}(\ln {1\over
\vert \mu _0\vert }){\rm dist (\mu _0,\Gamma _{1,2})}}\max (\vert
a_2\vert ,\vert a_4\vert ),$$ and we conclude that the zeros of
$G$ are exponentially close to those of $a_1+a_3$. Essentially the
same statement holds when $\mu _0$ belongs to the lower part of
$S'$, but here the size of $\Re \mu $ also matters, so near a
point $\mu _0\in \Gamma _{3,4^+}\cap S'$ we get a bijection $b$
between the zeros of $a_3+a_{4^+}$ and those of $G$ with
$$b(\mu )-\mu ={\cal O}(1)({h\over \ln{1\over \vert \mu _0\vert
}})(e^{-{\ln {1\over \vert \mu _0\vert }\over Ch}{\rm dist\,}(\mu ,\Gamma
_{1,3})}+e^{-{2\pi \vert \Re \mu \vert \over h}}),$$
and we have to assume both that ${\rm dist\,}(\mu _0,\Gamma _{1,3})\gg
h/\ln{1\over \vert \mu_0 \vert }$ and that $\vert \Re \mu \vert \gg h$.

\par The analogous statements hold in the left half-plane.
\smallskip
\par\noindent \it More refined analysis in the region $\vert \Re \mu
\vert ={\cal O}(h)$. \rm The study of the upper part ($\Gamma _{1,3}$ or
$\Gamma _{1,2}$) of the skeleton is unchanged in this region, while the
lower part requires more attention in view of the fact that $\vert
a_4\vert $ may be considerably smaller than $\max (\vert a_{4^+}\vert
,\vert a_{4^-}\vert )$ when we are close to a zero of $a_4$. In order to
fix the ideas we assume that $\gamma _{3,4^+}(0)\le \gamma _{2,4^+}(0)$.

\par After multiplication of $a_{4},a_{4^\pm},a_3$ by the same exponential
factor, we arrive at \ekv{Sk.22}
{\widetilde{a}_4=\widetilde{a}_{4^+}+\widetilde{a}_{4^-}=2\cosh
{\pi \mu \over h},\ \widetilde{a}_{4^\pm}=e^{\pm\pi \mu \over h},\
\widetilde{a}_3=e^{{i\over h}\phi (\mu ;h)},} and we shall drop
the tildes in the following discussion. Here
\ekv{Sk.23} {\phi
(\mu ;h)=S_{3,4}(\mu )-i{\pi \mu \over 2}+ih{\cal O}_-({h\over \mu
})+\mu \ln {\mu \over i}+\frac{\pi h}{4}-\mu,} with
\ekv{Sk.24}
{-\Im \phi (\mu
;h)=-\Im S_{3,4}(\mu )+\Re {\pi \mu\over 2 }-h\Re {\cal
O}_-({h\over \mu })+\Im \mu \ln {1\over \vert \mu \vert }-\Re \mu
{\rm arg\,}({\mu \over i})+\Im \mu.}

\par We have
\begin{eqnarray*}
\partial _{\Im \mu }(-\Im \phi )&=&\ln {1\over \vert \mu \vert }+{\cal
O}(1),\\
\partial _{\Re \mu }(-\Im \phi )&=&{\cal O}(1),\\
\nabla _\mu ^\alpha (-\Im \phi )&=&{\cal O}(\vert \mu \vert ^{1-\vert
\alpha \vert }),\ \vert \alpha \vert \ge 2,
\end{eqnarray*}
and recall that we work in the region $\vert \mu \vert \ge Ch$. Notice
that $h\nabla _\mu  \ln \vert a_3\vert =\nabla _\mu  (-\Im \phi )$.
Similarly, we look at
\ekv{Sk.25}
{
h\partial _\mu \ln  a_4  =\pi {\sinh {\pi \mu \over h}\over
\cosh {\pi \mu
\over h}}=\pi \sqrt{1-{1\over (\cosh {\pi \mu \over h})^2}},}
for a suitable branch of the square root. Also $h\partial _{\overline{\mu
}}\ln a_4 =0$, so
this relation gives a bound for $h\nabla _\mu \ln a_4$ and its real
part; $h\nabla _\mu \ln \vert a_4\vert $. We have
the general estimate
\ekv{Sk.26}
{\vert \cosh z\vert \ge {1\over C}{\rm dist\,}(z,\cosh ^{-1}(0))e^{\vert
\Re z \vert },\hbox{ for }|\Re z|\le {\rm Const.},}
so we get
\ekv{Sk.27}{h\partial _\mu \ln a_4 =\pi {\rm sgn\,}(\Re \mu
)+{\cal O}(1){e^{-{\pi \vert \Re \mu \vert \over h}}\over {\rm dist\,}({\pi
\mu \over h},\cosh ^{-1}( 0) )}.
}

\par Assuming
\ekv{Sk.28}
{{\rm dist\,}(\mu ,{h\over \pi }\cosh^{-1}( 0) )\ge {Ch\over \ln {1\over
\vert \mu \vert }},\hbox{ with }C\gg 1 ,}
we conclude that
\ekv{Sk.29}
{
\vert h\nabla _\mu  \ln \vert a_4\vert \vert \ll h\vert \nabla _\mu \ln \vert
a_3\vert\vert ,
}
and consequently,
\ekv{Sk.30}
{
h\partial _{\Im \mu}\ln {\vert a_3\vert \over \vert a_4\vert
}=(1+o(1))\ln {1\over \vert \mu \vert} +{\cal O}(1),\ h\partial _{\Re
\mu }\ln {\vert a_3\vert \over \vert a_4\vert }={\cal O}(1)+o(1)\ln
{1\over \vert \mu \vert },
}
where ${\cal O}(1)$ denotes terms that are \ufly{} \bdd{} and $o(1)$
denotes terms that tend to 0, when
$${{\rm dist\,}(\mu ,{h\over \pi} \cosh^{-1}\{ 0\} )\over (h/\ln {1\over
\vert \mu \vert })}\to \infty .$$

\par For each zero $\mu _j$ of $\cosh {\pi \mu \over h}$, we introduce the
diamond shaped \neigh{}

\ekv{Sk.31}
{D_j=\{ \mu ;\, \vert \Re \mu \vert +\vert \Im \mu - \Im \mu _j\vert \le
{Ch\over \ln{1\over \vert \mu \vert }}\} ,}
with $C$ large enough so that the preceding estimates apply away from the
union of all the $D_j$. Define $\Gamma _{3,4}$ to be the set of points
with $\vert a_3\vert /\vert a_4\vert =1$ away from the union of all the $D_j$,
with $D_{j_0}$ added, if $D_{j_0}$ has the property that the
distance from this diamond to the points just defined, is zero. $D_{j_0}$
is unique if it exists, since the other points of $\Gamma _{3,4}$ form a
curve $\Im \mu =\gamma _{3,4}(\Re \mu )$, with $\vert \gamma '_{3,4}\vert
\ll 1$. From the above estimates we get
\ekv{Sk.32}
{
\gamma _{3,4^\pm}(\Re \mu )-{\cal O}(h){\ln\ln {1\over \vert \mu \vert
}\over \ln {1\over \vert \mu \vert }}\le \gamma _{3,4}(\Re \mu )\le \gamma
_{3,4^\pm}(\Re \mu )+{{\cal O}(h)\over \ln {1\over \vert \mu \vert }},\
\pm \Re \mu \ge 0.
}
(If $\Gamma _{3,4}$ stays away from an $h/C$-\neigh{} of the zeros of
$\cosh {\pi \mu \over h}$, then the agreement is better:
$$\gamma _{3,4}(\Re \mu )=\gamma _{3,4^\pm}(\Re \mu )+{{\cal O}(h)\over
\ln {1\over \vert \mu \vert }},\ \pm \Re \mu \ge 0.)$$
In fact, we can get an even more precise estimate for the distance
between $\Gamma _{3,4}$ and $\Gamma _{3,4^\pm}$: Let $\mu _0\in \Gamma
_{3,4}$, and put
$$d(\mu _0)=\max ({h\over \ln {1\over \vert \mu _0\vert }},{\rm dist\,}(\mu
_0,a_4^{-1}(0))).$$
Then
$${d(\mu _0)\over Ch}\le {a_4(\mu _0)\vert \over \vert a_{4^\pm}(\mu
_0)\vert }\le C,$$ away from the diamonds. We therefore get the
following estimate for the vertical distance from $\mu _0$ to
$\Gamma _{3,4^\pm}$: \ekv{Sk.33} { \vert \gamma _{3,4}(\Re \mu _0
)-\gamma _{3,4^\pm}(\Re \mu _0)\vert \le {Ch\over \ln {1\over
\vert \mu _0\vert }}\ln {h\over d(\mu _0 )},} (assuming for
simplicity that $d(\mu _0)\le h/2$). This is a refinement of the
lower bound in (\ref{Sk.32}), and the argument also gives the
upper bound there.

\par We reach the following conclusion about the location of the zeros in
the region $\vert \Re \mu \vert \le Ch$:
\begin{prop}\label{PropSk.3}
\par\noindent -- Above $S'$ and at distance $\ge Ch/\ln {1\over \vert \mu \vert }$ from
$S'$, $a_1(\mu )$ is dominating.
\par\noindent -- $a_4$ is dominating if $\mu $ is below $S'$, at distance
$\ge Ch/\ln {1\over \vert \mu \vert }$ from $a_4^{-1}(0)$
and at distance $\ge Ch {\ln\ln \over\ln }({1\over \vert \mu \vert })$
from $S'$.
\par\noindent -- In between (for instance below $\Gamma _{1,3}$ but above
$\Gamma _{3,4^\pm}$), $a_3$ (or $a_2$) is dominating if the
distance to the skeleton is $\ge Ch/\ln {1\over \vert \mu \vert
}$.
\end{prop}
\medskip
\par\noindent \it Improvement in the region $\Re \mu \gg h$. \rm Let us recall
that \ekv{Sk.34} {a_1a_{4^+}=a_2a_3.} Therefore,  \ekv{Sk.35} {
a_1+a_2+a_3+a_{4^+}=a_{4^+}\left(1+{a_2\over
a_{4^+}}\right)\left(1+{a_3\over a_{4^+}}\right) .} The zeros of
$1+(a_2/a_{4^+})$ are situated on $\Gamma _{2,4^+}$ and are given
by the explicit quantization condition \ekv{Sk.36} { \mu \ln \mu
-\mu +{\pi h\over 4}+S_{1,2}+ih{\cal O}_-({h\over \mu })=2\pi
h(k+{1\over 2}),\ k\in {\bf Z}. } The distance between successive
zeros is $\sim h/\ln {1\over \vert \mu \vert }$. If $\mu _0$ is
such a zero, then in a disc $D(\mu _0,r)$ with $r\ll h/\ln {1\over
\vert \mu _0\vert }$, we have
$$\abs{1+{a_2\over a_{4^+}}}\sim \vert \mu -\mu _0\vert {\ln{1\over
\vert \mu _0\vert }\over h}.$$
Away from the union of all such discs, we have
$$\abs{1+{a_2\over a_{4^+}}} \ge {1\over {\cal O}(1)}.$$
Similarly, the zeros of $1+(a_3/a_{4^+})$ are situated on the curve $\Gamma
_{3,4^+}$ and given by the quantization condition
\ekv{Sk.37}
{
\mu \ln \mu  -\mu +{\pi h\over 4}+S_{3,4}+ih{\cal
O}_-({h\over \mu })=2\pi h(k+{1\over 2}),\ k\in {\bf Z},
}
and the other statements about $1+(a_2/a_{4^+})$ carry over to
$1+(a_3/a_{4^+})$.

\par Now consider
\ekv{Sk.38} { G(\mu ;h)=a_{4^+}\left[\left(1+{a_2\over
a_{4^+}}\right)\left(1+{a_3\over a_{4^+}}\right)+{a_{4^-}\over
a_{4^+}}\right]=a_{4^+}\left[\left(1+{a_2\over
a_{4^+}}\right)\left(1+{a_3\over a_{4^+}}\right)+e^{-{2\pi \mu
\over h}}\right].} We get
\begin{prop}\label{PropSk.4}
In the region $\Re \mu \gg h$, there is a bijection $b$ from the union of
the zeros of $1+a_2/a_{4^+}$ and of $1+a_3/a_{4^+}$ to the zeros of $G$
with
\ekv{Sk.39}
{b(\mu )-\mu ={\cal O}(1){h\over \ln {1\over \vert \mu \vert }}e^{-{\pi \Re \mu \over h}}.}
\end{prop}

So in the region $\Re \mu \gg h$, and modulo an exponentially
small error, we can identify the zeros of $G$ with the union of the zeros of
$1+a_2/a_{4^+}$ and of $1+a_3/a_{4^+}$.

This finishes the analysis of the skeleton in the first case
(\ref{Sk.6}).

\smallskip
\noindent {\it Case 2:}  Assume that \ekv{Sk.40} {\abs{{\rm
arg\,}\mu +{\pi \over 2}} \le \pi -{1\over C},\ h\ll \vert \mu
\vert \ll 1. } In this case from (\ref{as.6}) and (\ref{as.7}),
respectively, we recall that
$$a_{2,3}=2 \cosh \left(\frac{\pi \mu}{h}\right)\exp \left[{\cal O}_+({h\over \mu })+{i\over
h}(\mu \ln (i\mu )-\mu +{\pi h\over 4})+{\pi \mu \over
2h}\right],$$ and
$$a_{1,4}=\exp \left[-{\cal O}_+({h\over \mu })-{i\over
h}(\mu \ln (i\mu )-\mu +{\pi h\over 4})+{\pi \mu \over
2h}\right].$$ Using (\ref{Sk.1}), we get \ekv{Sk.41} {F(\mu
;h)=e^{\pi \mu \over 2h}G(\mu ;h),\ G(\mu ;h)=a_1+a_2+a_3+a_4,}
where $a_j$, $j=1,2,3,4$ are the same as in case 1, but now  with
a partially different representation:
\begin{eqnarray}
a_1&=&a_{1^+}+a_{1^-},\\
a_{1^\pm}&=&e^{{i\over h}(S_{1,2}+S_{3,4})+{\cal
O}_+({h\over \mu })+{i\over h}(\mu \ln (i\mu )-\mu +{\pi h\over 4})\pm {\pi
\mu \over h}},
\nonumber\\
a_2&=&e^{{i\over h}S_{1,2}+{\pi \mu \over 2h}},\nonumber\\ a_3&=&e^{{i\over
h}S_{3,4}+{\pi \mu \over 2h}},\nonumber\\
a_4&=&e^{-{\cal O}_+({h\over \mu })-{i\over h}(\mu \ln
(i\mu )-\mu +{\pi h\over 4})}\nonumber
\end{eqnarray}

\par Again we consider $h$ times the real parts of the different
exponents of the $a_j$, $j=1^{\pm},2,3,4$:
\begin{eqnarray}\label{Sk.41.2}
r_{1^\pm}&=&-\Im S_{1,2}-\Im S_{3,4}+(\Im \mu )\ln {1\over \vert \mu \vert }
-\widetilde{Y}(\mu )\pm \pi \Re \mu ,\nonumber \\
r_2&=&-\Im S_{1,2}+{\pi \over 2}\Re \mu ,\nonumber\\ r_3&=&-\Im S_{3,4}+{\pi \over 2}\Re
\mu \nonumber\\
r_4&=&-(\Im \mu )\ln {1\over \vert \mu \vert }+\widetilde{Y}(\mu ),\nonumber
\end{eqnarray}
where
\ekv{Sk.41.3}
{\widetilde{Y}(\mu )=(\Re \mu ){\rm arg\,}(i\mu )-\Im \mu
-h\Re {\cal O}_+({h\over \mu }).}

\par By a symmetry argument, we shall now see that
$(r_{1^\pm},r_2,r_3,r_4)$ plays the same role in the present case
2 as $(r_{4^\pm},r_3,r_2,r_1)$ in case 1, provided that we perform
the following transformations on
$(r_{1^\pm},r_2,r_3,r_4)$:\smallskip
\par\noindent 1) Add $\Im (S_{1,2}+S_{3,4})$ to each of the five terms.
\smallskip
\par\noindent 2) Replace $\mu $ by $\overline{\mu }$. Then we get
$\widetilde{r}_j(\overline{\mu })=(\Im (S_{1,2}+S_{3,4})+r_j)(\mu )$:
\begin{eqnarray}
\widetilde{r}_{1^\pm}(\overline{\mu })&=& -(\Im \overline{\mu
})\ln {1\over \vert \overline{\mu } \vert }+(\Re \overline{\mu
}){\rm arg\,}({\overline{\mu}\over i})-\Im \overline{\mu} +h\Re
{\cal O}_+({h\over \mu })\pm \pi \Re
\overline{\mu }, \\
\widetilde{r}_2(\overline{\mu })&=&\Im S_{3,4}+{\pi \over 2}\Re
\overline{\mu },\nonumber\\ \widetilde{r}_3(\overline{\mu })&=&\Im
S_{1,2}+{\pi \over
2}\Re \overline{\mu },\nonumber\\
\widetilde{r}_4(\overline{\mu })&=&\Im S_{1,2}+\Im S_{3,4}+(\Im
\overline{\mu } )\ln {1\over \vert \overline{\mu } \vert }-(\Re
\overline{\mu } ){\rm arg\,}({\overline{\mu}\over i} )+\Im
\overline{\mu }-h\Re {\cal O}_+({h\over \mu }).\nonumber
\end{eqnarray}

This is analogous with $(r_{4^\pm},r_3,r_2,r_1)$ in case 1 except
for the fact that $\Re {\cal O}_+({h\over \mu })$ here corresponds
to $\Re {\cal O}_-({h\over \mu })$ in case 1.
\medskip
\par\noindent \it Remark. \rm Using (\ref{as.4}), it is easy to check
that \ekv{Sk.42} { \widetilde{Y}(\mu )-Y(\mu )=\pm \pi \Re \mu
+h{\cal O}(e^{-2\pi \vert \Re \mu \vert /h}), } when \ekv{Sk.43}
{\vert {\rm arg\,}(\pm \mu )\vert \le {\pi \over 2}-{1\over C},\
\vert \mu \vert \ge h.} It follows that \ekv{Sk.44} { r_j(\mu
)=r_{j^\pm}(\mu )+h{\cal O}(e^{-2\pi \vert \Re \mu \vert /h}),\
j=1,4, } when $\mu $ satisfies (\ref{Sk.43}), and hence if $\mu $
belongs to the skeleton $S'$ defined according to Case 1, the
distance from $\mu $ to the corresponding skeleton $S'$ defined
according to Case 2 is
$${\cal O}\left({h\over \ln {1\over \vert \mu \vert }}e^{-2\pi \vert \Re \mu
\vert /h}\right).$$

\bigskip
\par We end this section by some general considerations that will be useful
in Section \ref{SectionBt}.
We see from (\ref{Sk.8.5}) that the spectrum will have a genuinely
2-dimensional structure if
\ekv{Sk.45}
{
\vert \Im S_{3,4}(0)-\Im S_{1,2}(0)\vert \gg h,
}
or if $\Im S_{3,4}(0)$ and $\Im S_{1,2}(0)$ have the same sign and
\ekv{Sk.46}
{
\min (\vert \Im S_{3,4}(0)\vert , \vert \Im S_{1,2}(0)\vert )\gg
h\ln{1\over h}.
}
In the latter case, we even have some \ev{}s on the imaginary $\mu $-axis,
related to 1-dimensional barrier top \res{}s. It is therefore important to
have a \sufly{} invariant and direct description of $\Im
S_{3,4}(0)$, $\Im S_{1,2}(0)$.

\par The final definition of $S_{3,4}$ in the beginning of
section 6 is simply that we start with the null-solution $f_4$ of
$Q$ near $\alpha _4$ and extend it along the exterior part of
$K_{0,0}$ until we reach a \neigh{} of $\alpha_3$, where we get
$\exp ({i\over h}S_{3,4})f_3$. (Here we neglected the real Floquet
parameter $\theta _1$, since we are only interested in the
imaginary part of $S_{3,4}$). The definition of $S_{1,2}$ is
similar.

\par Now take $\mu =0$ (cf. (\ref{1D.33}) and (\ref{1D.34})) and represent the \op{} $Q$ as acting in $H_\Phi
^{\rm loc}(\Omega )$, where $\Phi $ is \st{} \plsh{}, with
$\Lambda _\Phi \simeq$ a \neigh{} of $K_{0,0}$ in $T^*{\bf R}$.
From the construction of $e_j$, $f_j$, we see that $f_j$ is near
$\alpha _j$ a normalized null-solution of $Q$ in $H_{\Phi _0}^{\rm
loc}(\Omega )$, where $\Phi -\Phi _0$ is small and $\Phi _0$ is
defined in a \sufly{} large \neigh{} of the projection of the
branching point. Here if $\kappa_T$ is the canonical transformation
associated to some standard FBI-Bargmann transform. then
$\Lambda_{\Phi}=\kappa_T(\real^2)$,
$\Lambda_{\Phi_0}=\kappa_T\circ \kappa_U(\real^2)$, with $U$ as in
(\ref{1D.33}). Since $\Lambda _{\Phi
_0}=\{\xi=\frac{2}{i}\frac{\partial \Phi_0}{\partial x}\}$ is
naturally identified with $T^*{\bf R}$, where $p_0=x\xi $, so
(since $\mu =0$), we know that the null set of $Q_0$ intersects
$\Lambda _{\Phi _0}$ along two crossing "real" curves that we we
can identify with "the interior part" of $K_{0,0}$.

\par Extend $\Phi _0$ to be defined in  $\Omega $ with $\Phi -\Phi _0$
small. If $S_{3,4}\sim S_{3,4}^0+hS_{3,4}^1+...$, then \ekv{Sk.47}
{ -\Im S_{3,4}^0=\int_{\gamma _{3,4}}(-\Im (\xi \cdot dx)-d\Phi
_0), } where $\gamma _{3,4}$ now (cf (\ref{1D.44})) is a real
curve from $\alpha _4$ to $\alpha _3$ in $Q_0^{-1}(z(0))$ close to
the exterior part of the "left loop" of $K_{0,0}$.  Here we let
$z(0)$ denote the $z$-value in \nr{1D.34} corresponding to $\mu
=0$.  Let this left loop be denoted by $\gamma ^1$ and let us
consider it (after slight deformation) as a closed curve in
$Q_0^{-1}(z(0))$ joining the critical point of $Q_0$ to itself
staying close to the left loop of $K_{0,0}$. Here, we may assume
that the interior part of $\gamma^1$ (joining $\alpha _3$ to
$\alpha _4$) is contained in $\Lambda _{\Phi _0}$, so $-\Im (\xi
\cdot dx)=d\Phi _0$ there. Hence (\ref{Sk.47}) becomes \ekv{Sk.48}
{\Im S_{3,4}^0=\int_{\gamma^1}(\Im (\xi \cdot dx)+d\Phi
_0)=\int_{\gamma^1}(\Im \xi \cdot dx ) .} Here $\Im (\xi \cdot
dx)$ in the last integral can be replaced by $\Im (\xi \cdot
dx)+d\Phi $, which by Stokes' formula can be further replaced by
any other real 1-form $\omega  $ with $d\omega =\Im \sigma $,
${\omega _\vert}_{\Lambda _\Phi }=0$. This means that we can
reinterpret the last integral in (\ref{Sk.48}) as the
corresponding one along the corresponding closed curve in the
complexification of ${\bf R}^2$.

To simplify things further, recall that $\epsilon $, $h^2/\epsilon
$ are small perturbative parameters for $Q_0$ and that
$Q_0^{-1}(z(0))$ is real when $\epsilon =h^2/\epsilon =0$. In
general, if $q=q_s$ depends smoothly on a real parameter $s$, if
$E$ is not a critical value and $\gamma =\gamma (s,E)$ is a simple
closed curve in $q_s^{-1}(E)$, then for $E$ fixed; \ekv{Sk.49} {
{\partial \over
\partial s}\int_\gamma \xi dx=-\int_0^{T(E,s)}{\partial q\over
\partial s}(x(t),\xi (t))dt, } where $[0,T(E,s)]\ni t\mapsto \exp
tH_{q_s}(\rho (0))$ is a natural parametrization. This can be
applied to the case $q_s=q-E(s)$, so if $E$ also depends on $s$,
we get
$$
{\partial \over \partial s}\int_\gamma \xi \cdot
dx=\int_0^{T(E,s)}\left({\partial E(s)\over \partial s}-{\partial q\over
\partial s}(x(t),\xi (t))\right)dt.
$$
In this form, we can treat a loop like $\gamma _1(s)\subset
q_s^{-1}(E(s)),$ starting and ending at the critical points $\rho
_c(s)$ of $q_s$, parametrized by $]-\infty ,+\infty [\mapsto \exp
(tH_{q_s})(\rho (0))$, provided that we take $E(s)$ equal to the
critical value $q_s(\rho _c(s))$: \ekv{Sk.50} { {\partial \over
\partial s}\int_{\gamma _1(s)}\xi dx=\int_{-\infty }^{+\infty
}\left({\partial \over \partial s}(q_s(\rho _c(s)))-({\partial
\over \partial s}q_s)(x(t),\xi (t))\right)dt. } This can be proved
by a limiting procedure, approaching $\gamma _1(s)$ by closed
curves at non-critical levels.

\par Taking the imaginary parts, this means that we have a fairly simple
way of computing $\Im S_{1,2}$, $\Im S_{3,4}$ perturbatively. From
this computation, we see that it is of interest to compute the
$\epsilon ^2$ contribution to the averaged principal symbol.  This
computation was carried out in Section 2 of \cite{HiSj2} under the
assumption that $\langle q\rangle =0$ and it works essentially the
same way without that assumption: We start with the principal
symbol \ekv{Sk.51} { p_\epsilon =p+i\epsilon q+\epsilon ^2r+{\cal
O}(\epsilon ^3). } The function \ekv{Sk.52} { G_0={1\over
T(E)}\int_0^{T(E)} \left(t-{T(E)\over 2}\right)q\circ \exp
(tH_p)dt,\,\,\wrtext{on}\,\, p^{-1}(E) } introduced in Proposition
\ref{PropRed1}, satisfies \ekv{Sk.53} {H_pG_0=q-\langle q\rangle
.} Put $G=G_0+i\epsilon G_1+{\cal O}(\epsilon ^2)$, where $G_1$
remains to be determined. As in \cite{HiSj2} we get at a general
point $\exp (i\epsilon H_G)(\rho )\in \Lambda _{\epsilon G}$,
($\rho \in T^*M$):

\begin{eqnarray*}
{{p_\epsilon }_\vert}_{\Lambda _\epsilon }&\simeq& p_\epsilon (\exp
(i\epsilon H_G)(\rho ))=\sum_0^\infty {(i\epsilon H_G)^k\over
k!}p_\epsilon (\rho )\\
&=& p+i\epsilon \langle q\rangle +\epsilon ^2(r+H_pG_1-H_{G_0}({1\over
2}(q+\langle q\rangle )))+{\cal O}(\epsilon ^3),
\end{eqnarray*}
where we used that $H_{G_0}^2p=-H_{G_0}(q-\langle q\rangle )$.

\par Letting $G_1$ solve
$$
H_pG_1=H_{G_0}({1\over 2}(q+\langle q\rangle ))-\langle H_{G_0}({1\over
2}(q+\langle q\rangle ))\rangle -(r-\langle r\rangle ),
$$
we get with $G=G_0+i\epsilon G_1$,
\ekv{Sk.54}
{
p_\epsilon (\exp (i\epsilon H_G)(\rho ))=p+i\epsilon \langle q\rangle
+\epsilon ^2\langle s\rangle +{\cal O}(\epsilon ^3).
}

\par Now assume for simplicity that $T(E)=T$
is constant. Then
\begin{eqnarray*}
\langle s\rangle &=&\langle r\rangle -{1\over 2T}\int_0^{T}\{
G_0,q+\langle q\rangle \}\circ \exp (tH_p)dt\\
&=& \langle r\rangle -{1\over 2T^2}\int_0^{T}\int_0^{T}
(s-{T\over 2})\{ q\circ \exp (t+s)H_p,(q+\langle q\rangle )\circ \exp
(tH_p)\} dtds \\
&=&\langle r\rangle -{1\over 2T}\int_0^{T}(s-{T\over 2})\langle \{
q\circ \exp (sH_p),q+\langle q\rangle \} \rangle ds.
\end{eqnarray*}
Here, we notice that
$$
\langle \{ q\circ \exp (sH_p),\langle q\rangle \} \rangle =\{ \langle
q\circ \exp (sH_p), \langle q\rangle \} =\{ \langle q\rangle ,\langle
q\rangle \}=0,
$$
so finally:
\ekv{Sk.55}
{
\langle s\rangle =\langle r\rangle -{1\over 2T}\int_0^{T}(s-{T\over 2})\langle \{
q\circ \exp (sH_p),q\}\rangle ds.}

\par The formulas (\ref{Sk.54}), (\ref{Sk.55}) will be used in Section \ref{SectionBt}
together with the following remark: If we put
\ekv{Sk.56}
{
{\rm Cor\,}(q_1,q_2;s)=\langle \{ q_1\circ \exp (sH_p),q_2\} \rangle ,
}
then a simple computation shows that
\ekv{Sk.57}
{
{\rm Cor\,}(q_1,q_2;s)=-{\rm Cor\,}(q_2,q_1;-s).
}
If we put
\ekv{Sk.58}
{
C(q_1,q_2)={1\over T}\int_0^{T}(s-{T\over 2}){\rm Cor\,}(q_1,q_2;s)ds,
}
then combining (\ref{Sk.57}) and the $T$ periodicity of ${\rm
Cor\,}(q_1,q_2;s)$ with the change of variables $T/2-s=\widetilde{s}-T/2$,
we get
\ekv{Sk.59}
{
C(q_1,q_2)=C(q_2,q_1).
}

\section {Skeleton for $\vert \mu\vert  \le {\cal O}(h)$}\label{Sectionsk}
\setcounter{equation}{0}

\par In this section we shall consider the case $\vert \mu \vert \le
{\cal O}(1)h$. In doing so, we will use (\ref{tr.12}) more directly.

\smallskip
\par\noindent \it Case 1. \rm
We will  first work in a region $\{ \mu \in{\bf C};\, \vert \mu
\vert <rh\}\cup \{ \mu \in{\bf C}\setminus \{0\};\, \vert {\rm
arg\,}\mu -{\pi \over 2}\vert <\pi  -1/C\}$, where $0<r<1/2$,
$C>0$. (The corresponding region with $\vert {\rm arg\,}\mu +{\pi
\over 2}\vert \le \pi -1/C$, can be treated with a symmetry
argument as in the end of Section \ref{SectionSk}, and this
argument will be given later.) It follows from (\ref{tr.12}) that here $a_{2,3}\ne 0$ and $\ln \Gamma
({1\over 2}-i{\mu \over h})$ is well-defined, while $\Gamma
({1\over 2}+i{\mu \over h})^{-1}$ may have zeros. Consequently we
use the reflection identity, to get
\begin{eqnarray*}
a_{1,4}&=&{\sqrt{2\pi
}\over \Gamma ({1\over 2}+i{\mu \over h})}h^{-i{\mu \over h}}e^{{\pi\mu
\over 2h}-{i\pi \over 4}}\\
&=& {\Gamma ({1\over 2}-i{\mu \over
h})\over \sqrt{2\pi }}h^{-i{\mu \over h}}e^{{\pi\mu  \over
2h}-{i\pi \over 4}}2\cosh ({\pi \mu \over h}).
\end{eqnarray*}

\par Now using (\ref{Sk.1}) we get
\begin{eqnarray*}
F(\mu ;h)&=&e^{-i{\mu \over h}\ln{1\over h}-\ln {\Gamma ({1\over 2}-i{\mu
\over h})\over \sqrt{2\pi }}+{\pi \mu \over 2h}+{i\pi \over 4}
+{i\over h}(S_{1,2}+S_{3,4})}\\
&&+e^{{i\over h}S_{1,2}+\pi {\mu \over h}}+e^{{i\over h}S_{3,4}+\pi {\mu
\over h}}\\
&&+
e^{\ln {\Gamma ({1\over 2}-i{\mu
\over h})\over \sqrt{2\pi }}+i{\mu \over h}\ln {1\over h}-{i\pi \over
4}+{\pi \mu \over 2h}+\ln 2\cosh {\pi
\mu \over h}}\\
&=&e^{\pi \mu \over 2h}G(\mu ;h),
\end{eqnarray*}
where
\begin{eqnarray}\label{sk.1}
G(\mu ;h)&=&e^{{i\over h}(S_{1,2}+S_{3,4})-i{\mu \over h}\ln{1\over h}-
\ln {\Gamma ({1\over 2}-i{\mu \over h})\over \sqrt{2\pi }}+
{i\pi \over 4}}\\
&&+e^{{i\over h}S_{1,2}+\pi {\mu \over 2h}}+e^{{i\over h}S_{3,4}+\pi {\mu
\over 2h}}\nonumber\\
&&+
e^{i{\mu \over h}\ln {1\over h}+\ln {\Gamma ({1\over 2}-i{\mu
\over h})\over \sqrt{2\pi }}-{i\pi \over
4}}2\cosh {\pi
\mu \over h}.\nonumber\\
&=&a_1+a_2+a_3+a_4,\quad\quad\quad a_4=a_{4^+}+a_{4^-},\nonumber
\end{eqnarray}
where the terms are the same as in (\ref{Sk.1.5}), although we shall now use
different \asy{} approximations.

\par Again we introduce $h$ times the real parts of the different exponents:
\begin{eqnarray}\label{sk.2}
r_1&=& -\Im S_{1,2}-\Im S_{3,4}+(\Im \mu )\ln {1\over h} -h\Re
\ln {\Gamma ({1\over 2}-i{\mu \over h})\over \sqrt{2\pi }},\\
r_2&=&-\Im S_{1,2}+{\pi \over 2}\Re \mu , \nonumber\\
r_3&=&-\Im S_{3,4}+{\pi \over 2}\Re \mu ,\nonumber\\
r_{4^\pm}&=&-(\Im \mu )\ln {1\over h}+h\Re \ln {\Gamma ({1\over 2}-i{\mu \over
h})\over \sqrt{2\pi }}\pm \pi\Re \mu.\nonumber
\end{eqnarray}
As before, we have
\ekv{sk.3} {r_2+r_3=r_1+r_{4^+}.}

\par Again, we define the different curves $\Gamma _{j,k}$ by $\vert
a_j\vert =\vert a_k\vert $, for $j\ne k\in\{ 1,2,3,4^\pm\}$ with
the exception of $(j,k)=(2,3)$ and $(j,k)=(4^+,4^-)$.  (The
segment $\Gamma _4$ is now defined to be the segment of the
positive imaginary axis joining 0 to the lower part of $S'$,
provided that this lower part is not hidden in the forbidden
region, in which case we let $\Gamma _4$ be empty).  More
explicitly, we get: \ekv{sk.4} {(\Im \mu )\ln {1\over h}=\cases{
\Im S_{3,4}\pm \pi \Re\mu -{\pi \over 2}\Re \mu +X\hbox{ on
}\Gamma _{3,4^\pm}\cr \Im S_{1,2}\pm \pi \Re\mu -{\pi \over 2}\Re
\mu +X\hbox{ on }\Gamma _{2,4^\pm}\cr \Im S_{1,2}+{\pi \over 2}\Re
\mu  +X\hbox{ on }\Gamma _{1,3},\cr \Im S_{3,4}+{\pi \over 2}\Re
\mu  +X\hbox{ on }\Gamma _{1,2},\cr {1\over 2}(\Im S_{1,2}+\Im
S_{3,4})\pm {\pi \mu \over 2}+X,\hbox{ on }\Gamma _{1,4^\pm} ,}}
with \ekv{sk.5} {X=h\Re \ln \left({\Gamma ({1\over 2}-i{\mu \over
h})\over \sqrt{2\pi }}\right).} The function $X$ now differs from that
of Section \ref{SectionSk} by a term $-{\pi \over 2}\Re \mu $. The
definition of $\Gamma _{1,3}=\Gamma _{2,4^+}$, $\Gamma
_{1,2}=\Gamma _{3,4^+}$ coincides with that in Section
\ref{SectionSk} in the overlap region.

\par We shall also define a set $\Gamma _{j,4}$ for $j=1,2,3$ as in the
preceding section. To do so, we
check that
\ekv{sk.6}
{h\abs{{\nabla \cosh {\pi \mu \over h}\over \cosh {\pi \mu \over h}}}\ll
\ln {1\over h},}
if
\ekv{sk.7}
{{\rm dist\,}(\mu ,{h\over \pi }\cosh^{-1}(0))\gg {h\over \ln {1\over h}}.}
In this region, we also have
\ekv{sk.8}
{{1\over {\cal O}(1)\ln{1\over h}}\le \abs{a_4\over a_{4^\pm}} \le
{\cal O}(1).}
In the region (\ref{sk.7}) we can define $\Gamma _{j,4}$ by $\vert
a_j\vert =\vert a_4\vert $, and see that we get a curve
\ekv{sk.9}
{\Im \mu =\gamma _{j,4}(\Re \mu ),}
with $|\gamma '_{j,4}\vert \ll 1$. Using (\ref{sk.8}), we also see that if
we represent $\Gamma _{j,4^\pm}$ by $\Im \mu = \gamma _{j,4^\pm}(\Re \mu )$,
then
\ekv{sk.10}
{(\gamma _{j,4}-\gamma _{j,4^\pm})(\Re \mu )={\cal O}(1)h({\ln\ln\over
\ln})({1\over h}).} Actually the upper bound can here be improved to
${\cal O}(h)/\ln {1\over h}$.  In analogy with Section \ref{SectionSk}, we
define a diamond shaped \neigh{} of each zero $\mu _j$ of $a_4$ by
\ekv{sk.11}
{D_j=\left\{ \mu ;\, \vert \Re \mu \vert +\vert \Im \mu - \Im \mu _j\vert \le
{Ch\over \ln{1\over h}}\right\} ,}
The previously defined $\Gamma _{j,4}$ can hit at most one of the $D_\nu $
and if that happens, we add that diamond to the set $\Gamma _{j,4}$.

\par Now define the skeleton as before:
\ekv{sk.12}{S=S'\cup\Gamma _4,}
and as before we can describe the regions of dominance:
\smallskip
\par $a_4$ dominates at distance $\gg h({\ln\ln\over \ln})({1\over h})$
below $\inf_{j=1,2,3}\gamma _{j,4^\pm}$, for $\pm \Re \mu \ge 0$,
intersected with the complement of the union of the diamonds.
\smallskip
\par The other $a_j$ dominate according to the earlier rules in their
respective regions at a distance $\gg h/\ln{1\over h}$ from the skeleton.
\smallskip

\par\noindent \it Case 2. \rm
 We now consider the case
when $\mu $ belongs to the set
\ekv{sk.13}
{
\{ \mu \in{\bf C};\, \vert \mu \vert <rh\} \cup\{ \mu \in {\bf
C}\setminus \{ 0\} ;\, \abs{\arg\,
\mu +{\pi \over 2}} <\pi -{1\over C}\} ,}
where $0<r<{1\over 2}$, $C>0$. From (\ref{tr.12}), we get
$$a_{1,3}=-a_{2,4}=e^{{\pi \mu \over h}+i{\pi \over 2}},$$
$$a_{1,4}={\sqrt{2\pi }\over
\Gamma ({1\over 2}+i{\mu \over h})}h^{-i{\mu \over h}}e^{{\pi \mu \over
2h}-{i\pi \over 4}},$$
where $a_{1,4}$ is non-vanishing, while
$$a_{2,3}={\sqrt{2\pi }\over
\Gamma ({1\over 2}-i{\mu \over h})}h^{i{\mu \over h}}e^{{\pi \mu \over
2h}+{i\pi \over 4}}$$
may have zeros, so we use the reflection identity to write
$$a_{2,3}=
{ \Gamma ({1\over 2}+i{\mu \over h})\over \sqrt{2\pi }}h^{i{\mu \over h}}e^{{\pi \mu \over
2h}+{i\pi \over 4}}2\cosh {\pi \mu \over h}.$$
We then use (\ref{Sk.1}) to get
$$F(\mu ;h)=e^{\pi \mu \over 2h}G(\mu ;h),$$
with
$$G(\mu ;h)=a_1+a_2+a_3+a_4,\quad a_1=a_{1^+}+a_{1^-},$$
$$a_{1^\pm}=\exp\left[
{i\over h}(S_{1,2}+S_{3,4})+\ln \left({\Gamma ({1\over
2}+i{\mu \over h})\over \sqrt{2\pi }}\right)+i{\mu \over h}\ln h+{i\pi \over
4}\pm {\pi \mu \over h}\right],$$
$$a_2=\exp\left[{i\over h}S_{1,2}+{\pi \mu \over 2h}\right],\quad a_3=\exp\left[{i\over
h}S_{3,4}+{\pi \mu \over 2h}\right],$$
$$a_4=\exp\left[-\ln ({\Gamma ({1\over 2}+i{\mu \over h})\over
\sqrt{2\pi }})-i{\mu \over h}\ln h -{i\pi \over 4}\right]$$

\par Again, we introduce $h$ times the real parts of the different
exponents:
$$r_{1^\pm}=-\Im S_{1,2}-\Im S_{3,4}+(\Im \mu )\ln {1\over h}+h\Re\ln
({\Gamma ({1\over 2}+i{\mu \over h})\over \sqrt{2\pi }})\pm\pi \Re \mu ,
$$
$$r_2=-\Im S_{1,2}+{\pi \over 2}\Re \mu ,\quad r_3=-\Im S_{3,4}+{\pi \over
2}\Re \mu ,$$
$$
r_4=-(\Im \mu )\ln{1\over h}-h\Re\ln ({\Gamma ({1\over 2}+i{\mu \over
h})\over \sqrt{2\pi }}).
$$

We shall now make the same symmetry transformations as in Section 8, to see that
the functions $r_{1^\pm}$, $r_2$, $r_3$, and $r_4$ play the same role
as $r_{4^\pm}$, $r_3$, $r_2$, and $r_1$ respectively, in the
previously considered case:
\smallskip
\par\noindent 1) Add $\Im (S_{1,2}+S_{3,4})$ to each of the $r_j$.\smallskip

\par\noindent 2) Consider the $r_j$ as \fu{}s of $\overline{\mu }$. Then we
get
$$\widetilde{r}_{1^\pm}(\overline{\mu })=-(\Im \overline{\mu })\ln{1\over
h}+h\Re\ln ({\Gamma ({1\over 2}-i{\overline{\mu }\over h})\over
\sqrt{2\pi }})\pm \pi \Re \overline{\mu },$$
$$\widetilde{r}_2(\overline{\mu })=\Im S_{3,4}+{\pi \over 2}\Re
\overline{\mu },\quad \widetilde{r}_3=\Im S_{1,2}+{\pi \over 2}\Re
\overline{\mu },$$
$$
\widetilde{r}_4=\Im S_{1,2}+\Im S_{3,4}+(\Im \overline{\mu })\ln{1\over h}
-h\Re \ln ({\Gamma ({1\over 2}-i{\overline{\mu }\over h})\over
\sqrt{2\pi }}).
$$
Thus apart from a change of sign in $\Im S_{1,2},\, \Im S_{3,4}$, we see that
$(\widetilde{r}_{1^\pm},\widetilde{r}_2,\widetilde{r}_3,\widetilde{r}_4)$
has the same properties as $(r_{4^\pm},r_3,r_2,r_1)$ in the previously
considered case.
\smallskip

\par\noindent \it Remark. \rm In the overlap region
$$D(0,rh)\cup\left \{ \mu ;\, \vert \mu \vert \le Ch,\, \vert {\rm arg\,}\mu
\vert \le {\pi\over 2} -{1\over C},\hbox{ or }\vert {\rm arg\,}(-\mu )\vert
\le {\pi \over 2}-{1\over C}\right \}, $$
where both cases apply, we notice that trivially $r_j=r_{j^\pm}+{\cal O}(h)$
for $\pm\Re \mu \ge 0$, $j=1,4$, (and these estimates improve by
(\ref{Sk.44}) when $\vert \mu \vert /h$ increases). As in the remark at
the end of Section \ref{SectionSk},the distance between the two skeletons,
defined according to Case 1 and according to Case 2, is therefore ${\cal
O}(h/\ln {1\over h})$.

\section{Eigenvalue counting}\label{SectionEv}
\setcounter{equation}{0}

\par In each of the cases 1 and 2 of Sections \ref{SectionSk},
\ref{Sectionsk}, we defined a skeleton $S$ consisting of a
horizontal part $S'$ possibly with a vertical part ($\Gamma _4$ in Case 1
and $\Gamma _1$ in Case 2) added. We notice that the definitions in the two
sections agree for each of the cases 1 and 2 in the overlap regions
for the two sections, and we saw in the remarks at the end of the
sections, that if we compare the skeletons for the two cases in the
overlap region
$$\{ \vert \Re \mu \vert >{1\over C}\vert \Im \mu \vert \}\cup D(0,rh),$$
then the distance between the corresponding skeletons is
$${\cal O}({h\over \ln \langle \mu \rangle _h}
e^{-2\pi \vert \Re \mu \vert /h}),\ \langle \mu \rangle
_h:=\sqrt{h^2+\vert \mu \vert ^2}.$$

\par Now define the body by widening the skeleton:
\ekv{ev.1}
{B=\left(\bigcup_{\mu \in S'}D\left(\mu ,{Ch\over \ln {1\over \langle \mu \rangle
_h}}\right)\right)\cup B_v\cup B_e.} Here $B_v,B_e$ may be empty and will now be
defined.  They are non-empty if $S'$ stays entirely in the admissible
regions for one of the cases, and in order to fix the ideas, we assume that
this is Case 1, and $S'$ does not intersect the negative imaginary half
axis.  If so, we have a non-empty segment $\Gamma _4$ in the imaginary
axis, joining 0 to the closest imaginary point of $S'$.  Recall that we
have defined the diamonds $D_j$ around the zeros of $a_4$ in $\Gamma _4$,
by (\ref{Sk.31},\ref{sk.11}) .  We define $B_v$ to be the union of $\Gamma
_4$ (in Case 1, and $\Gamma _1$ in Case 2) and the corresponding diamonds.
(In Case 2 we do the corresponding definition with "4" replaced by "1" and
$\Gamma _1$ is then the segment in the negative imaginary axis, joining 0
to the closest part of $S'$.  ) $B_e$ is non-empty precisely when $\Gamma
_4$ or $\Gamma _1$ is.  In Case 1, it is defined to be the set of points
$\mu $ below $S'$ at distance at most $$ Ch{\ln\ln\over \ln}\big({1\over
\langle \mu \rangle _h}\big) $$ from $S'$ with $C>0$ \sufly{} large and
with $\vert \Re \mu \vert <h$.  Here the upper bound $h$ in the last
estimate, may be replaced by $h/C_0$ for any fixed $C_0>0$, and we could
decrease $B_e$ further by a more detailed discussion.  In Case 2 we have
the analogous definition.

\par We next define what we mean by an \it admissible curve. \rm It should
be a piecewise $C^1$-curve $\gamma :[a,b]\to{\bf C}$ without self-intersections,
parametrized by arc-length. It is tacitly assumed that we consider a family
of such curves, which is \ufly{} \bdd{} in the sense that we have uniform
bounds on the number of jump discontinuities of $\dot{\gamma }$, the
continuity of $\dot{\gamma }$ between the discontinuities, and on the
length $b-a$. It is also required that $\gamma (t)$ may belong to $B$
only for $t\in I_j$, $j=1,2,...,M$, where $I_j$ are disjoint intervals of length
$\le Ch/\ln {1\over \langle \mu  \rangle _h}$, for some $\mu =\gamma (t)\in
I_j$, if $\gamma (I_j)\cap B_e =\emptyset $ and of length
$\le Ch {\ln\ln\over \ln}({1\over \langle \mu \rangle _h})$ otherwise. We
also assume that we have a uniform bound on the number $M$ of such
intervals.

\par Assume for simplicity that $a,b\notin \cup I_j$ and let us
partition $[a,b]$ into intervals in increasing order:
\ekv{ev.2}
{[a,b]=J_0\cup I_1\cup J_1\cup I_2\cup...\cup I_M\cup J_M.}
For each $J_k$, let $a_{\nu (k)}$, be the corresponding dominant term along
$\gamma (J_k)$. For
simplicity we shall assume that the image of $\gamma $ is entirely
contained in the admissible region for one of the cases 1 or 2, so
that $\nu (k)$ is either in $\{ 1,2,3,4^\pm \}$ (Case 1), or in $\{1^\pm
,2,3,4\}$ (Case 2). Let
$$\mu _{k,e}=\mu _{k+1,s}=\gamma (t_{k+1})\hbox{ for some }t_{k+1}\in
I_{k+1},\ k=0,..,M-1,$$
and put
$$\mu _{0,s}=\gamma (a),\, \mu _{M,e}=\gamma (b).$$
Then we have, with $a_j=e^{i\phi _j/h}$ (cf.~(\ref{Sk.1.5}),
(\ref{Sk.41}), (\ref{sk.1})):
\begin{theo}\label{ThEv.1} Let $\gamma $ be an admissible curve as above.
Then
\begin{eqnarray}\label{ev.3}
& & \Re {1\over 2\pi i}\int_\gamma  {G'\over G}d\mu = \Re {1\over 2\pi
h}\big(\big( \phi _{\nu (M)}(\mu _{M,e})+\\& & \sum_{k=0}^{M-1} (\phi _{\nu (k)}(\mu
_{k,e})-\phi _{\nu (k+1)}(\mu _{k+1,s}))-\phi _0(\mu _{0,s})\big)
+{\cal O}(1)+{\cal O}(\max\ln\ln {1\over \langle \mu
_{k,e}\rangle _h})\big),\nonumber
\end{eqnarray}
where the maximum is taken over all $k$ with $\gamma (I_k)\cap
B_e\ne\emptyset $, so if $\gamma $ never meets $B_e$, we only have
the remainder
${\cal O}(1)$.
\end{theo}

\begin{proof}
Notice that the first term of the \rhs{} of (\ref{ev.3}) can also be
written
\ekv{ev.4}
{
\Re {1\over 2\pi h}\sum_{k=0}^{M}(\phi _{\nu (k)}(\mu _{k,e})-\phi
_{\nu (k)}(\mu _{k,s})).
}
Consider an interval $J_k$. If we first assume $\nu (k)\ne 4^\pm$ (if we
are in Case 1), then for $t\in J_k$:
\ekv{ev.5}
{G(\gamma (t))=b_{\nu (k)}(\gamma (t))a_{\nu (k)}(\gamma (t)),\ \vert
b_{\nu (k)}(\gamma (t))-1\vert <{1\over 2}.}

\par
Let $\widetilde{\mu }_{k,s}$ and $\widetilde{\mu }_{k,e}$ be the start and
the end points of ${\gamma _\vert}_{J_k}$, so that
\begin{eqnarray}\label{ev.6}
\widetilde{\mu}_{k,s}&=&\mu _{k,s}+{\cal O}(h{\ln\ln\over \ln}({1\over
\langle \mu_{k,s} \rangle_h })),\\
\widetilde{\mu}_{k,e}&=&\mu _{k,e}+{\cal O}(h{\ln\ln\over \ln}({1\over
\langle \mu_{k,e} \rangle_h })),\nonumber
\end{eqnarray}
with the ${\ln\ln\over \ln}$ improving to ${1\over \ln}$ if the
corresponding neighboring interval $I_{..}$ does not meet $B_e$. Using
(\ref{ev.5}), we see that
\ekv{ev.7}
{{1\over 2\pi i}\int_{{\gamma _\vert}_{J_k}}{G'\over G}d\mu ={\cal
O}(1)+\phi _{\nu (k)}(\widetilde{\mu }_{k,e})-\phi _{\nu
(k)}(\widetilde{\mu }_{k,s}).}

\par In the case $\nu (k)=4^\pm$, we know that $a_4$ is dominating along
${\gamma _\vert}_{J_k}$ and (\ref{ev.7}) holds with $\phi _{\nu (k)}$
replaced by ${h\over i}\ln a_4$. Now we also know that along
${\gamma _\vert}_{J_k}$, we have
$$a_4(\gamma (t))=c(\gamma (t))a_{\nu (k)}(\gamma (t)),$$
with
$$1/(C\ln {1\over \langle \mu \rangle _h})\le \vert c(\gamma
(t))\vert \le C,$$
and with ${\rm arg\,}c(\gamma (t))$ of \bdd{} variation. It follows that
the real part of the equation (\ref{ev.7}) still holds in this case.

\par We next estimate the integral along ${\gamma _\vert}_{I_k}$, and let
us consider the worst case, when $\gamma (I_k)\cap B_e \ne
\emptyset $. Let
$$r={\cal O}(h{\ln\ln\over \ln}({1\over \langle \mu \rangle _h}))$$
be such that
$$\gamma (I_k)\subset D(\widetilde{\mu }_{k-1,e},{r\over 2}).$$
On this disc, we write
\ekv{ev.8}
{
G(\mu ;h)=a_{\nu (k-1)}(\mu )b_k(\mu ),
}
where
$b_k$ is \hol{}, and
\begin{eqnarray}\label{ev.9}
&C\ge \vert b_k(\widetilde{\mu }_{k-1,e})\vert \ge  1/(C\ln{1\over \langle
\widetilde{\mu }_{k-1,e}\rangle _h}),&\\
&\vert b_k(\mu )\le  \exp {\cal O}(1)[{\ln {1\over \langle \mu \rangle _h}
\over h}h{\ln\ln\over \ln} ({1\over \langle \mu \rangle _h})]=\exp {\cal
O}(\ln\ln {1\over \langle \mu \rangle _h}).&\nonumber
\end{eqnarray}
Using (\ref{ev.9}) and the elementary arguments recalled in the second part
of Section \ref{SectionZe}, we get
\ekv{ev.10}
{
{1\over 2\pi i}\int_{{\gamma _\vert}_{I_k}}{b_k'\over b_k}d\mu ={\cal
O}(1)\ln\ln {1\over \langle \mu \rangle _h}.
}

\par On the other hand,
\ekv{ev.11}
{
{1\over 2\pi i}\int_{{\gamma _\vert}_{I_k}}{a'_{\nu (k-1)}\over a_{\nu
(k-1)}}d\mu ={1\over 2\pi h}(\phi _{\nu (k-1)}(\widetilde{\mu
}_{k,s})-\phi _{\nu (k-1)}(\widetilde{\mu }_{k-1,e})).
}
Combining the real parts of (\ref{ev.7}), (\ref{ev.8}), (\ref{ev.10})
(\ref{ev.11}), we get
\begin{eqnarray}\label{ev.12}
&\Re {1\over 2\pi i}\int_\gamma {G'\over G}d\mu =&\\
&
\Re (\phi _{\nu (M)}(\widetilde{\mu }_{M,e})+\sum_{k=0}^{M-1}(\phi _{\nu
(k)}(\widetilde{\mu }_{k+1,s})-\phi _{\nu (k+1)}(\widetilde{\mu
}_{k+1,s}))-\phi _{\nu (0)}(\widetilde{\mu }_{0,s}))
&\nonumber\\
&
+{\cal O}(1)\ln\ln {1\over \langle \mu \rangle _h},\nonumber
&
\end{eqnarray}
with the remainder improving to ${\cal O}(1)$ if we do not encounter
$B_e$. Now, $\widetilde{\mu }_{M,e}=\mu _{M,e}$, $\widetilde{\mu
}_{0,s}=\mu _{0,s}$, and
$$\mu _{k,e}=\mu _{k+1,s}+{\cal O}(h{\ln\ln\over \ln}({1\over \langle \mu
\rangle _h})),$$
with the last remainder improving to ${\cal O}(h/\ln{1\over \langle \mu
\rangle _h})$, if we avoid $B_e$, and (\ref{ev.3}) follows.
\end{proof}

\par We end this section by some rough estimates on the location of the
skeleton and the corresponding distribution of \ev{}s for the reduced
\op{}s constructed in Sections \ref{SectionRed}, \ref{SectionExp}. Our
starting point is the reduced symbol in (\ref{Red.8}) and the corresponding
1-dimensional symbol
\ekv{ev.13}
{Q\left(\tau ,x,\xi,\epsilon ,{h^2\over \epsilon } ;h\right)=
\langle
q\rangle (\tau ,x,\xi )+{\cal O}(\epsilon )+{h^2\over \epsilon }p_2(\tau
,x,\xi )+h\widetilde{p}_1+h^2\widetilde{p}_2+...\ .}
Here we shall take $\tau $ real (and eventually  of the form
$h(k-{k_0\over 4})-{S_0\over 2\pi })$). If $z$ is the original spectral parameter,
we introduce the new
spectral parameter $w$, by
\ekv{ev.14}
{
z=g(\tau )+i\epsilon w,
}
and we will work under the assumption $h^2\ll \epsilon \ll h^{1/2}$.

\par Recall the \ml{} normal form for $Q$ near the branch point, given by
Proposition \ref{Prop1D2} and in particular (\ref{1D.31}):
\ekv{ev.15}
{ U^{-1}QU=K_{\epsilon ,h^2/\epsilon }(\tau ,I;h)+{\cal O}(e^{-{1\over
Ch}}),\ I={1\over 2}(x\circ hD_x+hD_x\circ x), }
where the leading symbol in $K_{\epsilon ,h^2/\epsilon }$ is
$k_{\epsilon ,h^2/\epsilon }(\tau ,\iota )$ with $\iota =x\xi $, given in Proposition
\ref{Prop1D1}. Correspondingly, we replace $w$ by the new spectral
parameter $\mu $, given by
\ekv{ev.16}
{
K_{\epsilon ,h^2/\epsilon }(\tau ,\mu ;h)=w.
}

\par We next estimate the location of the skeleton in the $\mu $-plane, and
start with the case $\vert \mu \vert \ge Ch$. Assume for simplicity that
we are in the case 1: $\Im \mu \ge -C\vert \Re \mu \vert $. We will only
be concerned with the horizontal part $S'$ of the skeleton. When $\Re \mu
\ge 0$, it is given by the curves $\Gamma _{3,4^+}=\Gamma _{1,2}$, $\Gamma
_{2,4^+}=\Gamma _{1,3}$ in (\ref{Sk.4}), where
\begin{eqnarray*}
X(\mu )&=&Y(\mu )+{\pi \over 2}\Re \mu ,\\
Y(\mu )&=& (\Re \mu ){\rm arg\,}({\mu \over i})-\Im \mu
+h\Re {\cal O}_-({h\over \mu }).
\end{eqnarray*}
Clearly $X(\mu )$ is \ufly{} Lipschitz continuous and for $\mu >0$, we get
$$
X(\mu )=h\Re {\cal O}_-({h\over \mu }).
$$
According to (\ref{as.4.1}), we have $\Re {\cal O}_-({h\over \mu
})={\cal O}(e^{-2\vert \mu \vert /h})$.

\par When $\epsilon =0$, $h^2/\epsilon =0$, we know that the leading part
of $Q$ in (\ref{ev.13}) is real-valued (assuming that $\langle q\rangle $
is real for simplicity), so it follows in this case that when $\mu $ is
real, then $\Im S_{j,k}={\cal O}(h)$. Since $S_{j,k}$ depends \hol{}ally
on $\mu $, we conclude that in general
\ekv{ev.17}
{
\Im S_{j,k}(\mu )={\cal O}(\epsilon +h^2/\epsilon ),\ \mu \in{\bf R}.
}
Now combine this with (\ref{Sk.4}), the estimate $X(\mu )={\cal
O}(he^{-2\pi \vert \mu \vert /h})$ and Proposition \ref{PropA3.1} to
conclude that in the region $\vert \mu \vert \ge Ch$, $\Re \mu \ge 0$, the
horizontal part $S'$
of the spectrum is given by the union of two curves of the form
\ekv{ev.18}
{
\Im \mu =f(\Re \mu ),
}
with $f'$
satisfying (\ref{A3.21}), and further,
\ekv{ev.19}
{
\vert f(x)\vert \le C\left(\epsilon +{h^2\over \epsilon }\right)\max \left({1\over \ln
{1\over \vert x\vert }},{1\over \ln {1\over (\epsilon +{h^2\over
\epsilon })}}\right).
}

\par In the left half-plane, we recall that $S'$ has a more complicated
structure. Assume, to fix the ideas, that $\Re \mu _A\le 0$. Then $S'$ is
the union of the curves (defined in (\ref{Sk.5})):
\begin{eqnarray*}
\Gamma _{1,3}:&& \Im\mu = \gamma _{1,3}(\Re \mu ) \Leftrightarrow (\Im
\mu )\ln{1\over \vert \mu \vert }=\Im S_{1,2}+(X-\pi \Re \mu )+\pi \Re
\mu \\
\Gamma _{3,4^-}:&& \Im\mu = \gamma _{3,4^-}(\Re \mu ) \Leftrightarrow (\Im
\mu )\ln{1\over \vert \mu \vert }=\Im S_{3,4}+(X-\pi \Re \mu )-\pi \Re
\mu ,
\end{eqnarray*}
in the region $\Re \mu _A\le \Re \mu \le 0$. Here $\gamma _{3,4^-}(\Re
\mu )\le \gamma _{1,3}(\Re \mu )$ and the two curves cross at $\mu _A$. In
the region $\Re \mu \le \Re \mu _A$ $S'$ is given by
\begin{eqnarray*}
\Gamma _{1,4^-}:&& \Im\mu = \gamma _{1,4^-}(\Re \mu ) \Leftrightarrow (\Im
\mu )\ln{1\over \vert \mu \vert }={1\over 2}(\Im S_{1,2}+\Im S_{3,4})+
(X-\pi \Re \mu ),
\end{eqnarray*}
and this curve also contains $\mu _A$.

\par When $\mu <0$, we have
$$X-\pi \Re\mu =h\Re {\cal O}_-({h\over \mu })=h{\cal O}(e^{-2\pi \vert \mu \vert /h}),$$
so again $f:=\gamma _{1,4^-}$ satisfies (\ref{ev.19}), while
\begin{eqnarray*}
\widetilde{\gamma }_{3,4^-}(\Re \mu ):=\Im S_{3,4}+(X-\pi \Re \mu ) &
\le & \gamma _{3,4^-}(\Re \mu )\le \gamma _{1,3}(\Re \mu ) \\
& \le & \Im S_{1,2}+(X-\pi \Re \mu )=:\widetilde{\gamma }_{1,3}(\Re
\mu ),
\end{eqnarray*}
for $\Re \mu _A\le \Re \mu \le 0$, where $f=\widetilde{\gamma
}_{1,3},\widetilde{\gamma }_{3,4^-}$ satisfy (\ref{ev.19}).

\par In the region $\vert \mu \vert \le Ch$ the horizontal part of the
spectrum is a union of curves  $\Gamma _{j,k}$ given in (\ref{sk.4}),
(\ref{sk.5}). Here the new \fu{} $X$ is \ufly{} Lipschitz and ${\cal
O}(h)$, so the skeleton is here contained in a region

\ekv{ev.20}
{
\vert \Im \mu \vert \le {\cal O}(1){(\epsilon +{h^2\over \epsilon })\over
\ln {1\over h}}.
}

\par The \it overall conclusion \rm is that the skeleton is contained in a
region
\ekv{ev.21}
{
\vert \Im \mu \vert \le {\cal O}(1)(\epsilon +{h^2\over \epsilon })\max
\left({1\over\ln {1\over \langle \Re \mu \rangle _h}},{1\over
\ln {1\over \epsilon +{h^2\over
\epsilon }}}\right),
}
where we recall that $\langle \Re \mu \rangle _h=(h^2+(\Re \mu )^2)^{1\over
2}.$

\par We end this section by establishing a simplified statement, to be
used in Theorem 1.1.  We shall simply remove a small rectangle
around $\mu =0$ where we have seen that the description of the spectrum is
more intricate.

\par Start by recalling the definition of $\mu _A$, $\mu _B$ prior to
(\ref{Sk.18}). For instance $\mu _A$ is the intersection of the curves
\par\noindent
$A$: $-2\pi \Re \mu =\Im S_{1,2}-\Im S_{3,4}$,
\par\noindent and
\par\noindent $\Gamma _{1,4^-}$: $\Im \mu =\gamma _{1,4^-}(\Re \mu )$,
where $f=\gamma _{1,4^-}$ satisfies (\ref{ev.19}). Using that $\Im
S_{1,2},\Im S_{3,4}={\cal O}(\epsilon +h^2/\epsilon +\vert \Im \mu
\vert )$, we get
$$
\Re \mu _A={\cal O}\left(\epsilon +{h^2\over \epsilon }+\vert \Im \mu
_A\vert \right),\ \Im \mu _A={\cal O}\left(\epsilon +{h^2\over \epsilon }\right)\max
\left({1\over \vert \ln \vert \Re \mu _A\vert \vert },{1\over \vert \ln
(\epsilon +{h^2\over \epsilon })\vert }\right),
$$
implying,
\ekv{ev.22}
{
\Re \mu _A={\cal O}(\epsilon +{h^2\over \epsilon }),\ \Im \mu _A={{\cal
O}(\epsilon +{h^2\over \epsilon })\over \vert \ln (\epsilon +{h^2\over
\epsilon })\vert }.
}
We have of course the same estimates for $\mu _B$.

\par Choose $C>0$ \sufly{} large so that the "black box"
\ekv{ev.23}
{
{\cal B}=[-a,a]+i[-b,b],\hbox{ with }a=C(\epsilon +{h^2\over \epsilon }),\
b=C{(\epsilon +{h^2\over \epsilon })\over \vert \ln (\epsilon +{h^2\over
\epsilon })\vert },
}
contains $\mu _A$, $\mu _B$. Then we have
\begin{prop}\label{PropEv.2}
The number of \ev{}s in ${\cal B}$ is ${\cal O}({\epsilon \over
h}+{h\over \epsilon })\vert \ln (\epsilon +{h^2\over \epsilon })\vert $.
The \ev{}s outside ${\cal B}$
are exponentially close to $\Gamma _{1,4^-}\cup \Gamma _{1,2}\cup \Gamma
_{1,3}$. More precisely introduce
\begin{eqnarray*}
E_{1,4^-}&=&\{ \mu \in \Gamma _{1,4^-}\setminus{\cal B};\, a_1+a_{4^-}=0,\
\Re \mu <0\},\\
E_{1,2}&=&\{ \mu \in \Gamma _{1,2}\setminus{\cal B};\, a_1+a_2=0,\
\Re \mu >0\},\\
E_{1,3}&=&\{ \mu \in \Gamma _{1,3}\setminus{\cal B};\, a_1+a_3=0,\
\Re \mu >0\}.
\end{eqnarray*}
Then there is a bijection $b$ (possibly after a slight modification of ${\cal
B}$) between the set of \ev{}s outside ${\cal B}$ and $E_{1,4^-}\cup
E_{1,2}\cup E_{1,3}$, such that $b(\mu )-\mu ={\cal O}(e^{-\pi \vert \Re
\mu \vert /h}h/\vert \ln \vert \mu \vert \vert )$.
\end{prop}
\begin{proof}
We may first notice that we can replace the index $4^-$ by $4$ without
changing the validity of the statement of the proposition, since
$a_4-a_{4^-}={\cal O}(e^{-2\pi \vert \Re \mu \vert /h})$, $\Re \mu \ll -h$.
In view of (\ref{ev.21}), we know that there are no \ev{}s outside ${\cal
B}$ with $\vert \Re \mu \vert \le a $ and the discussion in Section
\ref{SectionSk} then shows that the \ev{}s outside ${\cal B}$ have to be
exponentially close to  $\Gamma _{1,4}\cup \Gamma _{1,2}\cup \Gamma
_{1,3}$ and that there is a bijection $b$ as stated. To estimate the number
of \ev{}s inside ${\cal B}$, we simply apply Theorem \ref{ThEv.1}.
with $\gamma $ a rectangular contour containing ${\cal B}$ but contained
in $2{\cal B}$ and working directly with $a_4$ instead of $a_{4^\pm}$.
\end{proof}

\par Consider $E_{1,4^-}$ of the preceding proposition. In view of
(\ref{Sk.1.5}), it is given by the quantization condition
\ekv{ev.23.5} { S_{1,2}+S_{3,4}+2\mu (\ln (-\mu )-1)+{\pi h\over
2}+2hi{\cal O}_-({h\over \mu })=2\pi (k+{1\over 2})h,\ k\in{\bf Z}
} Here we recall from the beginning of Section \ref{SectionAs},
that the term ${\cal O}_-({h\over \mu })$ is $\sim C_1{h\over \mu
}+C_2({h\over \mu })^2+...$, as ${h\over \mu }\to 0$. We also know
that if $\alpha =(\epsilon ,h^2/\epsilon )$ denote the small
additional parameters in the \pb{}, then $S_{j,k}\sim
\sum_0^\infty S_{j,k}^\nu (\mu ,\alpha )h^\nu $, for $(j,k)=(1,2),
(3,4)$, where $S_{j,k}^\nu $ are smooth in $\alpha$ and \an{} in
$\tau$. Hence the condition \nr{ev.23.5} takes the form
\ekv{ev.24} { b_{1,4^-}(\mu ,\alpha ;h)=2\pi (k+{1\over 2})h, }
where \ekv{ev.25} { b_{1,4^-}(\mu ,\alpha ;h)\sim\sum_{\nu
=0}^\infty b_{1,4^-}^\nu (\mu ,\alpha )h^\nu , } in the space of
\bdd{} \hol{} \fu{}s defined in truncated sector:
$$
\Re \mu \le -Ch,\ \vert \Im \mu \vert \le {1\over C}(-\Re \mu ),
$$
with \ekv{ev.26} {b_{1,4^-}^0(\mu ,\alpha )-2\mu \ln(-\mu ),\
b_{1,4^-}^1\hbox{ \hol{} in a full neighborhood of }\mu =0,\
\alpha =0, } and \ekv{ev.27} { b_{1,4^-}^\nu (\mu ,\alpha )={\cal
O}(\mu ^{1-\nu }),\ \nu \ge 2. } Notice that the singularity
structure \nr{ev.25}--\nr{ev.27} of $b_{1,4^-}$ is essentially
unchanged if we replace $\mu $ by $\widetilde{\mu }=a(\mu ,\alpha
;h)(\mu +hd(\mu ,\alpha ;h))$, where $a,d$ are classical symbols
of order 0 in $h$ with \coef{}s that are \an{} near $\mu =0$,
$\alpha =0$ and with $a$ elliptic, $\Re a>0$, $\vert \Im a\vert
\ll \Re a$.

\par On the other hand, in the region $\Re \mu <-1/C$, $C\gg 0$, we know
(and that was done for instance in \cite{HiSj1}), that the \ev{}s
sit on a curve and are given by a Bohr--Sommerfeld condition
\ekv{ev.28} { \widetilde{b}(\mu ,\alpha ;h)=2\pi (k+{1\over 2})h,\
k\in{\bf Z}, } where $\widetilde{b}$ is a classical analytic
symbol of order 0: $\widetilde{b}\sim \sum_0^\infty
\widetilde{b}^\nu (\mu ,\alpha )h^\nu $, and where \ekv{ev.29} {
\widetilde{b}^0(\mu ,\alpha )=\int_{\gamma _{\rm ext }(\mu ,\alpha
)}\xi dx. } Here $\gamma _{\rm ext}(\mu ,\alpha )$ denotes a
closed loop in the energy surface $Q^0(\mu ,\alpha ,x,\xi )=w$
with $w$ and $\mu $ related by \nr{ev.16}, that can be obtained
from the real energy curve we get by taking $\mu $ real and
putting $\alpha =0$. Clearly $\widetilde{b}=b_{1,4^-}$, so our
discussion gives detailed description about how one can push the
standard WKB-construction to the limit $\vert \mu \vert \gg h$ in
the region $\Re \mu <0$.

\par The same discussion applies to $E_{1,2}$, $E_{1,3}$. We get the
conditions
\ekv{ev.30}
{
b_{1,2}(\mu ,\alpha ;h)=2\pi kh\hbox{ and }b_{1,3}(\mu ,\alpha ;h)=2\pi
kh\hbox{ respectively ,}}
where $b_{j,k}$, $(j,k)=(1,2), (1,3)$ are defined in the truncated
sector $\Re \mu \ge Ch$, $\vert \Im \mu \vert \le \Re \mu /C$, and
$b_{j,k}^\nu $ have the analogous properties to those in \nr{ev.26},
\nr{ev.27}, for $\nu \ge 1$, while the first part of \nr{ev.26} should be
replaced by the condition that
\ekv{ev.31}
{
b_{j,k}^0(\mu ,\alpha )-\mu \ln \mu \hbox{ is \hol{} near }\mu =0,\ \alpha
=0.
}
$b_{1,2}^0$ is the action along a closed loop inside the appropriate
complex energy curve, that can be obtained by deformation from the case
$\mu >0$, $\alpha =0$ where we take the left real component, close to the
left loop in the $\infty $-shaped set $K_{0,0}$. For $b_{3,4}^0$ we deform
from the right real component.

\section{Justification by means of a global Grushin \pb}\label{SectionGr}
\setcounter{equation}{0}

\it One dimensional Grushin \pb{}s. \rm We may assume here without loss of
generality, that $\langle q\rangle $ is real-valued. Then we know that $f$
in (\ref{1D.33}) satisfies
\ekv{Gr.1}
{f(w;h)\sim\sum_0^\infty  f_k(w;h),}
where $f_0$ is real-valued when $\epsilon ,h^2/\epsilon =0$, $\tau \in {\bf
R}$. Recall that $f$ and $f_k$ depend \an{}ally on $\tau $ and smoothly
on $\epsilon , h^2/\epsilon $.

\par Also recall that the spectrum of $Q$ is localized to the region
(\ref{1D.8}): $\vert \Im w\vert ={\cal O}(h+\epsilon +h^2/\epsilon )$ (as
follows also from the more refined estimate (\ref{ev.21}) and
in view of (\ref{1D.34}): $f(w;h)=-\mu $, it follows that $\mu $ is
localized to a domain of the same type.

\par We shall now introduce three different Grushin \pb{}s for $Q-w$ in
the spirit of \cite{HeSj2}, \cite{SjZw}. Let $\chi \in C_0^\infty
(({\rm neigh\,}(0,0),{\bf R}^2))$ be equal to one near (0,0). We
realize $\chi $ as an $h$-\pop{}, that we also denote by $\chi $,
using a Gaussian resolution of the identity
(see~\cite{MeSj},~\cite{Sj3}), so that our calculus errors will be
exponentially small rather than just ${\cal O}(h^\infty )$.
Assume, in order to fix the ideas, that the support of $\nabla
\chi $ is a thin annulus around (0,0), containing the points
$\alpha _j$, $j=1,..,4$ (see Section \ref{Section1D}). Recall the
definition of $g_j$ in (\ref{1D.38}). Let $g_j^*$ be the
corresponding \fu{}s defined for $Q^*-\overline{w}$, depending
anti-\hol{}ally on $w$. Define \ekv{Gr.2} { R_+^ju=(-1)^j(u\vert
{i\over h}[Q^*,\chi ]g_j^*)_{W_j}. } Here $W_j$ is a small
\neigh{} of $\alpha _j$ and $(u\vert v)_{W_j}=(\chi _ju\vert v) $
where $\chi _j\in C_0^\infty (W_j)$ is equal to 1 near $\alpha _j$
and we also let $\chi _j$ denote the corresponding Gaussian
quantization. We may normalize the choice of $g_j^*$ so that
\ekv{Gr.3} { R_+^jg_j=1. } Expanding the commutator, we see that
the definition of $R_+^ju$ (up to an exponentially small error) is
\indep{} of the choice of $\chi $, provided that $u$ is a \ml{}
null-solution of $Q-w$.

\par Assume for simplicity that $\chi $ is real-valued and that the
corresponding quantization is \sa{}. Put
\ekv{Gr.3.5}
{
R_-^ju_-=(-1)^ju_-{i\over h}[Q,\chi ]g_j.
}
\smallskip

\par Our first Grushin \pb{} is directly adapted to the derivation of the
quantization condition (\ref{1D.49}) in Section \ref{Section1D}. It is the
form
\ekv{Gr.4}
{
\cases{(Q-w)u+R_-u_-=v\cr R_+u=v_+,}
}
with
\ekv{Gr.5}
{
R_+=R_+^4: L^2_{\theta '}\to {\bf C},\quad R_-=R_-^1:{\bf C}\to L^2_{\theta '}.
}
Using Section \ref{Section1D} we see as in \cite{HeSj2, SjZw} that it is
well posed for $w$ in some fixed complex \neigh{} of $0$ with a solution of
the form
\ekv{Gr.6}{\cases{u=Ev+E_+v_+,\cr u_-=E_-v+E_{-+}v_+.}
}
Here we get
\ekv{Gr.7}
{
E_{-+}v_+=ihv_+(c_{2,3}e^{2\pi i(\theta _1+\theta _2)+{i\over
h}(S_{3,4}+S_{1,2})}+c_{2,4}e^{2\pi i\theta _2+{i\over
h}S_{1,2}}-c_{1,3}e^{2\pi i\theta _1+{i\over h}S_{3,4}}-c_{1,4}),
}
where the parenthesis is the same as in the quantization condition
(\ref{1D.49}). As usual, we read off the approximate \ev{}s as the zeros of
E$_{-+}$.
\def\Grpb{Grushin problem}
\par The drawback with this first \Grpb{} is  that the solution \op{}
(\ref{Gr.6}) will grow exponentially when $\mu >0$. This can be seen either
directly from the explicit formulae for $a_{j,k}$ and the slightly less
explicit expression for $c_{j,k}$, or from the fact that for $\mu >0$, we
have approximately a double well problem and with $R_+$ we prescribe the
solution $u$ in (\ref{Gr.4}) in one of the wells, and hence $u$ will in
general be exponentially large in the other well. Of course, we will have
to accept some exponential growth with a rate ${\cal O}(\vert \Im w\vert
+\epsilon +h^2/\epsilon )$ but certainly would like to avoid exponential
growth with a fixed rate when $\mu $ is real.

\par It seems impossible to cover a full \neigh{} of $w=0$ with a single
\Grpb{} whose solution \op{} does not exhibit exponential growth in some
region, so we shall use 2 \Grpb{}s where one will be nice roughly in the
upper half-plane and the other in the lower half plane.

\par The second \Grpb{} is designed to cover the region $\Im \mu \ge 0$
with some margin. It is of the form (\ref{Gr.4}) with
\ekv{Gr.8}
{ R_+u=(R_+^2u,R_+^4u),\ R_-u_-=R_-^1u_-^1+R_-^3u_-^3, } so that
$R_+=L^2_{\theta '}\to{\bf C}^2$, $R_-:{\bf C}^2\to L^2_{\theta '}$.  For
the corresponding model \pb{} for $P_0-\mu $ in Section \ref{SectionTr}, we
get (cf(\ref{tr.11}), (\ref{tr.12}):
\ekv{Gr.9}
{
\pmatrix{u_1\cr u_3}=U\pmatrix{u_2\cr u_4},
}
\ekv{Gr.10}
{ U=\pmatrix{u_{1,2} &u_{1,4}\cr u_{3,2} & u_{3,4}}={\Gamma ({1\over
2}-i{\mu \over h})h^{^{-{i\mu \over h}}}\over \sqrt{2\pi }}
\pmatrix{e^{{\pi \over 2}{\mu \over h}+i{\pi \over 4}} &e^{-{\pi \over
2}{\mu \over h}-i{\pi \over 4}}\cr e^{-{\pi \over 2}{\mu \over h}-i{\pi
\over 4}} &e^{{\pi \over 2}{\mu \over h}+i{\pi \over 4}}}.
}
This is basically the approach of \cite{HeSj2} and as there, we see by
direct calculation or by a more general normalization argument that $U(\mu )$
is unitary when $\mu $ is real. We also see that $U$ is \ufly{} \bdd{} in
any disc $D(0,Ch)$ for
\ekv{Gr.11}
{
\Im \mu \ge -{Ch\over \ln{1\over h}}.
}

\par For $\mu $ outside an angle around $-i[0,+\infty [$ with $\vert \mu
\vert \gg h$, we apply Stirling's formula (\ref{as.2_-}) and get
\ekv{Gr.12} { u_{j,k}=\exp \big( {i\mu \over h}-{i\mu \over h}\ln
(-i\mu )\pm ({\pi \over 2}{\mu \over h}+i{\pi \over 4})+{\cal
O}_-({h\over \mu })  \big) } with the $+$ sign valid for
$u_{1,2}$, $u_{3,4}$, and the $-$ sign for $u_{1,4}$, $u_{3,2}$.
As in Section \ref{SectionSk}, we get \ekv{Gr.13} { \vert u
_{j,k}\vert =\exp \big( -{1\over h}(1+\ln {1\over \vert \mu \vert
})\Im \mu +{\Re \mu \over h}({\rm arg\,}({\mu \over i})\pm {\pi
\over 2})+\Re {\cal O}_-({h\over \mu }) \big) . } We now also
assume that $\vert \mu \vert \le C^{-1}\ll 1$, so that $\ln (\vert
\mu \vert ^{-1})\gg 1$. We shall estimate the exponent from above.
When $\Re \mu \ge 0$, the worst exponent is the one with $+\pi /2$
in the middle term and we approximate
$$
{1\over h}\Re \mu ({\rm arg\,}({\mu \over i})+{\pi \over 2})={1\over
h}(\Re \mu )( {\rm arg\,}\mu) \sim {1\over h}{(\Re \mu) (\Im \mu )\over \vert \mu
\vert },
$$
which is dominated by the first term and hence
\ekv{Gr.14}
{
\vert u _{j,k}\vert \le \exp (-{1\over h}(\ln {1\over \vert \mu \vert
}+{\cal O}(1))\Im \mu ) ,
}
when $\Re \mu \ge 0$. When $\Re \mu \le 0$, the worst case is the one
with $-\pi /2$ in the middle term and we approximate
$$
{1\over h}\Re \mu ({\rm arg\,}({\mu \over i})-{\pi \over 2})={-\Re \mu \over
h}(\pi - {\rm arg\,}\mu ) \sim {1\over h}{\vert \Re \mu\vert  (\Im \mu )\over \vert \mu
\vert },
$$
leading to (\ref{Gr.14}) also in this case. We conclude that $U$ is \bdd{}
in a domain of the form
\ekv{Gr.15}
{
\{\mu ;\, \Im \mu \ge -{Ch\over \ln{1\over h+\vert \mu\vert  }}\} .
}
As a consequence, we get
\begin{prop}\label{PropGr1}
The \pb{} {\rm (\ref{Gr.4})} with $R_{\pm}$ given by {\rm
(\ref{Gr.8})} is \ml{}ly well-posed (with errors ${\cal
O}(e^{-1/(Ch)}))$ for $\vert \mu \vert $ small with $\mu $ in
$D(0,{h\over 4})$ or away from a small conic \neigh{} of $i{\bf
R}_-$. If we write the solution as in {\rm (\ref{Gr.6})}, then for
$\mu $ as in {\rm (\ref{Gr.15})}, we have \ekv{Gr.16} {h\Vert
E\Vert ,\,\Vert E_\pm\Vert ,\,\vert E_{-+}\vert ={\cal O}(1)\exp
{C\over h}(\vert \Im \mu \vert +\epsilon +{h^2\over \epsilon }).}
\end{prop}

Let us also compute $E_{-+}$. Near the branching point, we recall that we
have the relation (\ref{1D.41}) for null-solutions to $P-z$, equal to
$u_jg_j$ near $\alpha _j$. To determine $E_{-+}$, we consider (\ref{Gr.4})
with $v=0$, so that $u_j=v_j^+$ for $j=2,4$. We then want to express
$u_j$, $j=1,3$ in terms of $u_2,u_4$. We can do this using (\ref{Gr.9}),
(\ref{Gr.10}) above and redo some of the work in Section \ref{Section1D},
but it is easier to to use the work already done and "solve"
(\ref{1D.41}). We get
\ekv{Gr.17}
{
\pmatrix{u_1\cr u_3}=\pmatrix{{c_{1,3}\over
c_{2,3}}&c_{1,4}-{c_{1,3}c_{2,4}\over c_{2,3}}\cr {1\over c_{2,3}}
&-{c_{2,4}\over c_{2,3}}}\pmatrix{u_2\cr u_4}.
}
Notice that $u_1^-$, $u_3^-$ in our \Grpb{} (with $v=0$) are the
discontinuities we obtain at $\alpha _1$, $\alpha _3$ when trying to extend
a null-solution near (0,0) with prescribed $u=v_j^+g_j$ near $\alpha
_j$, $j=2,4$, to a global null-solution near $K_{0,0}$. We get
\ekv{Gr.18}
{
\pmatrix{u_1^-\cr u_3^-}=E_{-+}\pmatrix{v_2^+\cr v_4^+}=\pmatrix{
e^{2\pi i\widetilde{\theta }_2}-{c_{1,3}\over c_{2,3}}
& {c_{1,3}c_{2,4}\over c_{2,3}}-c_{1,4}
\cr
-{1\over c_{2,3}}
& e^{2\pi i\widetilde{\theta }_1}+{c_{2,4}\over c_{2,3}}
}\pmatrix{v_2^+\cr v_4^+},
}
where
\ekv{Gr.19}
{
2\pi \widetilde{\theta }_1=2\pi \theta _1+{1\over h}S_{3,4},\ 2\pi
\widetilde{\theta }_2=2\pi \theta _2+{1\over h}S_{1,2}.
}
It follows that
\begin{eqnarray}\label{Gr.20}
\det E_{-+}&=&{1\over c_{2,3}}(e^{2\pi i(\widetilde{\theta
}_1+\widetilde{\theta
}_2)}c_{2,3}+c_{2,4}e^{2\pi i\widetilde{\theta }_2}-c_{1,3}e^{2\pi
i\widetilde{\theta }_1}-c_{1,4})\\
&=&e^{2\pi i(\widetilde{\theta }_1+\widetilde{\theta
}_2)}+{c_{2,4}\over c_{2,3}}e^{2\pi i\widetilde{\theta }_2}-{c_{1,3}\over c_{2,3}}e^{2\pi
i\widetilde{\theta }_1}-{c_{1,4}\over c_{2,3}}.\nonumber
\end{eqnarray}
From the middle expression, we see that this determinant is equal to the
expression in the quantization condition (\ref{1D.49}) times a
non-vanishing factor.
\smallskip

\par The third \pb{} is designed to cover the region $\Im \mu \le 0$. It is
of the form (\ref{Gr.4}) with
\ekv{Gr.21}
{R_+u=(R_+^1 u,R_+^3u),\quad R_-u_-=R_-^2u_-^2+R_-^4u_-^4.}
For the corresponding model \pb{} for $P_0-\mu $ in Section \ref{SectionTr}
we get (cf (\ref{tr.11}), (\ref{tr.12})):
\ekv{Gr.22}
{
\pmatrix{u_2\cr u_4}=V\pmatrix{u_1\cr u_3}=U^{-1}\pmatrix{u_1\cr u_3},
}
with $U$ as in (\ref{Gr.9}), unitary for real $\mu $, so that
\ekv{Gr.23}
{
V(\mu )=U(\overline{\mu })^*.
}
Thus from (\ref{Gr.14}) we see that the matrix elements $v_{j,k}$ satisfies
\ekv{Gr.24}
{
\vert v_{j,k}\vert \le \exp {1\over h}(\ln {1\over \vert \mu \vert }+{\cal
O}(1))\Im \mu ,
}
and $V$ is \bdd{} in a domain of the form
\ekv{Gr.25}
{
\{\mu ;\, \Im \mu \le {Ch\over \ln{1\over h+\vert \mu\vert  }}\} .
}
\begin{prop}\label{PropGr2}
The \pb{} {\rm (\ref{Gr.4})} with $R_{\pm}$ given by {\rm
(\ref{Gr.21})} is \ml{}ly well-posed (with errors ${\cal
O}(e^{-1/(Ch)})$) for $\vert \mu \vert $ small with $\mu $ in
$D(0,{h\over 4})$ or away from a small conic \neigh{} of $i{\bf
R}_+$. If we write the solution as in {\rm (\ref{Gr.6})}, then for
$\mu $ as in (\ref{Gr.25}), we have the estimate {\rm
(\ref{Gr.16})}.
\end{prop}

\par We now compute the corresponding $E_{-+}$ so we put $v=0$ in
(\ref{Gr.4}) and repeat the arguments above. Near the branching point $u$
is a null-solution, $=u_jg_j$ near $\alpha _j$, now with $v_j^+=u_j$,
$j=1,3$. By "solving" (\ref{1D.41}), we express $u_2,u_4$ in terms of
$u_1,u_3$ and find
\ekv{Gr.26}
{
\pmatrix{u_2\cr u_4}=\pmatrix{\displaystyle {c_{2,4}\over c_{1,4}} & \displaystyle
c_{2,3}-{c_{1,3}c_{2,4}\over c_{1,4}}\cr\displaystyle {1\over c_{1,4}}
&\displaystyle -{c_{1,3}\over c_{1,4}}}\pmatrix{v_1^+\cr v_3^+}.  }
We then get
\ekv{Gr.27}
{
\pmatrix{u_2\cr u_4}= E_{-+}\pmatrix{v_1^+\cr v_3^+}
=\pmatrix{\displaystyle {c_{2,4}\over c_{1,4}}-e^{-2\pi i\widetilde{\theta
}_2}&\displaystyle c_{2,3}-{c_{1,3}c_{2,4}\over c_{1,4}} \cr
\displaystyle -{1\over c_{1,4}}
&\displaystyle {c_{1,3}\over c_{2,4}}+e^{-2\pi i\widetilde{\theta }_1}
}
\pmatrix{v_1^+\cr v_3^+},
}
and
\begin{eqnarray}\label{Gr.28}
\det E_{-+}&=& {e^{-2\pi i(\widetilde{\theta }_1+\widetilde{\theta
}_2)}\over c_{1,4}}(c_{2,3}e^{2\pi i(\widetilde{\theta }_1+\widetilde{\theta
}_2)}+ c_{2,4}e^{2\pi i\widetilde{\theta }_2}-c_{1,3}e^{2\pi
i\widetilde{\theta }_1}-c_{1,4})
\\
&=& {c_{2,3}\over c_{1,4}}+{c_{2,4}\over c_{1,4}}e^{-2\pi
i\widetilde{\theta }_1}-{c_{1,3}\over c_{1,4}}e^{-2\pi
i\widetilde{\theta }_2}-e^{-2\pi i(\widetilde{\theta }_1+\widetilde{\theta
}_2)}.
\nonumber
\end{eqnarray}

\smallskip
\par\noindent \it The global Grushin \pb{}. \rm  We first explain what the
natural range will be for $\epsilon $. Our global Grushin \pb{} will be
built from a direct sum of \pb{}s for the \op{}s $\widehat{P}_\epsilon -z$
in (\ref{Red.5}). These \op{}s can be written
\ekv{Gr.29}
{
g(h(k-{k_0\over 4})-{S_0\over 2\pi })+i\epsilon Q(h(k-{k_0\over
4}),x,hD_x,\epsilon ,{h^2\over \epsilon };h),
}
where $Q$
is the \op{} appearing earlier in this section and in Section
\ref{Section1D}.

\par Clearly, we add conditions "$R_{\pm}$" only for such $k$ for which
$$
w_k={z-g(h(k-{k_0\over 4})-{S_0\over h})\over i\epsilon }
$$
is close the spectrum of $Q$, i.e.~ for which
\ekv{Gr.30}
{
\vert \Im w_k \vert \le {\cal O}(\epsilon +{h^2\over \epsilon }),
}
(cf. (\ref{1D.8})), assuming for simplicity that $\langle q\rangle $ is
real-valued. Then according to (\ref{Gr.16}) we can expect that our \Grpb{}
will have an inverse ${\cal E}$ with norm $\Vert {\cal E}\Vert ={\cal
O}(1)e^{C({\epsilon \over h}+{h\over \epsilon })}$.

\par Now we cannot expect to have a complete decomposition into a direct sum
of \op{}s \nr{Gr.29}, but Proposition \ref{PropExp3} shows that so
is possible up to an error ${\cal O}(1)\exp(-1/(C(\epsilon +h)))$.
Thus in order to absorb the error by a standard perturbation
argument, we need
$$
e^{C({\epsilon \over h}+{h\over \epsilon })-{1\over C(\epsilon +h)}}\ll 1,
$$
which would follow from
\ekv{Gr.30.5}
{
{1\over \epsilon +h}\gg {\epsilon \over h}+{h\over \epsilon },
}
or equivalently
$$
\epsilon ^3+h^3\ll \epsilon h,
$$
Here we already know that $h^3\ll \epsilon h$, since we assume
$h^2/\epsilon \ll 1$, so the new constraint is $\epsilon ^2\ll h$. From
now on we work in the range
\ekv{Gr.31}
{
h^2\ll \epsilon \ll h^{1/2}.
}

\par Since $g$ is real-valued with $g'\ne 0$, it follows from \nr{1D.8}
(or \nr{1D.9} in the general case, when $\langle q\rangle $ is not assumed
to be real) that the \op{}s \nr{Gr.29} have disjoint spectra. When \nr{Gr.30}
is never satisfied for any $k$, it is straight forward to see that $z$ is
not in the spectrum of our original \op{} $P_\epsilon $. Assume now that
\nr{Gr.30} holds for (at most) one $k=\widetilde{k}$. Let
$R_\pm^{(\widetilde{k})}$ be the corresponding \op{}s $R_{\pm}$ defined
earlier in this section. Using the notation of Proposition \ref{PropRed1},
we define
\begin{eqnarray}
R_+u&=&R_+^{(\widetilde{k})}((e^{{i\over h}A}U^{-1}e^{-{\epsilon \over
h}G}u\vert e_{\widetilde{k}-{k_0\over 4}-{S_0\over 2\pi
h}})_{L^2(S^1)})\label{Gr.32}\\
R_-u_-&=& e^{{\epsilon \over h}G}Ue^{-{i\over
h}A}(e_{\widetilde{k}-{k_0\over 4}-{S_0\over 2\pi h}}\otimes R_-^{(\widetilde{k})}u_-),
\nonumber\end{eqnarray}
for $u\in L^2(M)$, $u_-\in{\bf C}$. Here $e_k(t)=e^{ikt}$.

\par Repeating the arguments from \cite{HiSj1, HiSj2} we see that we get a
well-posed \pb{} with
\ekv{Gr.33}
{
E_{-+}(z)=E_{-+}^{(\widetilde{k})}(z)+{\cal O}(e^{-1\over C(\epsilon +h)}).
}
Here $E_{-+}^{(\widetilde{k})}$ is the "$E_{-+}$" of the approximate
1-dimensional \Grpb{} treated in the appropriate one of the Propositions
\ref{PropGr1}, \ref{PropGr2}, depending on the sign of the corresponding
parameter $\Im \mu $.

\par In the expressions \nr{Gr.20}, \nr{Gr.28}, we have the term $e^{\pm
2\pi i(\widetilde{\theta }_1+\widetilde{\theta }_2)}$, where
$\widetilde{\theta }_j$ are given by \nr{Gr.19}. We see that $\Im
\widetilde{\theta }_j={\cal O}({\epsilon\over h}+{h\over \epsilon })$, so
\ekv{Gr.34}
{
e^{\pm 2\pi i(\widetilde{\theta }_1+\widetilde{\theta }_2)}=e^{{\cal
O}({\epsilon \over h}+{h\over \epsilon })}.
}
We conclude from this and \nr{Gr.30.5} that the remainder in \nr{Gr.33}
is ${\cal O}(e^{-{1\over C(\epsilon +h)}})$ times the dominating term in the
expression for $E_{-+}^{(\widetilde{k})}$ in \nr{Gr.20}, or \nr{Gr.28}
respectively. {\it This implies that if we pass to the $\mu $-variable (for
$k=\widetilde{k}$) and define the skeleton as in Sections \ref{SectionSk},
\ref{Sectionsk} and the corresponding body $B$ as in the beginning of Section
\ref{SectionEv}, then the zeros of $\det E_{-+}$ are confined to $B$ and
Theorem \ref{ThEv.1} still applies to give the number of \ev{}s (in the
$\mu $-plane) inside an admissible curve.}

\section{Improved parameter range for barrier top resonances in the
resonant case}\label{SectionPara}
\setcounter{equation}{0}

\par We start the discussion in this section with the following
general observation. Let $P_\epsilon $ be a smooth family of operators, satisfying all the
assumptions of the introduction, and in particular (\ref{H1}). In
Proposition 2.1 we have seen that the operator $P_{\eps=0}$ can be reduced by
successive averaging procedures to \ekv{Para.1} {
\widehat{P}_{\eps=0}=g(hD_t)+hp_1(hD_t,x,hD_x)+h^2p_2(hD_t,x,hD_x)+....,\
t\in S^1,\,\tau,x,\xi \approx 0. }

\begin{prop} Assume that the subprincipal symbol of $P_{\eps=0}$ vanishes and that
the spectrum of $P_{\eps=0}$ clusters into bands of size $\le {\cal O}(1)h^{N_0}$,
for some integer $N_0\ge 2$.  Then $p_j(\tau ,x,\xi )=p_j(\tau )$ are
\indep{} of $(x,\xi )$ for $1\le j\le N_0-1$ in {\rm (\ref{Para.1})}.
\end{prop}

\begin{proof} \rm Since the subprincipal symbol vanishes, we already know
that $p_1=0$.  Suppose that the conclusion of the proposition does not hold and
let $N_1\in \{ 2,3,...,N_0-1\}$ be the smallest $N$ with $p_N(\tau
_0,x,\xi )$ non-constant for some $\tau _0\approx 0$.
Take a \fy{} of Gaussian quasimodes $e_\alpha (x)$, $\alpha =(\alpha
_x,\alpha _\xi )\in {\rm neigh\,} (0,{\bf R}^2)$ with $$
\Vert e_\alpha \Vert =1,\ p_{N_1}^w(\tau _0,x,hD_x)e_\alpha =p_{N_1}(\tau
_0,\alpha )(\alpha )e_\alpha
+{\cal O}(h^{1/2})\hbox{ in }L^2.
$$
See, e.g.~\cite{CoRaRo} for the standard construction of such
quasimodes. Then put
$$
f_{\alpha ,h}=(2\pi
)^{-1/2}e^{\frac{i}{h}\left(h\left(k-\frac{k_0}{4}\right)-\frac{S_0}{2\pi}\right)t}e_\alpha
(x),
$$
with $k=k(h)$ such that $h(k(h)-k_0/4)-S_0/2\pi \to
\tau _0$, so that
\ekv{Para.2}
{
\widetilde{P}_0^wf_{\alpha
,k}=g\left(h\left(k-\frac{k_0}{4}\right)-\frac{S_0}{2\pi}\right)+
h^{N_1}p_{N_1}\left(h\left(k-\frac{k_0}{4}\right)-\frac{S_0}{2\pi},\alpha
\right) f_{\alpha ,k}+{\cal O}(h^{N_1+{1\over 2}}) \hbox{ in } L^2.
}
Hence, since we are dealing with \sa{} \op{}s,
\ekv{Para.3}
{
{\rm
dist\,}\left(g\left(h\left(k-\frac{k_0}{4}\right)-\frac{S_0}{2\pi}\right)+
h^{N_1}p_{N_1}\left(h\left(k-\frac{k_0}{4}\right)-\frac{S_0}{2\pi},\alpha
\right),\sigma (P_0)\right)\le {\cal O}(h^{N_1+{1\over
2}}),
}
and varying $\alpha $, so that the values $p_{N_1}(h(k-\frac{k_0}{4})-\frac{S_0}{2\pi},\alpha )$ fill up a whole
interval, we get a contradiction to the clustering
assumption.\end{proof}

\Remark. Proposition 12.1 remains to hold in the case described in
Section 4 of~\cite{HiSj1}, where the operator $P_{\eps=0}$ is
conjugated into a normal form in a \neigh{} of a Lagrangian torus,
rather than near a closed $H_p$--trajectory.

From now on we shall assume that $P_{\eps=0}$ satisfies the
assumptions of Proposition 12.1. Let us now switch on $\eps$. An
application of Proposition 2.1 together with Proposition 12.1 then
shows that microlocally, near a closed $H_p$--trajectory,
$P_{\eps}$ can be reduced to the form \ekv{Para.5} {
\widehat{P}_\epsilon =g(hD_t)+\epsilon \left(i\langle q\rangle
(hD_t,x,hD_x)+{\cal O}(\epsilon )+{\cal O}\left({h^{N_0}\over
\epsilon }\right)+h\widetilde{p}_1+...\right). } It follows
therefore that in the results of~\cite{HiSj1} (Theorems 6.4 and
6.7 there), we can replace the exponent 2 by the exponent $N_0$ in
the parameter range for $\eps$. Thus for the study of the spectrum
of $P_{\eps}$ in a region where $\abs{\Re z}<1/{\cal O}(1)$ and
$\abs{\Im z/\epsilon -F_0}\leq 1/{\cal O}(1)$, when $F_0$ is a
non-critical value or a \nondeg{} maximum or minimum of $\Re
\langle{q}\rangle$ along $p^{-1}(0)$, it suffices to assume that
\ekv{Para.6} { h^{N_0}\ll \epsilon \le h^\delta , } for some
$\delta >0$. In the case when $F_0$ is a saddle point value of
$\Re \langle{q}\rangle$, from Theorem 1.1 we get the condition
\ekv{Para.7} { h^{N_0}\ll \epsilon \ll h^{1/2}. } Indeed, in the
latter case we still have Proposition \ref{PropExp3} and the
decoupling condition analogous to \nr{Gr.30.5} becomes $\epsilon
^3+h^{N_0+1}\ll \epsilon h$, which is fulfilled by \nr{Para.7}.

We shall now apply these observations to improve the
result of Proposition 7.1 of~\cite{HiSj1}, giving a description of
the individual barrier top resonances of the semiclassical
Schr\"odinger operator in the resonant case. Before doing so, and also
for the future use in Section 13, we shall first briefly recall the
general setup in Section 7 of~\cite{HiSj1}, as well as in Section 5
of~\cite{HiSj2}.

As in~\cite{HiSj1,HiSj2}, let us consider
\ekv{Bt.00}
{
P=-h^2\Delta +V(x),\ P(x,\xi )=\xi ^2+V(x),\ (x,\xi )\in T^*{\bf R}^2,
}
where $V$ satisfies the general assumptions of Section 7 in \cite{HiSj1},
allowing us to define the resonances of $P$ in the lower half-plane inside
some fixed \neigh{} of $E_0>0$, where $V(0)=E_0$, $V'(0)=0$,
$V''(0)<0$. As in \cite{HiSj1,HiSj2}, we assume that $\{ (0,0)\}$ is
the only trapped $H_P$-trajectory in $P^{-1}(E_0)$. After a linear
symplectic change of coordinates, we may write
\ekv{Bt.01}
{
P(x,\xi )-E_0=\sum_{j=1}^2 {\lambda _j\over 2}(\xi
_j^2-x_j^2)+p_3(x)+p_4(x)+...,\ (x,\xi )\to 0,
}
where $\lambda _j>0$ and $p_j(x)$ is a homogeneous \pol{} of degree{}
$j\ge 3$. Recall further from \cite{HiSj1} that the study of \res{}s
of $P$ near $E_0$ can be reduced to an \ev{} \pb{} for $P$ after
applying some variant of
the method of complex scaling, and that near $x=0$ this simply amounts to
working in the new real coordinates $(\widetilde{x},\widetilde{\xi })$, given
by $x=e^{i\pi /4}\widetilde{x}$, $\xi =e^{-i\pi /4}\widetilde{\xi }$.

Performing the scaling and dropping the tildes from the notation, we
see that the problem reduces to studying the \ev{}s close to $0$ of the \op{} $i(P-E_0)$, now
elliptic outside a small \neigh{} of $(0,0)$, with symbol
\ekv{Bt.03}
{P(x,\xi )=p(x,\xi )+ie^{3i\pi \over 4}p_3(x)+ie^{i\pi
}p_4(x)+ie^{5i\pi \over 4}p_5(x)+...,}
where
\ekv{Bt.1}
{
p(x,\xi)=\sum_{j=1}^2 {\lambda _j \over 2}(x_j^2+\xi _j^2).
}
Here we continue to write $P$ to denote
the scaled operator.

We assume that $\lambda _j>0$ in (\ref{Bt.03}) are rationally dependent,
\ekv{Bt.2}
{
\exists k^0=(k_1^0,k_2^0)\in {\bf Z}^2\setminus \{ 0\},\ k_1^0\lambda
_1+k_2^0\lambda _2=0,}
which implies that the $H_{p}$-flow is periodic.

As in~\cite{HiSj1},~\cite{HiSj2}, we are interested in eigenvalues $E$
of $P$ with $\abs{E}\sim \eps^2$, $0<\eps \ll 1$. After a
rescaling $x=\eps \tilde{x}$, and dropping the tildes over the new
variables, we get an operator $P_{\eps}=\frac{1}{\eps^2}P$ that we
view as an $\widetilde{h}$-pseudodifferential operator with
the symbol
$$
P_{\eps}(x,\xi)=\frac{1}{\eps^2}
P(\eps (x,\xi))=p(x,\xi)+i\eps e^{3\pi i/4}
p_3(x)-i\eps^2 p_4(x)+{\cal O}(\eps^{3}).
$$
Here $\widetilde{h}=h/\eps^2$. Now the spectrum of $P_{\eps=0}$ is
that of the harmonic oscillator, and hence it clusters into sets
of diameter $0$ and separation of order $h$. An application of
Proposition 12.1 shows that all the $p_j$ in (\ref{Para.1}) are
constant. Moreover, since in this case all the \ev{}s depend
linearly on $h$, we see from the proof that the $p_j$ have to
vanish. It follows from (\ref{Para.6}) that in the zone
corresponding to non-critical values $F_0$ or non-degenerate
maxima or minima, the range of energies that we get is \begeq
\label{range1} h^{2N_0/(1+2N_0)}\ll \abs{E-E_0}\leq
h^{2\delta/(1+2\delta)},
\endeq
for all $N_0=2,3,\ldots$ and all $\delta>0$. When $F_0$
corresponds to a branching level, we get from (\ref{Para.7})
\begeq \label{range2} h^{2N_0/(1+2N_0)}\leq \abs{E-E_0}\leq
h^{1/2}. \endeq

We summarize the discussion above in the following proposition, which
is an improvement of Proposition 7.3 in~\cite{HiSj1}. Clearly, in a
similar fashion, we also obtain an improvement of Theorem 5.1
in~\cite{HiSj2}.

\begin{prop}
Assume that the principal symbol $P(x,\xi)$ in {\rm (\ref{Bt.00})}
has an asymptotic expansion {\rm (\ref{Bt.01})}, and assume that
{\rm (\ref{Bt.2})} holds. Assume furthermore that the function
$\langle{p_3}\rangle$, defined as the average of $p_3$ along the
Hamilton flow of $p$ in {\rm (\ref{Bt.1})} does not vanish
identically. Then the resonances of the operator $P$ in the domain
\begeq \label{Para5} \left\{ z\in \comp; h^{2N_0/(1+2N_0)}\ll
\abs{z-E_0}\leq h^{\delta}\right\} \backslash \bigcup \left\{z\in
\comp; \abs{\Re z-E_0-A\abs{\Im z}^{3/2}}<\eta \abs{\Im
z}^{3/2}\right\},
\endeq
where $\eta>0$, $\delta>0$, and $N_0=2,3,\ldots$ are arbitrary but fixed, are
given by
\begeq
\sim E_0 - i\left(h(k_1-\alpha/4)+\eps^3
\sum_{j=0}^{\infty} h^j \eps^{-2j}
r_j\left(\frac{h}{\eps^2}\left(k-\frac{k_0}{4}\right)-\frac{S}{2\pi},\eps,
\frac{h^{N_0}}{\eps^{1+2N_0}}\right)\right),
\endeq
with
$$
r_0=ie^{3\pi i/4}\langle{p_3}\rangle(\xi)+{\cal
O}\left(\eps+\frac{h^{N_0}}{\eps^{1+2N_0}}\right),\quad r_j={\cal
O}(1),\quad j\geq 1.
$$
We have $k=(k_1,k_2)\in \z^2$, $S=(S_1,S_2)$ with $S_1=2\pi$, and
$\alpha=(\alpha_1,\alpha_2)\in \z^2$ is fixed, and we choose $\eps>0$
with $\abs{E-E_0}\sim \eps^2$.
The union in {\rm (\ref{Para5})} is taken over the set of critical values of
$\langle{p_3}\rangle$, restricted to $p^{-1}(1)$, with $A$ varying
over this set.
\end{prop}

\Remark. If $\langle{p_3}\rangle$ restricted to $p^{-1}(1)$ has
precisely one non-degenerate saddle point with the critical value
$A$, then the results of the present paper apply and give a
description of the individual resonances in a half-cubic
neighborhood of the curve $\Re z=E_0+A\abs{\Im z}^{3/2}$. In the
following section, we shall consider explicit examples of
homogeneous polynomials for which the assumptions of Theorem 1.1
are satisfied.

\section{Examples in the barrier top case}\label{SectionBt}
\setcounter{equation}{0}

This section is devoted to a study of examples of potentials of the
Schr\"odinger operator (\ref{Bt.00}) to which Theorem 1.1 is applicable.

Let us recall from (\ref{Bt.03}) that we are interested in eigenvalues
close to $0$ of the operator $P$, elliptic outside a small \neigh{} of
$(0,0)$ with symbol
\begeq
\label{bt1}
P(x,\xi)=p(x,\xi)+ie^{\frac{3\pi
i}{4}}p_3(x)+ie^{i\pi} p_4(x)+\ldots,
\endeq
where the harmonic oscillator $p(x,\xi)$ has been defined in
(\ref{Bt.1}). As before, we make the resonant assumption
(\ref{Bt.2}).

\par Consider first a general perturbation of $p$ of
the form of a linear combination of terms $x^\alpha \xi ^\beta $ with
$\vert \alpha \vert +\vert \beta \vert =m$, for some $m\in\{ 3,4,5,...\}$.
Recall from \cite{HiSj1}
how to compute the corresponding trajectory average $\langle x^\alpha \xi
^\beta \rangle $: Basically we use action-angle coordinates, but to start
with, we can do things a little easier by introducing
\ekv{Bt.3}
{
z_j=x_j+i\xi _j\in{\bf C},
}
and notice that along a $H_p$-trajectory we get in the $z_1,z_2$
coordinates:
\ekv{Bt.4}
{
z_j(t)=e^{-i\lambda _jt}z_j(0).
}
Then write $x_j(t)=\Re z_j(t)$, $\xi _j(t)=\Im z_j(t)$, so that
\begin{eqnarray}
x(t)^\alpha \xi (t)^\beta &=&\prod_{j=1}^2 ((\Re z_j(t))^{\alpha _j}(\Im
z_j(t))^{\beta _j})\label{Bt.5}
\\&=&{1\over 2^{\vert \alpha \vert +\vert \beta \vert }i^{\vert \beta \vert }
}\prod_{j=1}^2 ((z_j(0)e^{-i\lambda _jt}+\overline{z_j}(0) e^{i\lambda
_jt})^{\alpha _j}(z_j(0)e^{-i\lambda _jt}-\overline{z_j}(0) e^{i\lambda
_jt})^{\beta _j})\nonumber\end{eqnarray} Then expand the product by means
of the binomial theorem.  The time average is equal to the time-\indep{}
term and since this average is constant along each trajectory we shall
replace the symbols $z_j(0)$ simply by $z_j$.

\par In this section we consider the case when $\lambda _1=\lambda _2=1$,
$m=4$, $\beta =0$. (In \cite{HiSj1} we noticed that in this case the average will
vanish when $m=3$ and in \cite{HiSj2} we made a more refined study of that
case taking into account one more term in the perturbative expansion).
This means that we take $p_3(x)=0$, and for simplicity we also assume
that $p_m=0$ for all odd $m$ in (\ref{bt1}), so that we
can concentrate on the perturbation $-ip_4$ in
(\ref{bt1}). Performing a rescaling as described in the previous
section, with $\eps$ replaced by $\eps^{1/2}$, (i.e. setting
$x=\eps^{1/2} \tilde{x}$ rather than $x=\eps \tilde{x}$), we get
\ekv{Bt.5.5}
{
p(x,\xi)-i\epsilon p_4(x)-\epsilon^2 p_6(x)+
i\epsilon ^3 p_8(x)+\epsilon ^4 p_{10}(x)+...,
}
and here, as before, we choose $\epsilon$ of the same order of magnitude as the
modulus of the \ev{}s for the \op{} $P(x,hD)$, that we want to study.

\par Now we continue the calculations of trajectory averages using
(\ref{Bt.5}). We have
\begin{eqnarray}\label{Bt.6}
\langle x_1^4\rangle &=&{1\over 2^4}\langle
(z_1+\overline{z}_1)^4\rangle ={1\over 2^4}\langle
z_1^4+4z_1^3\overline{z}_1+6z_1^2\overline{z}_1^2+4z_1\overline{z}_1^3+
\overline{z}_1^4\rangle
\\
&=&{6\over 16}\langle z_1^2\overline{z}_1^2\rangle ={3\over 8}\vert
z_1\vert ^4 \ (={3\over 8}(x_1^2+\xi _1^2)^2).
\nonumber
\end{eqnarray}
In the same way, we get
\ekv{Bt.7}
{
\langle x_2^4\rangle={3\over 8}\vert z_2\vert ^4.
}

Next look at the averages of mixed terms:
\begin{eqnarray}\label{Bt.8}
\langle x_1^3x_2\rangle &=&{1\over 2^4}\langle
(z_1+\overline{z}_1)^3(z_2+\overline{z}_2)\rangle ={1\over 2^4}\langle
(z_1^3
+3z_1^2\overline{z}_1+3z_1\overline{z}_1^2+\overline{z}_1^3)(z_2+\overline{z}_2)\rangle \\
&=&{3\over 16}(\vert z_1\vert ^2\overline{z}_1z_2+\vert z_1\vert
^2z_1\overline{z}_2)={3\over 8}\vert z_1\vert ^2\Re (z_1\overline{z}_2).
\nonumber
\end{eqnarray}
\ekv{Bt.9}
{
\langle x_1x_2^3\rangle ={3\over 8}\vert z_2\vert ^2\Re (z_2\overline{z}_1),
}
\begin{eqnarray}\label{Bt.10}
\langle x_1^2x_2^2\rangle &=&{1\over 2^4}\langle
(z_1^2+2z_1\overline{z}_1+\overline{z}_1^2)(z_2^2+2z_2\overline{z}_2+\overline{z}_2^2)\rangle \\
&=& {1\over 2^4}(z_1^2\overline{z}_2^2+4\vert z_1\vert ^2\vert z_2\vert
^2+\overline{z}_1^2z_2^2)\nonumber\\
&=& {1\over 8}\Re (z_1^2\overline{z}_2^2)+{1\over 4}\vert z_1\vert ^2\vert
z_2\vert ^2.\nonumber
\end{eqnarray}
Notice that our averages are invariant under the anti-symplectic
involution
\ekv{Bt.11}
{
j:(x,\xi )\mapsto (x,-\xi ).
}
This is necessarily the case since we stay in the framework of ordinary
Schr{\"o}dinger operators (without magnetic fields) whose symbols have this
invariance.

\par Now write our results in the action angle variables $(\rho _j,\theta
_j)$, given by
\ekv{Bt.12}
{
z_j=\sqrt{2\rho _j}e^{-i\theta _j},
}
so that
$${1\over 2}\vert z_j\vert ^2={1\over 2}(x_j^2+\xi _j^2)=\rho _j:$$
\ekv{Bt.13}
{
\langle x_1^4\rangle ={3\over 2}\rho _1^2,\ \langle x_2^4\rangle ={3\over
2}\rho _2^2,
}
\ekv{Bt.14}
{
\langle x_1^3x_2\rangle ={3\over 2}\rho _1^{3/2}\rho _2^{1/2}\cos (\theta
_1-\theta _2),\ \langle x_1x_2^3\rangle ={3\over 2}\rho _1^{1/2}\rho
_2^{3/2}\cos (\theta _2-\theta _1),
}
\ekv{Bt.15}
{
\langle x_1^2x_2^2\rangle =\rho _1\rho _2+{1\over 2}\rho _1\rho _2\cos
2(\theta _1-\theta _2).
}
It follows from the Hamilton equations that $\rho _j$ and $\theta :=\theta _1-\theta _2$ are constant along
every $H_p$-trajectory. The involution $j$ can also be described as
$(z_1,z_2)\mapsto (\overline{z}_1,\overline{z}_2)$, and hence in action--angle variables as
\ekv{Bt.16}
{(\rho _1,\rho _2,\theta _1,\theta _2)\mapsto (\rho _1,\rho _2,-\theta _1,-\theta _2)}

\par We shall study our averages as \fu{}s on the abstract symplectic
\mfld{}
\ekv{Bt.17}
{
\Sigma =p^{-1}(1)/\exp {\bf R}H_p.
}
Using the $(z_1,z_2)$-coordinates, we have
\ekv{Bt.18}
{
p^{-1}(1):\,\, {1\over 2}(\vert z_1\vert ^2+\vert z_2\vert ^2)=1,
}
and the equivalence relation induced by the $H_p$-flow is: $(z_1,z_2)\sim
(w_1,w_2)$ if and only if
$$
(w_1,w_2)=(e^{it}z_1,e^{it}z_2)
$$
for some $t\in{\bf R}$. Thus we see that $\Sigma $ can be identified with the complex projective
space $P({\bf C}^2)$. It is well known that this space is diffeomorphic
to ${\bf S}^2$. Indeed, $P({\bf C}^2)$ can be identified with the 1-point
compactification ${\bf C}\cup\{\infty \}$ via the map $(z_1,z_2)\mapsto
z_1/z_2$ and the one point compactification can be identified with the
Riemann sphere.

\par $\Sigma $ can be parametrized by $(\rho _1,\rho _2,\theta )$ with
$\rho _j\ge 0$, $\rho _1+\rho _2=1$, $\theta \in {\bf R}/2\pi {\bf Z}$,
with the convention that all the $(1,0,\theta )$ denote the same point and
similarly for $(0,1,\theta )$. The involution $j$ induces the
anti-symplectic involution
\ekv{Bt.19}
{
j:\,\, \Sigma \ni (\rho _1,\rho _2,\theta )\mapsto (\rho _1,\rho
_2,-\theta ).
}
Notice that the set of fixed points of $j$ is given by all points with
$\theta =0$ or $\theta =\pi $. These points form a (great) circle on
$\Sigma $ and can also be described as the set of trajectories in $p^{-1}(1)$
whose $x$-space projections hit the \bdy{} of the potential well
$\{x\in{\bf R}^2;\, \vert x\vert =1\}$.

\par We consider perturbations of the form
\ekv{Bt.20}
{
q(x)={2\over 3}a(x_1^4+x_2^4)+bx_1^2x_2^2+{2\over 3}c(x_1^3x_2+x_1x_2^3).
}
Then on $\Sigma $ we get with $\rho =\rho _1$, so that $\rho _2=1-\rho $:
\begin{eqnarray}\label{Bt.21}
\langle q\rangle &=&a(\rho _1^2+\rho _2^2)+b\rho _1\rho _2+{b\over 2}\rho
_1\rho _2\cos (2\theta )+c(\rho _1^{1\over 2}\rho _2^{3\over 2}+\rho
_1^{3\over 2}\rho _2^{1\over 2})\cos \theta
\\
&=& a(\rho ^2+(1-\rho )^2)+b\rho (1-\rho )(1+{1\over 2}\cos (2\theta ))+c(\rho ^{1\over
2}(1-\rho )^{3\over 2}+\rho ^{3\over
2}(1-\rho )^{1\over 2})\cos \theta
\nonumber
\\
&=& a+(b-2a)\rho (1-\rho )+{b\over 2}\rho (1-\rho )\cos (2\theta )+c\rho
^{1\over 2}(1-\rho )^{1\over 2}\cos \theta
\nonumber
\\
&=& a+({b\over 2}-2a)\rho (1-\rho )+b\rho (1-\rho )\cos ^2\theta +c\rho
^{1\over 2}(1-\rho )^{1\over 2}\cos \theta ,
\nonumber
\end{eqnarray}
where we used that $\rho ^{1\over 2}(1-\rho )^{3\over 2}+\rho ^{3\over
2}(1-\rho )^{1\over 2}=\rho ^{1\over 2}(1-\rho )^{1\over 2}$.

\par We are interested in the critical points of this \fu{} on $\Sigma $,
and the values $\rho =0,1$ will have to be treated separately. In
particular we are interested in the number of saddle points. If we have
only one saddle point we will be able to apply the results of this
paper.
This is still the case if there are two saddle points provided that the
corresponding critical values are different. We will also encounter the
case of two saddle points $S_1,S_2$ away
from the equator and then necessarily with $j(S_1)=S_2$. In that case the
critical vales will be equal and the results of this paper will not apply
directly. We plan to return to that case in a future paper, where the role
of symmetries will be studied.

\par Put
$$d={b\over 2}-2a,\ g=\rho ^{1\over 2}(1-\rho )^{1\over 2},\ y=\cos \theta .$$
Then, \ekv{Bt.22} { \langle q\rangle =a+dg^2+bg^2y^2+c g y. }
Notice that $y=\cos \theta $ is critical precisely when $\theta
=0,\pi $ and that $y\in [-1,1] $. When $y\ne \pm 1$, we may treat
$g$ as an \indep{} variable. The same observation is valid for
$g(\rho )\in ]0,{1\over 2}]$. It is non-critical in $[0,{1\over
2}[$ (i.e. for $\rho \ne {1\over 2}$). (As already mentioned, the
value $g=0$, corresponding to $\rho =0,1$, will require a
different treatment.)

\par In order to avoid various degenerations, we shall assume
\ekv{Bt.22.5}
{d\ne 0,}
\ekv{Bt.22.6}
{\hbox{When }c\ne 0,\hbox{ we have }b\ne 0,\,\, b+d\ne 0.}
\medskip
\par\noindent 1) {\it Critical points with}
\ekv{Bt.23}
{
 \theta \ne 0,\pi,\ \rho \ne 0,{1\over 2},1.
}
Here both $y$ and $g$ can be treated as \indep{} variables and the critical
points are determined by
$$\cases{2bg^2y+cg=0,
\cr
2dg+2bgy^2+cy=0.
}$$
This can also be written
$$\cases{g(2bgy+c)=0,\cr
2dg+y(2bgy+c)=0.
}$$
Under the assumption (\ref{Bt.23}) we have $g\ne 0$, so we get
$$
\cases{
2bgy+c=0,\cr 2dg=0.
}
$$
This is in contradiction with the assumption (\ref{Bt.22.5}), so we
conclude that there are no critical points away from the union of the
"vertical circle" given by $\theta \in\{ 0,\pi \}$ and the horizontal
circle: $\rho ={1\over 2}$.\medskip

\par\noindent 2) {\it Critical points on the horizontal circle away from
the vertical one:}
\ekv{Bt.24}
{
\theta \ne 0,\pi ,\quad \rho ={1\over 2}.
}
Then $g={1/2}$ and this is a critical value, so we only have to look for
critical points \wrt{} $y$, leading to
$$2b({1\over 2})^2y+c{1\over 2}=0,$$
\ekv{Bt.25} { y=-{c\over b}. } Recall that $\vert y\vert <1$ under
the assumption (\ref{Bt.24}), so we reach the conclusion that if
\ekv{Bt.26}
{\abs{c\over b} <1,}
then there are two distinct
critical points in the region (\ref{Bt.24}), given by \ekv{Bt.27}
{\rho ={1\over 2},\ \cos \theta =-{c\over b},} and otherwise there
are no such points. (In the remaining degenerate case $b=c=0$ the
whole horizontal circle is critical.)

\par We also study the nature of the critical points, by computing the
hessian of $\langle q\rangle $ \wrt{} $\rho ,y$. Using that $g'({1\over
2})=0$, $g''({1\over 2})=-2$, we get at both points
\ekv{Bt.26.1}
{
\partial _y^2\langle q\rangle ={b\over 2},\ \partial _y\partial _\rho
\langle q\rangle =0,\
\partial _\rho ^2\langle q\rangle =-2d.}
So both critical points are of signature $(b,-d)$ where the first
component corresponds to the horizontal ($\theta $) direction.
(We use the convention that a signature described by $(\alpha ,\beta )$ is
given by $({\rm sign\,}(\alpha ),{\rm sign\,}(\beta ))$.)
\medskip

\par\noindent 3) {\it Critical points on the vertical circle away from the
horizontal one and from the poles $\rho =0,1$:}
\ekv{Bt.27.1}
{
\theta \in\{ 0,\pi \} ,\ \rho \not\in \{ 0,\, {1\over 2},\, 1\} .
}
For $\theta =0$, we have $y=1$ and we look for critical points of
$g\mapsto (d+b)g^2+cg$, leading to
\ekv{Bt.27.2}
{
g=-{c\over 2(b+d).}
}
Hence we get two critical points in this region if
\ekv{Bt.27.3}
{-1<{c\over b+d}<0.}
and otherwise no point on this half of the vertical circle. (In the
degenerate case $b+d=0$, $c=0$ the whole vertical circle is critical.)

\par For $\theta =\pi $, we have $y=-1$ and we look for critical points
of $g\mapsto (d+b)g^2-cg$, leading to
\ekv{Bt.27.4}
{g={c\over 2(b+d)},}
so we get 2 critical points in this case if
\ekv{Bt.27.5}
{0<{c\over b+d}<1,}
and otherwise no critical points on this half of the vertical circle. We
will see shortly that we have critical points at the poles when $c=0$.

\par In both subcases, we get by a straight forward calculation:
$$
\langle q\rangle ''_{gg}=2(b+d),\ \langle q\rangle ''_{g\theta }=0,\
\langle q\rangle ''_{\theta \theta }={c^2\over 2(b+d)^2}d,
$$
so the signature is
\ekv{Bt.27.6}
{
(d+b,d)
}
where the first component refers to the direction of the vertical circle
through the critical point.\medskip

\par\noindent 4) {\it The two points of intersection of the two circles:}
Here $\rho ={1\over 2}$ and $\theta \in\{ 0,\pi \}$. Here both
$g$ and $y$ are critical, so our intersection points are both critical.
By straight forward calculation, we get for $C_f:$ $\theta =0$, $\rho ={1\over 2}$
\ekv{Bt.27.7}
{
\langle q\rangle ''_{\theta \theta }=-{1\over 2}(b+c),\ \langle q\rangle
''_{\theta \rho }=0,\ {1\over 2}\langle q\rangle ''_{\rho \rho }=-c-b-d,
}
and for $C_b:$
$\theta =\pi $, $\rho ={1\over 2}$
\ekv{Bt.27.75}
{
\langle q\rangle ''_{\theta \theta }={1\over 2}(-b+c),\ \langle q\rangle
''_{\theta \rho }=0,\ {1\over 2}\langle q\rangle ''_{\rho \rho }=c-b-d.
}
In particular, the signature is
\ekv{Bt.27.8}
{\cases{(-c-b-d,-b-c),\hbox{ when }\theta =0,\cr
(c-b-d,c-b),\hbox{ when }\theta =\pi .}}
Also notice that
\ekv{Bt.27.85}
{
\langle q\rangle (C_f)=a+{d+b\over 4}+{c\over 2},\quad \langle q\rangle
(C_b)=a+{d+b\over 4}-{c\over 2}
}
\medskip
\par\noindent 5) {\it It remains to study the "poles", given by
} $\rho =0,1$. Here the $\rho ,\theta $ coordinates degenerate and we
return to the $z$-coordinates. Using that
$$\Re (z_1^2\overline{z}_2^2)=2(\Re z_1\overline{z}_2)^2-\vert z_1\vert
^2\vert z_2\vert ^2,$$
we get
$$\langle q\rangle ={a\over 4}(\vert z_1\vert ^2+\vert z_2\vert
^2)^2+({b\over 8}-{a\over 2})\vert z_1\vert ^2\vert z_2\vert ^2+{b\over
4}(\Re z_1\overline{z}_2)^2+{c\over 4}(\vert z_1\vert ^2+\vert z_2\vert
^2)\Re (z_1\overline{z}_2).
$$
Make the change of variables $\zeta _j=z_j/\sqrt{2}$ and restrict to the
energy surface $p^{-1}(1)$, which now becomes $\vert \zeta _1\vert ^2+ \vert
\zeta _2\vert ^2=1$. Then we get
\ekv{Bt.27.9}
{
\langle q\rangle = a+d\vert \zeta _1\vert ^2\vert \zeta _2\vert ^2+b(\Re
(\zeta _1\overline{\zeta }_2))^2+c\Re (\zeta _1\overline{\zeta }_2),
}
again with $d={b\over 2}-2a$.

\par Recall that we work on the projective space, described as the
3-sphere $\vert \zeta _1\vert ^2+\vert \zeta _2\vert ^2=1$ modulo the
action of the rotations $t\mapsto (e^{it}\zeta _1,e^{it}\zeta _2)$.
Consider the case $\rho =0$. Correspondingly, we can choose the point
$(\zeta _1^0,\zeta _2^0)=(0,1)$. The $H_p$-integral curve through that point
is $t\mapsto (0,e^{-it})$ and locally, we can identify $\Sigma $ with the
transversal hypersurface $H$ in the 3-sphere which is given by $\Im \zeta
_2=0$. Thus $\zeta _2=1-w$ with $w\in{\rm neigh\,}(0,{\bf R})$, and we get
$$w=1-(1-\vert \zeta _1\vert ^2)^{1\over 2}={1\over 2}\vert \zeta _1\vert
^2+{\cal O}(\vert \zeta _1\vert ^4).$$
We can use the real and imaginary parts of $\zeta _1$ as local coordinates
on $H$. Then on $H$, we get the Taylor expansion:
$$\langle q\rangle =a+d\vert \zeta _1\vert ^2+b(\Re \zeta _1)^2+c\Re \zeta
_1+{\cal O}(\vert \zeta _1\vert ^3).$$

\par We conclude that the "pole" $\rho =0$ is a critical point iff $c=0$ and
when this point is critical, the signature  is $(d+b,d)$, where the first
component corresponds to the direction of the (vertical) circle through
the pole. By symmetry in the indices 1, 2, we have the identical conclusion
for the opposite pole, given by $\rho =1$. Notice finally that this case
together with the case 3 give a complete description of the critical
points on the vertical circle away from the crossings with the
horizontal one.\medskip

\par We observe that the critical points away from the intersection of the
two circles are non-degenerate and keep constant signatures under
small perturbations of the parameters (except in the degenerate
cases $c=b=0$ and $c=b+d=0$). These critical points can only be
killed or born by passing through one of the two crossing points.
This happens in the following 4 cases:
\medskip
\par\noindent 1) ${c\over b}=-1$: The critical points on the horizontal
circle coalesce into the crossing point $\theta =0,\, \rho ={1\over 2}$.
When $c/b$ goes from $-1+\epsilon $ to $-1-\epsilon $, the two critical
points disappear and the signature at $\theta =0,\, \rho ={1\over 2}$ goes
from
$(-d,-b\epsilon )$ to $(-d,b\epsilon )$\medskip

\par\noindent 2) ${c\over b}=1$. The two critical points on the horizontal
great circle coalesce into the crossing point $\theta =\pi ,\, \rho
={1\over 2}$. When ${c\over b}$ goes from $1-\epsilon $ to $1+\epsilon $,
the signature of that crossing point goes from
$(-d,-b\epsilon )$ to $(-d, b\epsilon )$.\medskip

\par\noindent 3) ${c\over b+d}=-1$: The two critical points on the
vertical circle coalesce into the crossing point $\theta =0,$ $\rho ={1\over
2}$. When ${c\over b+d}$ goes from $-1+\epsilon $ to $-1-\epsilon $,
the signature of that crossing point goes from $(c\epsilon ,d)$
to $(-c\epsilon ,d)$. \medskip

\par\noindent 4) ${c\over b+d}=1$: The two critical points on the
vertical circle coalesce into the crossing point $\theta =\pi ,\rho ={1\over
2}$. When ${c\over b+d}$ goes from $1-\epsilon $ to $1+\epsilon $,
the signature of that crossing point goes from $(-c\epsilon ,d)$
to $(c\epsilon ,d)$. \medskip

\par In the following, we may assume in order to fix the ideas, that
$d>0$. In the $b,c$-plane, we define the following open sets, separated
from each other by the 4 lines $c=\pm b$, $c=\pm (b+d)$, where all
the critical points will be \nondeg{}:
\smallskip
\par\noindent $A:\ b>0,\,\, -b<c<b$.
\smallskip
\par\noindent $B_+:\ \max (b,-b)<c<b+d$.
\smallskip
\par\noindent $B_-:\ -(b+d)<c<\min (b,-b)$.
\smallskip
\par\noindent $C_+:\ c>\max (b+d,-b)$.
\smallskip
\par\noindent  $C_-:\ c<\min (b, -b-d)$.
\smallskip
\par\noindent  $D:\ b<0,\,\, \max (b,-b-d)<c<\min (-b, b+d)$.
\smallskip
\par\noindent $E_+:\ \max (b+d,-b-d)<c<-b$.
\smallskip
\par\noindent $E_-:\ b<c<\min (-b-d,b+d)$.
\smallskip
\par\noindent $F:\ b<-d,\,\, b+d<c<-b-d$.

\medskip
\begin{figure}
\centering
    \psfrag{1}[l][l]{$c$}
    \psfrag{2}[l][l]{$b$}
    \psfrag{3}[1][1]{$c=b$}
    \psfrag{4}[1][1]{$c=b+d$}
    \psfrag{5}[1][1]{$c=-b$}
    \psfrag{6}[1][1]{$c=-(b+d)$}
    \psfrag{7}[1][1]{$C_+$}
    \psfrag{8}[1][1]{$A$}
    \psfrag{9}[1][1]{$B_+$}
    \psfrag{10}[1][1]{$D$}
    \psfrag{11}[1][1]{$E_+$}
    \psfrag{12}[1][1]{$F$}
    \psfrag{13}[1][1]{$E_-$}
    \psfrag{14}[1][1]{$C_-$}
    \psfrag{15}[1][1]{$B_-$}
\scalebox{0.9} {\includegraphics{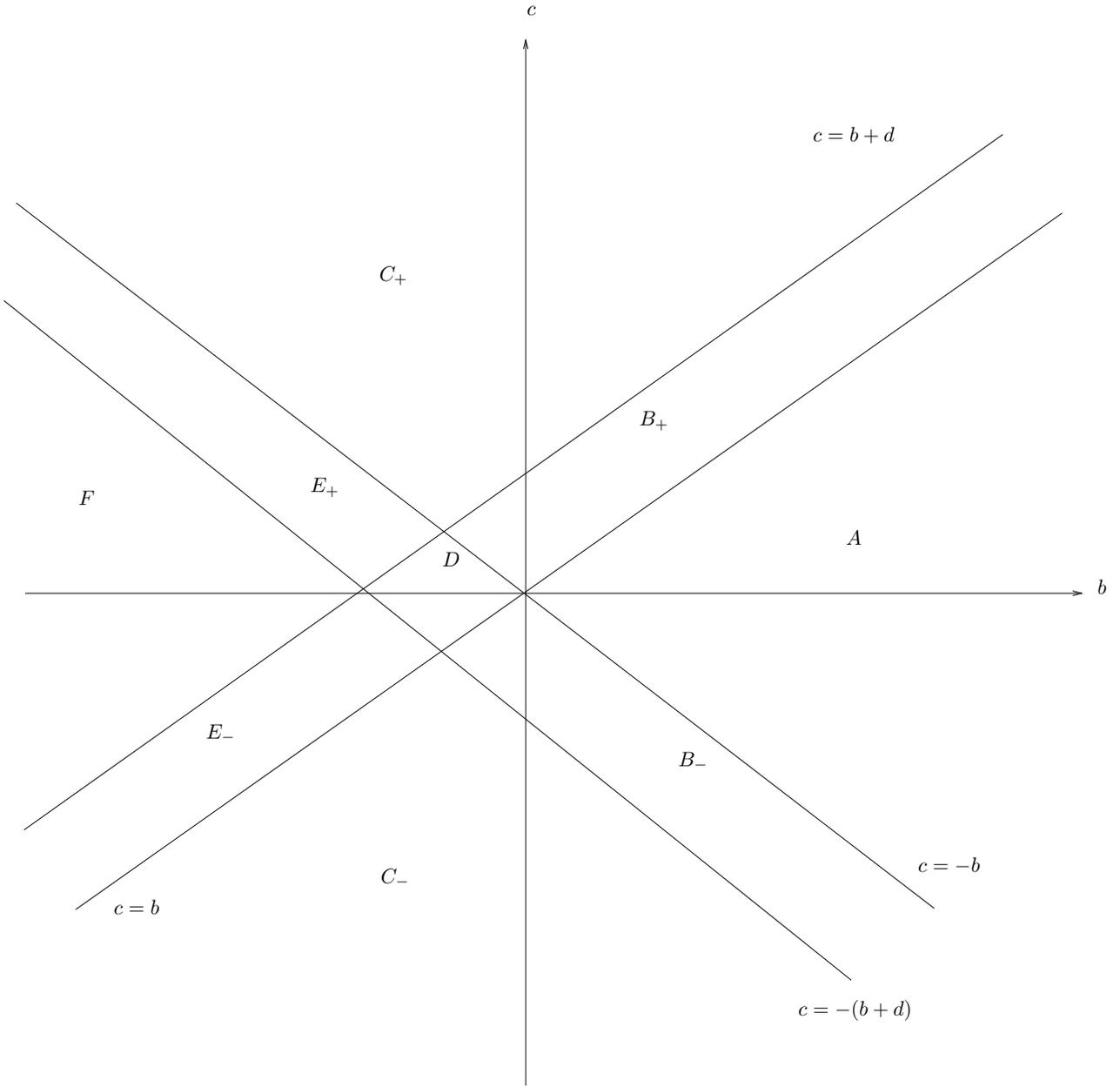}} \caption{When the
parameters are in the regions $B_+$, $B_-$, $E_+$, and $E_-$, we
have precisely one saddle point, and therefore the results of the
present paper apply. In the cases $C_+$ and $C_-$ there are no
saddle points at all, while in the case $D$ there two saddle
points. The corresponding critical values of $\langle{q}\rangle$
are separated provided that we assume that $c\neq 0$.}
\end{figure}

\bigskip
\par Then the earlier discussion gives the location and the signature of
the critical points in each of the cases. Let $C_f$ denote the "forward"
crossing point of the two circles, given by $\rho ={1\over 2}$, $\theta
=0$. Similarly let $C_b$ denote the "backward" crossing point, given by
$\rho ={1\over 2}$, $\theta  =\pi $.
\smallskip

\par\noindent $A$:
\par\noindent Signature at $C_f$: $(-,-)$
\par\noindent Signature at $C_b$: $(-,-)$
\par\noindent Away from the crossings:
\par\noindent On the horizontal circle: Two critical points with
signature $(+,-)$
\par\noindent On the vertical circle: Two critical points with
signature $(+,+)$
\smallskip

\par\noindent $B_+$:
\par\noindent Signature at $C_f$: $(-,-)$
\par\noindent Signature at $C_b$: $(-,+)$
\par\noindent Away from the crossings:
\par\noindent On the horizontal circle: No critical points
\par\noindent On the vertical circle: Two critical points with
signature $(+,+)$
\par\noindent Here $\langle q\rangle (C_b)$ is smaller than $\langle
q\rangle (C_f)$ but larger than the two other critical values.
\smallskip

\par\noindent $B_-$:
\par\noindent Signature at $C_f$: $(-,+)$
\par\noindent Signature at $C_b$: $(-,-)$
\par\noindent Away from the crossings:
\par\noindent On the horizontal circle: No critical points
\par\noindent On the vertical circle: Two critical points with
signature $(+,+)$.
\par\noindent Here $\langle q\rangle (C_f)$ is smaller than $\langle
q\rangle (C_b)$ but larger than the two other critical values.

\smallskip

\par\noindent $C_+$:
\par\noindent Signature at $C_f$: $(-,-)$
\par\noindent Signature at $C_b$: $(+,+)$
\par\noindent Away from the crossings:
\par\noindent On the horizontal circle: No critical points
\par\noindent On the vertical circle: No critical points
\smallskip

\par\noindent $C_-$:
\par\noindent Signature at $C_f$: $(+,+)$
\par\noindent Signature at $C_b$: $(-,-)$
\par\noindent Away from the crossings:
\par\noindent On the horizontal circle: No critical points
\par\noindent On the vertical circle: No critical points
\smallskip

\par\noindent $D$:
\par\noindent Signature at $C_f$: $(-,+)$
\par\noindent Signature at $C_b$: $(-,+)$
\par\noindent Away from the crossings:
\par\noindent On the horizontal circle: Two critical points with
signature $(-,-)$
\par\noindent On the vertical circle:  Two critical points with
signature $(+,+)$
\par\noindent Here $\langle q\rangle (C_f)$, $\langle q\rangle (C_b)$ are
larger than the values at the critical points on the vertical circle and
smaller than the values at the critical points on the horizontal circle.
From (\ref{Bt.27.85}) we also know that $\langle q\rangle (C_f)-\langle
q\rangle (C_b)=c$.
\smallskip

\par\noindent $E_+$:
\par\noindent Signature at $C_f$: $(-,+)$
\par\noindent Signature at $C_b$: $(+,+)$
\par\noindent Away from the crossings:
\par\noindent On the horizontal circle:  Two critical points with
signature $(-,-)$.
\par\noindent On the vertical circle:  No critical points.
\par\noindent In this case $\langle q\rangle (C_f) $ is larger than
$\langle q\rangle (C_b)$ but smaller than the two other critical values.
\smallskip

\par\noindent $E_-$:
\par\noindent Signature at $C_f$: $(+,+)$
\par\noindent Signature at $C_b$: $(-,+)$
\par\noindent Away from the crossings:
\par\noindent On the horizontal circle:  Two critical points with
signature $(-,-)$.
\par\noindent On the vertical circle:  No critical points.
\par\noindent
In this case $\langle q\rangle (C_b) $ is larger than
$\langle q\rangle (C_f)$ but smaller than the two other critical values.
\smallskip

\par\noindent $F$:
\par\noindent Signature at $C_f$: $(+,+)$
\par\noindent Signature at $C_b$: $(+,+)$
\par\noindent Away from the crossings:
\par\noindent On the horizontal circle: Two critical points with
signature $(-,-)$.
\par\noindent On the vertical circle: Two critical points with
signature $(-,+)$.
\smallskip
\par In the cases $B_+,B_-,E_+,E_-$ we have precisely one saddle point
(necessarily) situated on the vertical circle which is the fixed point
set of $j$. In these cases, the results of this paper
apply. The results also apply in the case $D$, provided that we assume
that $c\ne 0$ in order to separate the two saddle point values.
In these cases it is easy to understand the structure and the
shape of the level sets $\langle q\rangle =C$. In particular we
see that when we let $C$ be a saddle point value we get a
connected "$\infty $" shaped set (and no "circular" components on which
$\langle q\rangle $ is non-critical everywhere).

\par In the case $F$, we have two saddle points situated on the
vertical circle symmetrically \wrt{} the horizontal circle.  Since we have
chosen to use perturbations which are symmetric under permutation of
$x_1,x_2$, the function $\langle q\rangle $ is invariant under the map
$\rho \mapsto 1-\rho $, so the critical values are necessarily equal.  Here
we can break the symmetry by adding a small multiple of for instance
$x_1^4$ so that we still have precisely two saddle points, but with
different critical values.  Then the results of our paper
apply.

\par In the case $A$, we have two saddle points on the horizontal circle.
They are of course exchanged by application of $j$ and this symmetry
remains under perturbations within the class of \sop{}s without
magnetic field. We hope to analyze this case in a future work.

We shall next compute the $\epsilon ^2$-contribution to the averaging of
the principal symbol $p(x,\xi )-i\epsilon p_4(x)-\epsilon ^2p_6(x)+{\cal
O}(\epsilon ^3)$ appearing in in (\ref{Bt.5.5}), by applying the
calculations of the end of Section \ref{SectionSk}, with
\ekv{Bt.28}
{
q=-p_4(x),\ r=-p_6(x),\ T=2\pi .
}
Recall from there, that we have the averaged symbol
\ekv{Bt.29}
{
{{p_\epsilon }_\big\vert}_{\Lambda _{\epsilon G}}\simeq p+i\epsilon \langle
q\rangle +\epsilon ^2\left(\langle r\rangle -{1\over 2} C(q,q)\right)+{\cal O}(\epsilon ^3),
}
where $C(q_1,q_2)$ and ${\rm Cor\,}(q_1,q_2)$ were defined in
(\ref{Sk.58}), (\ref{Sk.56}).

\par A simple calculation gives
$$
\{ z^\alpha ,z^\beta \},\ \{\overline{z}^\alpha ,\overline{z}^\beta \}
=0,\ \{z^\alpha ,\overline{z}^\beta\}  =2i\left({\alpha _1\beta _1\over \vert
z_1\vert ^2}+{\alpha _2\beta _2\over \vert
z_2\vert ^2}\right)z^\alpha \overline{z}^\beta .
$$
More generally,
$$
\{z^\alpha \overline{z}^{\widetilde{\alpha }},
z^\beta  \overline{z}^{\widetilde{\beta }}\} =2i\left( {\sigma
(\widetilde{\alpha }_1,\alpha _1;\widetilde{\beta }_1,\beta _1)\over \vert
z_1\vert ^2}+ {\sigma
(\widetilde{\alpha }_2,\alpha _2;\widetilde{\beta }_2,\beta _2)\over \vert
z_2\vert ^2}\right) z^{\alpha +\beta }\overline{z}^{\widetilde{\alpha
}+\widetilde{\beta }},
$$
where $\sigma $ denotes the symplectic form, viewed as an alternate
bilinear form on $T^*{\bf R}^2\times T^*{\bf R}^2$.
Hence
$$
{\rm Cor\,}(z^\alpha \overline{z}^{\widetilde{\alpha }},
z^\beta  \overline{z}^{\widetilde{\beta }};s) =
2i\left( {\sigma
(\widetilde{\alpha }_1,\alpha _1;\widetilde{\beta }_1,\beta _1)\over \vert
z_1\vert ^2}+ {\sigma
(\widetilde{\alpha }_2,\alpha _2;\widetilde{\beta }_2,\beta _2)\over \vert
z_2\vert ^2}\right) z^{\alpha +\beta }\overline{z}^{\widetilde{\alpha
}+\widetilde{\beta }}e^{is(\vert \widetilde{\alpha }\vert -\vert \alpha
\vert )},
$$
when $\vert \widetilde{\alpha }\vert -\vert \alpha \vert
=\vert \beta \vert -\vert \widetilde{\beta }\vert $, and ${\rm Cor\,}(z^\alpha \overline{z}^{\widetilde{\alpha }},
z^\beta  \overline{z}^{\widetilde{\beta }};s)=0$ otherwise.
If $a=(a_1,a_2),\, b=(b_1,b_2)\in{\bf N}^2$ with $\vert a\vert =\vert
b\vert =4$ we get by multinomial expansion
\begin{eqnarray}\label{Bt.29.5}
& & {\rm Cor\,}(x^a,x^b) = {1\over 2^8}{\rm
Cor\,}((z+\overline{z})^a,(z+\overline{z})^b;s)\\
& &={1\over 2^8}\sum_{\alpha
+\widetilde{\alpha }=a\atop \beta +\widetilde{\beta }=b}\pmatrix{a\cr
\alpha }\pmatrix{b\cr \beta }{\rm Cor\,}(z^\alpha
\overline{z}^{\widetilde{\alpha }},z^\beta \overline{z}^{\widetilde{\beta
}};s)\nonumber\\
& & ={2i\over 2^8}\sum_{{\alpha
+\widetilde{\alpha }=a\atop \beta +\widetilde{\beta }=b}\atop \vert
\widetilde{\alpha }\vert -\vert \alpha \vert =\vert \beta \vert -\vert
\widetilde{\beta }\vert }\pmatrix{a\cr
\alpha }\pmatrix{b\cr \beta }\big( {\sigma
(\widetilde{\alpha }_1,\alpha _1;\widetilde{\beta }_1,\beta _1)\over \vert
z_1\vert ^2}+ {\sigma
(\widetilde{\alpha }_2,\alpha _2;\widetilde{\beta }_2,\beta _2)\over \vert
z_2\vert ^2}\big) z^{\alpha +\beta }\overline{z}^{\widetilde{\alpha
}+\widetilde{\beta }}e^{is(\vert \widetilde{\alpha }\vert -\vert \alpha
\vert )}.
\nonumber
\end{eqnarray}
When calculating this kind of expressions, it is useful to observe that
the relations $\vert \alpha \vert +\vert \widetilde{\alpha }\vert =\vert
\beta \vert +\vert \widetilde{\beta }\vert =4$, $\vert \widetilde{\alpha
}\vert -\vert \alpha \vert =\vert \beta \vert -\vert \widetilde{\beta
}\vert $ imply: $\vert \widetilde{\beta }\vert =\vert \alpha \vert $,
$\vert \beta \vert =\vert \widetilde{\alpha }\vert $.

\par We have the \F{} series expansion
$$
1_{[0,2\pi [}(s)(s-\pi )=\sum_{k\in{\bf Z}\setminus \{ 0\} }{i\over
k}e^{isk}.
$$
Combining this with (\ref{Bt.29.5}) and the Parseval identity, we
get
\begin{eqnarray}\label{Bt.30}
C(x^a,x^b)&=&{1\over 2\pi }\int_0^{2\pi }(s-\pi ){\rm Cor\,}(x^a,x^b;s)\\
&=& {2\over 2^8}\sum_{{{\alpha
+\widetilde{\alpha }=a\atop \beta +\widetilde{\beta }=b}\atop \vert
\widetilde{\alpha }\vert -\vert \alpha \vert =\vert \beta \vert -\vert
\widetilde{\beta }\vert}\atop \vert \widetilde{\alpha }\vert -\vert
\alpha \vert \ne 0 }{\pmatrix{a\cr
\alpha }\pmatrix{b\cr \beta }\over \vert \widetilde{\alpha }\vert -\vert
\alpha \vert }\left( {\sigma
(\widetilde{\alpha }_1,\alpha _1;\widetilde{\beta }_1,\beta _1)\over \vert
z_1\vert ^2}+ {\sigma
(\widetilde{\alpha }_2,\alpha _2;\widetilde{\beta }_2,\beta _2)\over \vert
z_2\vert ^2}\right) z^{\alpha +\beta }\overline{z}^{\widetilde{\alpha
}+\widetilde{\beta }}.
\nonumber
\end{eqnarray}
Using this formula we get after a few days of simple but tedious
calculations:
\ekv{Bt.31}
{
C(x_1^4+x_2^4,x_1^4+x_2^4)=-{17\over 16}(\vert z_1\vert ^6+\vert z_2\vert
^6),
}
\ekv{Bt.32}
{
C(x_1^4+x_2^4,x_1^2x_2^2)=-{3\over 2^6}(3\vert z\vert
^2(z_1^2\overline{z}_2^2+\overline{z}_1^2z_2^2)+16\vert z_1\vert ^2\vert
z_2\vert ^2),
}
\ekv{Bt.33}
{
C(x_1^4+x_2^4,x_1^3x_2+x_1x_2^3)=
{1\over 2^7}(2(z_1^3\overline{z}_2^3+\overline{z}_1^3z_2^3)-(51(\vert
z_1\vert ^4+\vert z_2\vert ^4)+36\vert z_1\vert ^2\vert z_2\vert
^2)(z_1\overline{z}_2+\overline{z}_1z_2)),
}
\ekv{Bt.34}
{
C(x_1^2x_2^2,x_1^2x_2^2)=-{1\over 2^6}\vert z\vert ^2(9\vert z_1\vert
^2\vert z_2\vert ^2+8(z_1^2\overline{z}_2^2+\overline{z}_1^2z_2^2)),
}
\ekv{Bt.35}
{
C(x_1^2x_2^2,x_1^3x_2+x_1x_2^3)=-{1\over 2^8}((17(\vert z_1\vert ^4+\vert
z_2\vert ^4)+90\vert z_1\vert ^2\vert z_2\vert
^2)(\overline{z}_1z_2+z_1\overline{z}_2)+12(\overline{z}_1^3z_2^3+z_1^3
\overline{z}_2^3)),
}
\ekv{Bt.36}
{
C(x_1^3x_2+x_1x_2^3,x_1^3x_2+x_1x_2^3)=
-{1\over 2^8}(17(\vert z_1\vert ^6+\vert z_2\vert ^6)+153\vert z_1\vert
^2\vert z_2\vert ^2\vert z\vert ^2+51\vert z\vert
^2(z_1^2\overline{z}_2^2+\overline{z}_1^2z_2^2)).
}

\par Now recall that $q$ is given by (\ref{Bt.20}), so that by
(\ref{Bt.29}), we have
\begin{eqnarray}\label{Bt.37}
{{p_\epsilon }_\vert}_{\Lambda _{\epsilon G}}&\simeq& p+i\epsilon (\langle
q\rangle +i\epsilon f(r,a,b,c)+{\cal O}(\epsilon ^2))=:p+i\epsilon
\widetilde{q}_\epsilon ,\\
f(r,a,b,c)&=&-\langle r\rangle +{1\over 2} \big({4\over 9}a^2C(x_1^4+x_2^4,x_1^4+x_2^4)
+b^2C(x_1^2x_2^2,x_1^2x_2^2)+\nonumber \\
&&{4\over 9}c^2C(x_1^3x_2+x_1x_2^3,x_1^3x_2+x_1x_2^3)+{4ab\over
3}C(x_1^4+x_2^4,x_1^2x_2^2)+\nonumber\\
&&{8ac\over 9}C(x_1^4+x_2^4,x_1^3x_2+x_1x_2^3)+{4bc\over
3}C(x_1^2x_2^2,x_1^3x_2+x_1x_2^3)\big)+{\cal O}(\epsilon ^2)\nonumber
\end{eqnarray}
According to (\ref{Para.5}) our reduced 1-dimensional \op{} has the symbol
$$
Q_\epsilon =\widetilde{q}_\epsilon +{\cal O}(h+{h^{N_0}\over \epsilon }).
$$
Put $$q_s=sQ_\epsilon +(1-s)\langle q\rangle .$$
If we assume that $\epsilon \gg h$, then we get according to (\ref{Sk.50}):
\ekv{Bt.38}
{
\int_{\gamma _1(Q_\epsilon )}\xi dx-\int_{\gamma _1(\langle q\rangle )}\xi
dx=-i\epsilon \int_{\gamma _1(\langle q\rangle )}[f(r,a,b,c)]_{\rho
_c}^{(x(t),\xi (t))}dt+{\cal O}(\epsilon ^2+h),
}
where we recall that $\int_{\gamma _1(\langle q\rangle )}\xi dx$ is the
(real) action along a loop in $\langle q\rangle ={\rm Const}=\langle
q\rangle (\rho _c)$ starting and ending at the saddle point $\rho _c$, and
that $\int_{\gamma (Q_\epsilon )}\xi dx$ is the corresponding perturbed action
for $Q_\epsilon $. From (\ref{Bt.31}), we see that
$C(x_1^4+x_2^4,x_1^4+x_2^4)$ is minimal precisely on the horizontal
circle $\rho =1/2$. In the cases $B_{\pm}$, $E_{\pm}$, $D$
the saddle points belong to $\{ C_f,C_b\}$ situated on that circle. Since
$\langle q\rangle $ is invariant under reflection in that circle,
either the loop $\gamma _1(\langle q\rangle )$ is entirely in the upper or
lower hemisphere intersecting the equator only at $\rho _c$ (and this is
happens in the cases $B_{\pm}$ and for one of the saddles in case $D$)
or $\gamma _1(\langle q\rangle )$ intersects the equator at one more point
and is symmetric around the equator (and this happens in the cases
$E_{\pm}$ and for one of the saddles in case $D$). In both cases we see
that
\ekv{Bt.39}
{
\int_{\gamma _1(\langle q\rangle )}[C(x_1^4+x_2^4,x_1^4+x_2^4)]_{\rho
_c}^{(x(t),\xi (t))}dt>0.
}
Taking into account the form of $f$ in (\ref{Bt.37}) we conclude that for
every $r$ the integral in the \lhs{} of (\ref{Bt.38}) is $\ne 0$ except
for $(a,b,c)$ in a set of measure $0$. For $(a,b,c)$ outside that
exceptional set, we conclude from the discussion at the end of Section
\ref{SectionSk} that the spectrum of the one dimensional localized
\op{}s has a genuinely two-dimensional structure.

\appendix
\section{Proof of Proposition \ref{Prop1D2}}\label{Appendix2}
\setcounter{equation}{0}

To get a complete normal form we shall do further conjugations
with analytic \pop{}s of order 0 in such a way that the complete
symbol also becomes a \fu{} of $\tau$, $\epsilon$, $h^2/\epsilon $
and $x\xi $. Moreover, we need to do so with errors that are
${\cal O}(e^{-1/(Ch)})$ (rather than merely ${\cal O}(h^\infty )$
as in \cite{HiSj1}.  $Q$ is not a \clas{} but it has a \hol{}
realization and becomes a \clas{}, if we allow some of the
$h$-dependence to appear as an \indep{} parameter in the \coef{}s
of the $h$-\asy{} expansion.  Thus, our starting point will be a
symbol of the form \ekv{1D.13} { Q=Q_0(\tau ,x\xi ,\epsilon
,h^2/\epsilon )+hQ_1(\tau ,x,\xi ,\epsilon ,h^2/\epsilon ;h), }
where $Q_1$ is \hol{} and ${\cal O}(1)$ in some fixed complex
\neigh{} of $\tau =0$, $x=\xi =0$.

\par We define the $\rho $-quasi-norm as above, but now it is important
that we work in the Weyl quantization. To an \an{} symbol $a$ we
then associate the infinite order \dop{} $A={\rm Op}_a(x,\xi
,D_{x,\xi };h)$ as in (\ref{Exp.05}). From the definition, we
verify the following metaplectic invariance property: If $\kappa
:{\bf C}^{2n}\to {\bf C}^{2n}$ is an affine linear \ctf{} and
$\kappa ^*$, $\kappa _*$ denote the usual operation of pull-back
and push-forward of \fu{}s on ${\bf C}^{2n}$, then \ekv{1D.14} {
\kappa _*{\rm Op}_a\kappa ^*={\rm Op}_{\kappa _*a}. } This implies
that if we define our quasi-norms with the help of a \fy{} of
opens sets $\Omega _t$ which are invariant under $\kappa $, then
\ekv{1D.15} { \3 \kappa _*a\3_\rho =\3 a\3_\rho . }

\par In the case of (\ref{1D.13}), we shall let $\Omega _t$ be of the
form $\vert x\vert ^2+\vert \xi \vert ^2\le r(t)$ for a suitable $r(t)$,
and we observe that these balls are invariant under $\exp isH_{x\xi }$
when $s$ is real. After applying the inverse \fu{} of $Q_0(\tau ,\cdot
,\epsilon ,h^2/\epsilon )$ to our \op{}, we may assume that the principal
symbol of $Q$ is $x\xi $, so (\ref{1D.13}) simplifies to
\ekv{1D.16}
{
Q=x\xi +hQ_1(\tau ,x,\xi ,\epsilon ,h^2/\epsilon ;h)
}
with a new $Q_1$ having the same properties as the previous one.

\par Using the same letters for \op{}s and their symbols, we let
$Q_0={1\over 2}(xhD+hDx)$ be the quantization of $x\xi $. Notice
that $\exp (2\pi Q_0/h)=-1$, so that $\exp (2\pi\, {\rm
ad}_{Q_0}/h)=1$. If $B$ is an analytic $h$-\pop{} of order 0, we
put \ekv{1D.17} { \langle B\rangle ={1\over 2\pi }\int_0^{2\pi
}e^{tQ_0/h}Be^{-tQ_0/h}dt, } and notice that on the symbol level,
\ekv{1D.18} { \langle B\rangle ={1\over 2\pi }\int_0^{2\pi }B\circ
\exp itH_{x\xi }dt. } Also notice that $[Q_0,\langle B\rangle
]=0$. Choosing the $\rho $-quasi-norms as above, we further have
that \ekv{1D.19} { \3 \langle B\rangle \3_\rho \le \3 B\3_\rho . }

\par The \e{},
\ekv{1D.20}
{
{\rm ad}_{Q_0}A=B-\langle B\rangle
}
has the solution
\ekv{1D.21}
{
A={1\over h}\int k(t) e^{tQ_0/h}Be^{-tQ_0/h}dt={1\over h}\int k(t)
e^{tQ_0/h}(B-\langle B\rangle )e^{-tQ_0/h}dt,
}
where $k(t)$ is the \fu{} with support in $[-\pi ,\pi ]$ which is affine
on $[-\pi ,0[$, $]0,\pi ]$ with $k(\pm \pi )=0$, $k(\pm 0)=\mp {1\over 2}$.
We have
\ekv{1D.22}
{
\3 A\3_\rho \le C\3 {1\over h}(B-\langle B\rangle )\3_\rho  .
}

\par As in Section \ref{SectionExp}, we see that the map
\ekv{1D.23}
{
A\mapsto {\rm Ad}_A(Q)
}
has the differential
\ekv{1D.24}
{
\delta A\mapsto {\rm ad}_{\delta A}Q_0+\widetilde{K}(A,\delta A),
}
where
\ekv{1D.25}
{
\3 \widetilde{K}(A,\delta A)\3_\rho \le C\rho (h+\3 A\3_\rho )\3 \delta
A\3_\rho ,
}
under the assumption that $\3 A\3_\rho ={\cal O}(1)$.

\par Now return to (\ref{1D.16}). After a first conjugation, we may reduce
ourselves to the case when $Q_1-\langle Q_1\rangle $ is ${\cal
O}(h)$, so that \ekv{1D.25.5} { \3 Q_1-\langle Q_1\rangle \3_\rho
\le Ch, } for some $C>0$, when $\rho \le \rho _0>0$. We look for
$A$ such that ${\rm Ad}_AQ$ commutes with $Q_0$, and we try
$A=\sum_0^\infty  A_j$ with convergence in some $\rho
$-quasi-norm. Start by solving \ekv{1D.26} {
[A_0,Q_0]+hQ_1=h\langle Q_1\rangle ,\hbox{ with }\3 A_0\3_\rho \le
C\3 Q_1-\langle Q_1\rangle \3_\rho \le {\cal O}(h). } From
(\ref{1D.24}), (\ref{1D.25}), we get \ekv{1D.27} { {\rm
Ad}_{A_0}Q_0=Q_0+h\langle Q_1\rangle +hQ_2, }
\begin{eqnarray}\label{1D.28}
\3 Q_2\3_\rho &\le& h^{-1}C\rho (h+C\3 Q_1-\langle Q_1\rangle \3_\rho )\3
Q_1-\langle Q_1\rangle \3_\rho  \\
&\le &C\rho (1+C^2)\3 Q_1-\langle Q_1\rangle \3_\rho \nonumber\\
&\le & {1\over 4}\3 Q_1-\langle Q_1\rangle \3_\rho ,\nonumber
\end{eqnarray}
where the last estimate holds for $0\le \rho \le \rho _0$, with $\rho _0>0$
small enough. Choose $A_1$ with
\ekv{1D.29}
{
[A_1,Q_0]+hQ_2=h\langle Q_2\rangle ,\ \3 A_1\3_\rho \le C\3
Q_2-\langle Q_2\rangle \3_\rho \le {C\over 2}\3 Q_1-\langle Q_1\rangle
\3_\rho .
}
Then
$$
{\rm Ad}_{A_0+A_1}Q=Q_0+h\langle Q_1\rangle +h\langle Q_2\rangle +hQ_3,
$$
with $\3 Q_3\3_\rho \le 2^{-2}\3 Q_1-\langle Q_1\rangle \3_\rho $.
Iterating the procedure, we get $A_j$ with
$$
\3 A_j\3_\rho \le C 2^{-j}\3 Q_1-\langle Q_1\rangle \3 _\rho ,
$$
such that if $A=\sum_0^\infty  A_j$, then
\ekv{1D.30}
{
{\rm Ad}_AQ=Q_0+h\langle Q_1\rangle +h\langle Q_2\rangle +...,\ \3 \langle
Q_j\rangle \3_\rho \le 2^{-j}\3 Q_1-\langle Q_1\rangle \3_\rho .}

\par The previous discussion shows how to find
$U$ so that modulo an error ${\cal O}(e^{-1/(Ch)})$, $U^{-1}QU$ commutes
with $xhD_x$. Moreover $U^{-1}QU$ is a classical analytic \pop{} (after
allowing $h$ as an \indep{} parameter in the \coef{}s in
the \asy{} expansions). Put $x=e^s$ and work near $x=r$ for some fixed
small $r>0$. Then $xhD_x=hD_s$ and since the class of \an{} \pop{}s is
conserved under \an{} changes of variables, we know that
$$
U^{-1}QU=K_{..}(\tau ,s,hD_s;h),
$$
where $K$ is an \an{} symbol. But $[K,hD_s]=0$, so
$$
K=K_{..}(\tau ,hD_s;h)
$$
and returning to the $x$-coordinates, we get the representation
(\ref{1D.31}).\hfill$\Box$

\section{Study of $\Gamma _{j,k}$}\label{Appendix3}
\setcounter{equation}{0}

\par For simplicity, we restrict the attention to the right half-plane,
$\Re \mu \ge 0$ and pick one of the equations in (\ref{Sk.4}), that we
write
\ekv{A3.1}
{
(\Im \mu )\ln {1\over \vert \mu \vert }=F(\mu ),
}
where $F(\mu )$ is \ufly{} Lipschitz in a \neigh{} of $0$. As we have
already observed,
\ekv{A3.2}
{\partial _{\Im \mu }(\Im \mu  \ln {1\over \vert \mu \vert })=\ln {1\over
\vert \mu \vert }-\left({\Im \mu \over \vert \mu \vert }\right)^2\gg 1,}
so (\ref{A3.1}) determines a curve of the form
\ekv{A3.3}
{
\Im \mu =f(\Re \mu ),\hbox{ where }f'(\Re \mu )={\cal O}\left({1\over \ln
1/\vert (\Re \mu ,f(\Re \mu ))\vert }\right)\ll 1.
}
We want to express $f$ in terms of $F(\Re \mu )$ up to small errors.

\par Let us first compare the solution $\mu $ of (\ref{A3.1}) with the
solution $\widetilde{\mu }$ of the simplified \e{}
\ekv{A3.4}
{
\Im \widetilde{\mu }\ln {1\over \vert \widetilde{\mu }\vert }=F(\Re \mu
),\hbox{ with }\Re \widetilde{\mu }=\Re \mu .
}
Using that $F(\mu )-F(\Re \mu )={\cal O}(\Im \mu )$ together with
(\ref{A3.2}), we see that
\ekv{A3.5}
{
\Im \mu -\Im \widetilde{\mu }={\cal O}({\Im \mu \over \ln {1\over \vert
\mu \vert }}),\hbox{ so }\Im \mu \sim \Im \widetilde{\mu },\ \ln {1\over
\vert \mu \vert }\sim \ln {1\over \abs{\widetilde{\mu}}.}}
With this
estimate in mind, we now concentrate on the simplified \e{} (\ref{A3.4}),
and we drop the tildes for simplicity.

\par Assume first that we are in the region
\ekv{A3.6}
{
\vert \Im \mu \vert \le {\cal O}(\Re \mu ).
}
Then
$$
\ln {1\over \vert \mu \vert }=\ln \left({1\over x}\right)\left(1+{\cal O}\left(({y\over
x})^2{1\over \ln 1/x}\right)\right),
$$
where we write $\mu =x+iy$. Thus, if $\mu =\widetilde{\mu }$ solves
(\ref{A3.4}) and (\ref{A3.6}) holds, then we first see that
$$
y\sim {F(x)\over \ln {1\over x}},
$$
and then that \ekv{A3.7} { y={F(x)\over \ln {1\over x}}\left(1+{\cal
O}(1)({F(x)\over x\ln {1\over x}})^2{1\over \ln {1\over x}})\right). }
So, if we assume \ekv{A3.8} { \vert F(x)\vert \le {\cal O}(1)x\ln
{1\over x}, } then we are in the region (\ref{A3.6}), and the
solution $\mu =\widetilde{\mu }=x+iy$ of (\ref{A3.4}) takes the
form (\ref{A3.7}). Combining with (\ref{A3.5}) we get under the
assumption (\ref{A3.8}), \ekv{A3.9} { f(x)=\left(1+{{\cal O}(1)\over
\ln {1\over x}}\right){F(x)\over \ln {1\over x}}. }

\par We next consider the region $x\ll \vert y\vert \ll 1$ and assume for
simplicity that we have $y>0$. Then,
\ekv{A3.10}
{
\ln {1\over \vert \mu \vert }=\ln ({1\over y})\left(1+{\cal O}(1)({x\over
y})^2{1\over \ln {1\over y}}\right),
}
and (\ref{A3.4}) takes the form
\ekv{A3.11}
{
y(\ln {1\over y})\left(1+{\cal O}(1)({x\over y})^2{1\over \ln {1\over y}}\right)=F(x).
}

\par Consider first the simplified \pb{}
\ekv{A3.12} {y\ln {1\over y}=z.} With $Y=\ln (1/y),\, Z=\ln (1/z)$
(both $\gg 1$) we get \ekv{A3.13} { Y-\ln Y=Z } Try the
approximate solution $Y_0=Z+\ln Z$. Then by a simple calculation,
$$
Y_0-\ln Y_0=Z+{\cal O}({\ln Z\over Z}).
$$
Since the derivative of the \lhs{} in (\ref{A3.13}) is close to 1, we see
that the solution $Y$ of that \e{} is of the form $Y=Y_0+{\cal O}((\ln Z)/Z)$;
\ekv{A3.14}
{
Y=Z+\left(1+{\cal O}({1\over Z})\right)\ln Z.
}
Hence the solution of (\ref{A3.12}) is of the form
\ekv{A3.15}
{
y=(1+{\cal O}\left({\ln\ln {1\over z}\over \ln {1\over z}})\right){z\over \ln {1\over
z}}.}
If we replace $z$ by $F(x)$, we get the order of magnitude of the solution
to (\ref{A3.11}):
\ekv{A3.16}
{
y\sim {F(x)\over \ln {1\over F(x)}},
}
and the assumption that $x\ll y$ reads:
\ekv{A3.17}
{
{F(x)\over \ln {1\over F(x)}}\gg x.
}
The earlier arguments show that this condition is equivalent to
\ekv{A3.18}
{
F(x)\gg x\ln {1\over x},
}
which indeed is complementary to (\ref{A3.8}). Using (\ref{A3.16}), we get
$$
({x\over y})^2{1\over \ln {1\over y}}\le {\cal O}(1)({x\ln {1\over
F(x)}\over F(x)})^2{1\over \ln {\ln {1\over F}\over F}}\le {\cal O}(1)
({x\ln {1\over F(x)}\over F(x)})^2{1\over \ln {1\over F(x)}},
$$
where we notice that
$$
{x\ln {1\over F}\over F}\ll 1,
$$
by (\ref{A3.17}). Hence (\ref{A3.11}) gives
\ekv{A3.19}
{
y\ln {1\over y}=(1+{\cal O}(1)({x\ln 1/F(x)\over F(x)})^2{1\over \ln
1/F(x)})F(x),
}and applying (\ref{A3.15}) with $z$ equal to the \rhs{} of (\ref{A3.19}),
we get
$$
y=(1+{\cal O}(1){\ln\ln 1/F\over \ln 1/F}){(1+{\cal O}(1)({x\ln 1/F\over
F})^2{1\over \ln 1/F})F(x)\over (\ln {1\over F}+{\cal O}(1)({x\ln 1/F\over
F})^2{1\over \ln 1/F})},
$$
which simplifies to
\ekv{A3.20}
{
y=(1+{\cal O}(1){\ln\ln 1/F\over \ln 1/F}){F(x)\over \ln 1/F(x)}.
}
Recall that here $\widetilde{\mu }=x+iy$ in the simplified \e{}. To get the
corresponding result for (\ref{A3.1}), we apply (\ref{A3.1}) and conclude
that (\ref{A3.20}) holds for $\mu =x+iy$ solving (\ref{A3.1}),
under the equivalent conditions (\ref{A3.17}), (\ref{A3.18}), and assuming
also $0\le F(x)\ll 1$, $0\le x\ll 1$.

\par Summing up, we have proved:
\begin{prop} \label{PropA3.1} Let $F$ be a \ufly{} Lipschitz \fu{} with $\vert
F\vert \ll 1$, defined in a \neigh{} of $0\in{\bf C}$. Let $\mu =x+if(x)$
be the solution of (\ref{A3.1}). Then for small $x$, we have
\ekv{A3.21}
{
\vert f'(x)\vert \le {\cal O}(1)/\ln (1/\vert x+if(x)\vert ).
}
Further,
\ekv{A3.22}
{
f(x)=(1+{{\cal O}(1)\over \ln (1/\vert x\vert )}){F(x)\over \ln (1/\vert
x\vert )},\hbox{ when }\vert F(x)\vert \le {\cal O}(1)\vert x\vert \ln
(1/\vert x\vert ),}

\ekv{A3.23}
{
f(x)=(1+{\cal O}(1){\ln\ln (1/\vert F(x)\vert ) \over
\ln (1/\vert F(x)\vert }){F(x)\over \ln (1/\vert
F(x) \vert )},\hbox{ when }\vert F(x)\vert \gg \vert x\vert \ln
(1/\vert x\vert ),
}
\end{prop}

\end{document}